%% file: main.tex
\pdfoutput=1
\documentclass[12pt,reqno]{amsart}
\usepackage{amsmath, fullpage, amssymb, amsthm, color, comment, mathtools, tikz, dsfont, tikz-cd}

\usepackage[maxbibnames=99,doi=false,isbn=false,url=true,style=alphabetic]{biblatex}

\usepackage{subfig}

\usepackage{hyperref}

\allowdisplaybreaks[4]

\usepackage{graphicx}

\usepackage{float}

\theoremstyle{plain}
\newtheorem{theorem}{Theorem}[section]
\newtheorem{proposition}[theorem]{Proposition}
\newtheorem{lemma}[theorem]{Lemma}
\newtheorem{corollary}[theorem]{Corollary}

\theoremstyle{definition}

\theoremstyle{remark}
\newtheorem{remark}[theorem]{Remark}

\newcommand{\selfnote}[1]{\noindent {\color{red} ********** \\ #1 \\ **********}}

\newcommand{\vol}[0]{\textnormal{vol}}
\newcommand{\gmodgamma}[0]{\Gamma \backslash G}

\newcommand{\para}[0]{\textnormal{para} }

\newcommand{\card}[0]{\textnormal{card} }

\newcommand{\ceil}[0]{\textnormal{ceil}}
\usepackage{subfiles}

\title{Quantum ergodicity on the Bruhat-Tits building for $\textnormal{PGL}(3, F)$ in the Benjamini-Schramm limit}

\author{Carsten Peterson}

\begin{document}
\sloppy

\begin{abstract}
We study joint eigenfunctions of the spherical Hecke algebra acting on $L^2(\Gamma_n \backslash G / K)$ where $G = \textnormal{PGL}(3, F)$ with $F$ a non-archimedean local field of arbitrary characteristic, $K = \textnormal{PGL}(3, \mathcal{O})$ with $\mathcal{O}$ the ring of integers of $F$, and $(\Gamma_n)$ is a sequence of  torsion-free lattices. We prove a form of equidistribution on average for eigenfunctions whose spectral parameters lie in the tempered spectrum when the associated sequence of quotients of the Bruhat-Tits building Benjamini-Schramm converges to the building itself. This result is a higher rank non-archimedean analogue of existing results for graphs and locally symmetric spaces.

A recurring theme in the proof is the reduction of many computations to computing the sum of an exponential function over lattice points in a polytope; such expressions can subsequently be simplified using Brion's formula. Along the way of proving our main result we prove several other results which may be of independent interest including a ``degenerate'' version of Brion's formula which ``interpolates'' between the usual Brion's formula and the Ehrhart polynomial, an effective rate of convergence for the distribution of spectral parameters to the Plancherel measure under Benjamini-Schramm convergence, and a classification of relative positions of triples of points in buildings of type $\tilde{A}_2$.
\end{abstract}

\maketitle
\tableofcontents

\section{Introduction}

\subfile{sections/intro}

\section{Preliminaries} \label{preliminaries}
\subfile{sections/prelims_redo}

\section{Outline of the proof of Theorem \ref{main_thm}} \label{proof_outline}
\subfile{sections/proof_outline_redo}

\section{Brion's formula} \label{sec_brion}
\subfile{sections/polytopes_redo}


\section{Spectral bound for polytopal ball averaging operators}
\label{spectral_bound_section}
\subfile{sections/spectral_bound}

\section{Benjamini-Schramm convergence and Plancherel convergence}
\label{bs_implies_plancherel_section}
\subfile{sections/BS_convergence_fixed}

\section{Lifting the kernel function to $G/K$}
\label{kernel_section}
\subfile{sections/kernel}

\section{Changing variables in the kernel function integral}
\subfile{sections/changing_variables}

\section{The Kunze-Stein phenomenon and an ergodic theorem in the style of Nevo}
\label{sec_nevo}
\subfile{sections/kunze_stein}

\section{Classification of primitive triples of vertices in $\mathcal{B}$}
\subfile{sections/classification_redo}

\section{Geometric bound on the size of the intersection of polytopal balls} \label{sec_geometric_bound}

\subfile{sections/geometric_bound}

\section{Bounding the sum over $H_M^{\Lambda}$}
\label{sec_final_brion}
\subfile{sections/final_brion}

\printbibliography

\vspace{10mm}

\textsc{Universität Paderborn, Warburgerstr. 100, 33098 Paderborn, Germany} \\
\indent {\it Email address:} \texttt{peterson@math.upb.de}

\end{document}

%% file: sections/intro.tex
\sloppy
\subsection{Quantum ergodicity in the large eigenvalue limit}
Originally {\it quantum ergodicity} concerned eigenfunctions of the Laplacian on a closed Riemannian manifold $(M, g)$ with ergodic geodesic flow. The geodesic flow being {\it ergodic} implies that a ``classical particle'' moving along a generic geodesic is equally likely to end up everywhere on $M$ in the long run. On the other hand, a ``quantum particle'' on $M$ with wave function $\psi$ has probability measure $|\psi|^2 d \textnormal{vol}_g$ of being ``observed'' in a given region of $M$. If the geodesic flow is ergodic, we expect the quantum particle is equally likely to be observed everywhere. 

Through a procedure known as {\it quantization}, one is led to studying the Laplacian on $M$ in place of the geodesic flow. The Laplacian has an orthonormal basis of $L^2$-eigenfunctions $\{\psi_j\}$. The associated eigenvalues $0 = \lambda_1 \leq \lambda_2 \leq ...$ are all non-negative and go to infinity as $j$ goes to infinity. Snirelman \cite{snirelman}, Zelditch \cite{zelditch}, and Colin de Verdi\`ere \cite{colin_de_verdiere} proved the {\it quantum ergodicity theorem} which says that for a density-one subsequence of the eigenfunctions $(\psi_{j_k})$, the associated measures $|\psi_{j_k}|^2 d \textnormal{vol}_g$ weak-* converge to $\frac{d \textnormal{vol}_g}{\textnormal{vol}(M)}$, the normalized volume measure. In fact this statement is a consequence of the following:
\begin{theorem} [Quantum ergodicity theorem \cite{snirelman, zelditch, colin_de_verdiere}] \label{qe}
Suppose $(M, g)$ is a closed Riemannian manifold with ergodic geodesic flow. Then for every smooth test function $a \in C^\infty(M)$,
\begin{gather}
\lim_{\lambda \to \infty} \frac{1}{N(\lambda)} \sum_{\lambda_j \leq \lambda} \Big| \langle \psi_j, a \psi_j \rangle - \frac{1}{\vol(M)} \int_M a \ d \vol_g \Big|^2 = 0, \label{classical_qe}
\end{gather}
where $N(\lambda) = \#\{i : \lambda_i \leq \lambda \}$. 
\end{theorem}
\noindent In a certain sense letting the eigenvalue go to infinity is analogous to letting Planck's constant go to zero, so we expect to recover ``classical mechanical'' results (e.g. equidistribution of generic orbits when the geodesic flow is ergodic) in such a limit.

\begin{remark} \label{remark_microlocal_lift}
In fact Theorem \ref{qe} can be strengthened to replacing the summands in \eqref{classical_qe} by
\begin{gather}
\Big| \langle \psi_j, \mathcal{T} \psi_j \rangle - \int_{S^*M} \sigma^0(\mathcal{T}) d L \Big|^2, \label{microlocal_lift}
\end{gather}
where $\mathcal{T}$ is an order 0 pseudodifferential operator on $M$, $S^* M$ is the unit (co)tangent bundle of $M$ (which is the space on which the geodesic flow occurs and is ergodic), $\sigma^0(\mathcal{T})$ is the principal symbol of $\mathcal{T}$ (which is a function on $S^* M$), and $dL$ is the Liouville measure (normalized so that the total measure of $S^* M$ is 1). 
\end{remark}

\begin{remark}
Rudnick-Sarnak \cite{rudnick_sarnak} conjectured {\it quantum unique ergodicity}, namely that, if $M$ has negative curvature, then the only weak-* limit of the eigenfunction measures (or their \textit{microlocal lifts} as in Remark \ref{remark_microlocal_lift}) is the normalized volume measure on $M$ (or the Liouville measure on $S^*M$). Lindenstrauss \cite{que} proved quantum unique ergodicity for joint eigenfunctions of the Laplacian and the Hecke operators on compact arithmetic hyperbolic surfaces (see also \cite{h_cross_h, sound_que, brooks_lindenstrauss2}).
\end{remark}

\subsection{Quantum ergodicity in the Benjamini-Schramm limit}
As opposed to the original quantum ergodicity theorem, which concerns eigenfunctions of the Laplacian on a fixed manifold, more recently many authors have considered \textit{quantum ergodicity in the Benjamini-Schramm limit}, which concerns eigenfunctions of Laplacian-like operators for sequences of spaces ``converging'' to their common universal cover. More specifically, we say that a sequence of spaces \textit{Benjamini-Schramm converges} to their common (contractible) universal cover if asymptotically almost every point has arbitrarily large injectivity radius \cite{benjamini_schramm}. Benjamini-Schramm convergence may be viewed as a probabilistic version of Gromov-Hausdorff convergence of metric spaces. 

Many authors have studied the relationship between Benjamini-Schramm convergence and the distribution of eigenvalues (for example, Kesten \cite{kesten} and Mckay \cite{mckay} for regular graphs, or the more recently Abert-Bergeron-Biringer-Gelander-Nikolov-Raimbault-Samet \cite{7_samurai} for locally symmetric spaces). Such results are often closely related to representation theory in which case they often appear in the literature as {\it limit multiplicity} theorems. On the one hand quantum ergodicity in the Benjamini-Schramm limit may be seen as an extension of such questions to eigenfunctions; on the other hand it arose out of attempts to extend quantum ergodicity ideas to different contexts.

One such precedent for quantum ergodicity in the Benjamini-Schramm limit is modular forms on congruence coverings of the modular surface. Modular forms have two natural parameters: weight and level. {\it Weight} is roughly analogous to Laplacian eigenvalue and {\it level} relates to the congruence covering. Holowinsky and Soundararajan \cite{holowinsky} proved quantum unique ergodicity for fixed level and letting weight go to infinity. Fixing weight and letting level go to infinity results in a sequence of hyperbolic surfaces which Benjamini-Schramm converge to the hyperbolic plane as proven by Fraczyk \cite{fraczyk}, and quantum unique ergodicity in this context was proven by Nelson \cite{nelson} and Nelson-Pitale-Saha \cite{nelson_et_al}. These latter results are referred to in the literature as {\it quantum ergodicity in the level aspect}, and consequently this terminology is also used in place of quantum ergodicity in the Benjamini-Schramm limit.

Another motivating example of quantum ergodicity in the Benjamini-Schramm limit comes from the work of Anantharaman-Le Masson \cite{nalini_le_masson} who proved that for sequences of regular graphs Benjamini-Schramm converging to the infinite regular tree, the eigenfunctions of the adjacency operator are, on average, equidistributed in a weak sense. In particular their results include the following: 
\begin{theorem} [Quantum ergodicity on large regular graphs \cite{nalini_le_masson}; see also \cite{brooks_le_masson_lindenstrauss, nalini, nalini_sabri}] \label{qe_large_graphs}
Suppose $(G_n)$ is a sequence of $(q+1)$-regular graphs which Benjamini-Schramm converge to the $(q+1)$-regular tree. Suppose also that there is a uniform spectral gap for the adjacency operator, namely all eigenvalues other than $q+1$ are uniformly bounded away from $\pm (q+1)$. Let $a_n$ be a function on the vertices of $G_n$ such that $||a_n||_\infty \leq 1$, and let $\card(G_n)$ denote the number of vertices of $G_n$. Let $I$ be a closed subinterval of $[-2 \sqrt{q}, 2 \sqrt{q}]$ with non-empty interior. Let $\psi_1^{(n)}, \dots, \psi_{\card(G_n)}^{(n)}$ be an orthonormal basis of eigenfunctions of the adjacency operator on $G_n$ with associated eigenvalues $\lambda_1^{(n)}, \dots, \lambda_{\card(G_n)}^{(n)}$. Then
\begin{gather}
\lim_{n \to \infty} \frac{1}{N(I, G_n)} \sum_{\psi_j^{(n)}: \lambda_j^{(n)} \in I} \Big| \langle \psi_j^{(n)}, a_n \psi_j^{(n)} \rangle - \frac{1}{\card(G_n)} \sum_{\textnormal{vertices } v \in G_n} a_n(v) \Big|^2 = 0, \label{qe_large_graphs_eqn}
\end{gather}
where $N(I, G_n) = \#\{j: \lambda_j^{(n)} \in I\}$. 
\end{theorem}

\begin{remark} \label{remark_bipartite}
The assumption that all eigenvalues for all of the graphs $G_n$ in the sequence are bounded away from $-(q+1)$ in particular implies that all of the $G_n$ are non-bipartite. In fact, having $-(q+1)$ as an eigenvalue is equivalent to being bipartite. See also Remark 1.4 in \cite{nalini_le_masson}
\end{remark}

Le Masson-Sahlsten \cite{le_masson_sahlsten} noted the similarities between quantum ergodicity in the level aspect for modular forms and quantum ergodicity on large regular graphs, and they proved that for sequences of compact hyperbolic surfaces which Benjamini-Schramm converge to the hyperbolic plane, eigenfunctions of the Laplacian with eigenvalue in a fixed compact subinterval of $[1/4, \infty)$ are, on average, equidistributed in a weak sense. Parts of their technique were adapted from that of Brooks-Le Masson-Lindenstrauss \cite{brooks_le_masson_lindenstrauss} who reproved Theorem \ref{qe_large_graphs}. This technique was further adapted by Abert-Bergeron-Le Masson \cite{abert_bergeron_le_masson} to prove analogous results for sequences of compact rank one locally symmetric spaces converging to the their common universal cover; these authors also connect the original quantum ergodicity theorem (Theorem \ref{qe}) with quantum ergodicity in the Benjamini-Schramm limit. More recently Le Masson and Sahlsten \cite{le_masson_sahlsten2} have further adapted the technique to handle finite-volume non-compact hyperbolic surfaces.

\subsection{Quantum ergodicity in higher rank}

Much of the research in quantum ergodicity has been focused on the case of hyperbolic surfaces. Such manifolds are very special: their universal cover is a {\it symmetric space}. More precisely, the {\it hyperbolic plane} $\mathbb{H}$ may be realized as $\textnormal{SL}(2, \mathbb{R})/\textnormal{SO}(2)$, and each hyperbolic surface $X$ may be realized as $\Gamma \backslash \textnormal{SL}(2, \mathbb{R}) / \textnormal{SO}(2)$ for some discrete subgroup $\Gamma$ (which must be a cocompact lattice in case $X$ is compact). Furthermore, the geodesic flow on the unit tangent bundle of $X$ may be identified with the right action of the subgroup $\bigl( \begin{smallmatrix}e^{t/2} & 0\\ 0 & e^{-t/2}\end{smallmatrix}\bigr)$ on $\Gamma \backslash \textnormal{SL}(2, \mathbb{R}) / \{\pm I\}$.

Now suppose $G$ is a non-compact semisimple Lie group with Lie algebra $\mathfrak{g}$. Using a Cartan involution, we may write $\mathfrak{g} = \mathfrak{k} \oplus \mathfrak{p}$ where $K = \exp(\mathfrak{k})$ is a maximal compact subgroup. Let $\mathfrak{a} \subset \mathfrak{p}$ be a maximal toral subalgebra, and let $A = \exp(\mathfrak{a})$. Then $G/K$ has a natural $G$-invariant metric with respect to which it is a contractible manifold with non-positive curvature; such spaces are called symmetric spaces of non-compact type, and $K$ is clearly the stabilizer of the point $1 K$. The orbit of $1 K$ under $A$ is a flat subspace; it is an example of a {\it maximal flat} in $G/K$. 

Let $\Gamma < G$ be a lattice. The manifold $\Gamma \backslash G /K$ is called a {\it locally symmetric space}. If the (real) rank of $G$, or equivalently the dimension of $\mathfrak{a}$, is greater than one, then the geodesic flow on any associated locally symmetric space is never ergodic \cite{bekka_mayer}. On the other hand, we may consider the $A$-action on the space $\Gamma \backslash G / M$ where $M = Z_K(A)$; this double coset space may be viewed as the ``bundle of oriented flats with basepoint'' over the underlying locally symmetric space (see, e.g., Section 5.3 of \cite{silberman_venkatesh}). This $A$-action is ergodic and reduces to the geodesic flow in rank one \cite{bekka_mayer}; in the literature this action is also referred to as the Weyl chamber flow.

The validity of the strengthened version of Theorem \ref{main_thm} as in Remark \ref{remark_microlocal_lift} is in fact known to be equivalent to the ergodicity of the geodesic flow \cite{zelditch2}. Hence to extend the ideas of quantum ergodicity to a higher rank setting, it no longer suffices to simply consider eigenfunctions of the Laplacian. Lindenstrauss \cite{h_cross_h} suggested eigenfunctions of $\mathcal{D}(G/K)$, the {\it algebra of invariant differential operators on $G/K$}, as the object of study for quantum ergodicity on higher rank locally symmetric spaces. This perspective was taken up in \cite{silberman_venkatesh, nalini_silberman, silberman_venkatesh2, brumley_matz}. In particular, Brumley-Matz \cite{brumley_matz} investigated quantum ergodicity in the Benjamini-Schramm limit result for sequences of locally symmetric spaces $\Gamma_n \backslash \textnormal{SL}(d, \mathbb{R}) / \textnormal{SO}(d)$ which Benjamini-Schramm converge to $\textnormal{SL}(d, \mathbb{R})/\textnormal{SO}(d)$.

In rank one, the algebra $\mathcal{D}(G/K)$ is generated by the Laplacian \cite{knapp}. More generally, via the Harish-Chandra isomorphism, one may identify $\mathcal{D}(G/K)$ with a polynomial ring whose associated variety may be identified with $\mathfrak{a}_{\mathbb{C}}^*/W$ where $W = N_G(\mathfrak{a})/Z_G(\mathfrak{a})$ is the Weyl group. Given a joint eigenfunction $\psi$ of $\mathcal{D}(G/K)$ acting on $\Gamma \backslash G /K$, we get an associated homomorphism $\chi \in \textnormal{Hom}_{\textnormal{$\mathbb{C}$-alg.}}(\mathcal{D}(G/K), \mathbb{C})$ via $D \psi = \chi(D) \psi$, and hence an associated point $\nu \in \mathfrak{a}_{\mathbb{C}}^*/W$ which is called the {\it spectral parameter} of $\psi$. The locus $i \mathfrak{a}^*/W \subset \mathfrak{a}_{\mathbb{C}}^*/W$ plays a distinguished role and is called the {\it tempered spectrum}. In case $G = \textnormal{SL}(2, \mathbb{R})$, the tempered spectrum may be identified under a natural mapping with $[1/4, \infty)$ which is often referred to as the tempered spectrum of the Laplacian on $\mathbb{H}$. 

We now state the main theorem of Brumley-Matz \cite{brumley_matz}, but we first recall that a sequence of lattices in $G$ is called {\it uniformly discrete} if there is a universal lower bound on the injectivity radii of the associated locally symmetric spaces.
\begin{theorem} [Brumley-Matz \cite{brumley_matz}] \label{brumley_matz}
Let $d \geq 3$. Let $G = \textnormal{SL}(d, \mathbb{R})$ and $K = \textnormal{SO}(d)$. Suppose $\Gamma_n < G$ is a uniformly discrete sequence of torsion-free cocompact lattices. Let $Y_n = \Gamma_n \backslash G / K$. Suppose $\vol(Y_n) \to \infty$. Let $a_n$ be a measurable function on $Y_n$ such that $||a_n||_\infty \leq 1$. Let $\{\psi_j^{(n)} \}$ be an orthonormal basis for $L^2(Y_n)$ of eigenfunctions of $\mathcal{D}(G/K)$ with associated spectral parameters $\{\nu_j^{(n)}\}$. Then, there is $\rho > 1$ such that for sufficiently regular $\nu \in i \mathfrak{a}^*$, we have
\begin{gather*}
\lim_{n \to \infty} \frac{1}{N(B_0(\nu, \rho), Y_n)} \sum_{j: \nu_j^{(n)} \in B_0(\nu, \rho)} \Big| \langle \psi_j^{(n)}, a_n \psi_j^{(n)} \rangle - \frac{1}{\vol(Y_n)} \int_{Y_n} a_n \ d\vol_{Y_n} \Big|^2 = 0,
\end{gather*}
where $B_0(\nu, \rho) = \{\lambda \in i \mathfrak{a}^*: ||\lambda - \nu||_2 \leq \rho \}$ is the ball of radius $\rho$ centered at $\nu$ in the tempered spectrum, and $N(B_0(\nu, \rho), Y_n) = \#\{j : \nu_j^{(n)} \in B_0(\nu,\rho) \}$.
\end{theorem}
The assumption that $\vol(Y_n) \to \infty$ is in fact known to be equivalent to Benjamini-Schramm convergence in rank at least 2 \cite{7_samurai}. Furthermore, in contrast to Theorem \ref{qe_large_graphs} and the work of \cite{le_masson_sahlsten, abert_bergeron_le_masson} for rank one locally symmetric spaces, there is no uniform spectral gap assumption needed as it is in fact automatic in rank at least 2 by property (T) \cite{property_t_book}.

\subsection{Main result: quantum ergodicity in the Benjamini-Schramm limit for the Bruhat-Tits building associated to $\textnormal{PGL}(3, F)$}

\textit{Bruhat-Tits buildings} are infinite simplicial complexes constructed from reductive algebraic groups over non-archimedean local fields \cite{bruhat_tits}. The simplest example is the Bruhat-Tits building associated to $\textnormal{SL}(2, \mathbb{Q}_p)$, which is the infinite $(p+1)$-regular tree. Bruhat-Tits buildings may be viewed as non-archimedean analogues of symmetric spaces of non-compact type. On the other hand, their quotients may be seen as ``higher rank'' generalizations of regular graphs. Such quotients have also been studied recently because, in certain cases, they provide examples of high-dimensional expanders known as Ramanujan complexes (see, e.g. \cite{lubotzky_samuels_vishne, ramanujan_complex2, ramanujan_complex3}). 

Suppose $F$ is a non-archimedean local field with $\mathcal{O}$ its ring of integers and $q$ the order of the residue field. Suppose $G$ is a reductive algebraic group over $F$, and $K$ is a hyperspecial maximal compact subgroup. Let $\mathcal{B}$ be the associated building. Then $K$ is the stabilizer of a unique special vertex $x_0 \in \mathcal{B}$. There is a correspondence between maximal $F$-split tori in $G$ and so-called {\it apartments} in $\mathcal{B}$; these apartments are the analogues for Bruhat-Tits buildings of maximal flats in symmetric spaces. Let $T < G$ be a maximal $F$-split torus whose associated apartment contains $x_0$. Let $\Gamma < G$ be a lattice. Let $M = Z_K(T)$ and let $A = T/(T \cap K)$. Similarly to the case of (locally) symmetric spaces, there is an ergodic right $A$-action on $\Gamma \backslash G / M$. This double coset space may be viewed as the ``bundle of oriented apartments with special vertex basepoint'' over $\Gamma \backslash \mathcal{B}$. It is this ergodic action which in some sense we are ``quantizing'' in this work.

We shall be particularly concerned with the Bruhat-Tits building associated to $G = \textnormal{PGL}(d, F)$ (specifically the case of $d = 3$). Then $K = \textnormal{PGL}(d, \mathcal{O})$, and $G/K$ may be identified with the vertices of $\mathcal{B}$. The \textit{spherical Hecke algebra} $H(G, K)$ is the analogue of $\mathcal{D}(G, K)$ in this context, and it acts on functions on $G/K$ (or quotients thereof). In the case of $d = 2$, the algebra $H(G, K)$ is generated by one element whose associated action on $G/K$ is equivalent to the adjacency operator on the infinite $(q+1)$-regular tree \cite{serre_trees}. If $\Gamma < G$ is a cocompact torsion-free lattice, then $Y = \Gamma \backslash G /K$ is a finite simplicial complex. The space $L^2(Y)$ has an orthonormal basis of eigenfunctions $\{\psi_j\}$ of $H(G, K)$. By the Satake isomorphism, the spherical Hecke algebra is isomorphic to the coordinate ring of a variety $\Omega$, and hence given an eigenfunction $\psi_j$, it makes sense to talk about its associated {\it spectral parameter} $\nu_j \in \Omega$. Inside of $\Omega$ we have $\Omega^+$ which is the collection of all spectral parameters that can arise from unitary representations of $H(G, K)$. Inside of $\Omega^+$ there is a distinguished sublocus called the \textit{tempered spectrum}, denoted by $\Omega^+_{\textnormal{temp}}$ which corresponds to the spectral parameters that arise from the action of $H(G, K)$ on $L^2(G/K)$. The space $\Omega^+$ carries a natural probability measure called the \textit{Plancherel measure} $\mu$ which is supported on $\Omega^+_{\textnormal{temp}}$. In the case of $d = 2$, the tempered spectrum is transformed into $[-2 \sqrt{q}, 2 \sqrt{q}]$ under a natural mapping; this locus is often referred to as the tempered spectrum of the adjacency operator on the $(q+1)$-regular tree, and it appeared in Theorem \ref{qe_large_graphs}.

There is a natural coloring on the vertices of $\mathcal{B}$ which is preserved by the $\textnormal{PSL}(d, F)$-action, but not by the $\textnormal{PGL}(d, F)$-action. If $t = [\Gamma:\Gamma \cap \textnormal{PSL}(d, F) \cdot K]$, then there are $d/t$ special spectral parameters which arise from the fact such a $\Gamma$-action on the vertices of $\mathcal{B}$ preserves a coloring with $d/t$ colors; we call these the {\it coloring spectral parameters}. Each coloring spectral parameter has an associated {\it coloring eigenfunction}. One of these coloring eigenfunctions is the constant function; we call this eigenfunction the {\it trivial eigenfunction} and all other coloring eigenfunctions the {\it non-trivial coloring eigenfunctions}. We wish to ignore the coloring eigenfunctions in the analysis. One way of doing this is to restrict to test functions which are orthogonal to the non-trivial coloring eigenfunctions. A stronger assumption one may make is that $t = d$ above, i.e. that there are no non-trivial coloring eigenfunctions. In case $d$ is prime such as $d = 3$, this is equivalent to requiring that $\Gamma$ not be contained in $\textnormal{PSL}(d, F) \cdot K$. In case $d = 2$, this is in turn equivalent to requiring that the associated graphs not be bipartite. See also Remark \ref{remark_bipartite}.

Our main result is the following:

\begin{theorem} \label{main_thm}
Let $G = \textnormal{PGL}(3, F)$ and $K = \textnormal{PGL}(3, \mathcal{O})$, where $F$ is a non-archimedean local field of arbitrary characteristic and $\mathcal{O}$ is its ring of integers. Let $\Gamma_n < G$ be a sequence of  torsion-free lattices. Let $Y_n = \Gamma_n \backslash G / K$. Suppose $\card(Y_n) \to \infty$. Suppose $\Theta \subset \Omega_{\textnormal{temp}}^+$ is such that the Plancherel measure of $\Theta$ is positive, the Plancherel measure of the boundary of $\Theta$ is zero, and the closure of $\Theta$ does not intersect a certain codimension one exceptional locus $\Xi$. Let $\{\psi_j^{(n)}\}$ denote an orthonormal basis of eigenfunctions of $H(G, K)$ acting on $L^2(Y_n)$. Let $a_n$ be a function on $Y_n$ such that $||a_n||_\infty \leq 1$ and $a_n$ is orthogonal to all non-trivial coloring eigenfunctions. Then
\begin{gather}
\lim_{n \to \infty} \frac{1}{N(\Theta, Y_n)} \sum_{\psi_j^{(n)}: \nu_j^{(n)} \in \Theta} \Big| \langle \psi_j^{(n)}, a_n \psi_j^{(n)} \rangle - \frac{1}{\card(Y_n)} \sum_{\textnormal{vertices } v \in Y_n} a_n(v) \Big|^2 = 0, \label{main_thm_eqn}
\end{gather}
where 
\begin{gather}
N(\Theta, Y_n) = \#\{j: \nu_j^{(n)} \in \Theta\}. \label{n_theta}
\end{gather}
\end{theorem}
\noindent The notation $\card(\cdot)$ refers to the cardinality of a set. The codimension one exceptional locus $\Xi$ is defined in Section \ref{sec_tempered_pgl_3}.

This theorem is really a consequence of the following quantitative estimates:

\begin{theorem} \label{main_thm2}
Suppose $\Gamma < G$ is a torsion-free lattice and $Y = \Gamma \backslash G / K$. Suppose $\Theta \subset \Omega^+_{\textnormal{temp}}$ is such that its closure does not intersect $\Xi$. Let $\{\psi_j\}$ denote an orthonormal basis of eigenfunctions of $H(G, K)$ acting on $L^2(Y)$. Let $a \in L^\infty(Y)$ be mean-zero and orthogonal to all coloring eigenfunctions. There exist universal constants $C_1, C_2$, and constants $M_0 = M_0(\Theta)$ and $C = C(\Theta)$ depending on $\Theta$ such that for all $M \geq M_0$,
\begin{gather*}
\sum_{\psi_j: \nu_j \in \Theta} \Big| \langle \psi_j, a \psi_j \rangle \Big|^2 \leq C \Big( \frac{C_1 ||a||_2^2}{M} + \frac{C_2 M^2 q^{4 M}||a||_\infty^2}{\textnormal{InjRad}(Y)^2} \card(\{y \in Y: \textnormal{InjRad}_Y(y) \leq M \}) \Big).
\end{gather*}
\end{theorem}

\begin{theorem} \label{thm_bs_conv}
Suppose $E \subset \Omega^+$ is a measurable subset such that $\mu(\partial E) = 0$. Then for every $\varepsilon > 0$ there exists an $R = R(\varepsilon)$ and a $C = C(\varepsilon)$ such that for all $Y = \Gamma \backslash G / K$ with $\Gamma$ a torsionfree lattice we have
\begin{gather}
\Big| \frac{N(E, Y)}{\card(Y)} - \mu(E) \Big| \leq C \frac{R q^{2 R}}{q^{2 \cdot \textnormal{InjRad}(Y)}} \frac{\card(\{y \in Y: \textnormal{InjRad}_Y(y) \leq R\})}{\card(Y)} + \varepsilon.
\end{gather}
\end{theorem}

\begin{remark} \label{fewer_hypotheses}
Theorem \ref{main_thm} may be viewed simultaneously as a higher rank analogue of Theorem \ref{qe_large_graphs} and a non-archimedean analogue of Theorem \ref{brumley_matz}. The assumption that $\card(Y_n) \to \infty$ is in fact known to be equivalent to Benjamini-Schramm convergence in this setting by Gelander-Levit \cite{gelander_levit}. No spectral gap assumption need be made as it is automatic by Property (T) \cite{property_t_book}. Since all our lattices are torsion-free, they are necessarily also cocompact and the associated quotients of the building are finite simplicial complexes. To the best of our knowledge, our result is the first quantum ergodicity result of any kind in a higher rank non-archimedean setting.
\end{remark}

\begin{remark}
The analogue of the exceptional locus (or, more precisely, of what is called $\Xi_1$ in Section \ref{sec_tempered_pgl_3}) for the $(q+1)$-infinite regular tree is the set $\{-2 \sqrt{q}, 2 \sqrt{q}\}$, for the hyperbolic plane is $\{1/4\}$, and for $\textnormal{SL}(d, \mathbb{R})/\textnormal{SO}(d)$ is those points in $i \mathfrak{a}^*/W$ whose stabilizer is non-trivial (i.e. the non-regular points of a Weyl chamber). Notice that in Theorem \ref{qe_large_graphs}, subsets of the tempered spectrum meeting this exceptional locus are allowed. However, in all of the aforementioned works regarding symmetric spaces \cite{le_masson_sahlsten, abert_bergeron_le_masson, brumley_matz}, one does not allow subsets of the tempered spectrum intersecting the exceptional locus. We suspect that Theorem \ref{main_thm} is still true without the assumption about avoiding $\Xi$, but we have not been able to remove this assumption.
\end{remark}

\begin{remark}
There are two natural topologies on $\Omega^+_{\text{temp}}$ and $\Omega^+$: one coming from the natural embedding of the space into a larger complex algebraic variety $\Omega$, and the other being the Fell topology. In Section \ref{sec_compare_topologies} we show that these two topologies are actually the same. 
\end{remark}

\begin{remark}
The assumption about the test functions being orthogonal to the coloring eigenfunctions can be removed if we place additional hypotheses either about the set $\Theta$ or the orthonormal basis $\{\psi_j^{(n)} \}$. This is discussed in Section \ref{coloring_and_que}. This assumption can be seen as arising from the fact that the Weyl chamber flow on the quotient of the Bruhat-Tits building which has a non-trivial coloring is not mixing; for example on a bipartite regular graph colored with colors 1 and 2, any geodesic starting at color 1 will always be at a vertex of color 1 after an even number of steps and at a vertex of color 2 after an odd number of steps.
\end{remark}

\begin{remark}
Our paper is significantly longer than previous papers proving analogous results in other homogeneous settings \cite{nalini_le_masson, brooks_le_masson_lindenstrauss, le_masson_sahlsten, brumley_matz}. This is for a number of reasons.

In order to be able to prove a result about arbitrary torsionfree lattices in $\textnormal{PGL}(d, F)$, one must deal with the coloring eigenfunctions which have no archimedean analogue. This in turn requires an analysis of the geometry of the tempered spectrum and an associated $\mathbb{Z}/d \mathbb{Z}$-action on the tempered spectrum and on the coloring eigenfunctions. We have not found a sufficient exposition on this elsewhere and thus have included one in Section \ref{sec_coloring_eigenfunctions}. The analysis of the coloring eigenfunctions also reveals an interesting form of ``quantum scarring''/failure of quantum unique ergodicity that can happen because of a large degeneracy of a specific eigenspace. We found this to be worth noting and thus included exposition about it in Section \ref{coloring_and_que}.

We use polytopal methods at many places in this work; in particular we make repeated use of Brion's formula involving summing exponential functions over lattice points in polytopes \cite{brion}. Brion's formula may be seen as a higher dimensional generalization of geometric series as well as a combinatorial version of the method of stationary phase/Laplace's method. However, Brion's formula has certain technical assumptions which are not satisfied in all of the situations in which we wish to apply it. For this reason we derive in Section \ref{sec_brion} a ``degenerate'' version of Brion's formula; this formula is particularly useful when dealing with a family of polytopes ``of the same type'' (defined in Section \ref{sec_polytope_type}). In the maximally degenerate case this formula reduces to the Ehrhart polynomial. We believe that this degenerate version of Brion's formula may have other applications.

At a certain point in the analysis we need to understand the relationship between Benjamini-Schramm convergence and convergence of the distribution of spectral parameters to the Plancherel measure. Using a result of Deitmar \cite{deitmar} together with the Sauvageot density principle \cite{sauvageot}, one would obtain the necessary result at least in characteristic zero. However, recently there has been concern about a gap in the proof of the Sauvageot density principle \cite{nelson_venkatesh} and as such we avoid using it. We instead prove a version of the Sauvageot density principle for spherical representations; the proof works in both zero and positive characteristic and allows us to get an effective bound on the rate of convergence to the Plancherel measure; this is discussed in Section \ref{bs_implies_plancherel_section}.

The work of \cite{le_masson_sahlsten} and \cite{brumley_matz} relied on the Nevo ergodic theorem which is stated in the literature only for real Lie groups \cite{nevo}. However, a key ingredient in the proof of this result is the Kunze-Stein phenomenon which has subsequently been proven for simply connected semisimple algebraic groups over non-archimedean local fields of arbitrary characteristic \cite{veca}. By first proving that the Kunze-Stein phenomenon is stable under compact group extensions, we can prove that $\textnormal{PGL}(d, F)$ (which is not simply connected) has the Kunze-Stein phenomenon. We can then retrace Nevo's proof method to prove the Nevo ergodic theorem for $\textnormal{PGL}(d, F)$. This is done in Section \ref{sec_nevo}.

We originally followed the proof strategy of \cite{brumley_matz} for Theorem \ref{brumley_matz} very closely, and in particular we also utilize \textit{polytopal balls} (discussed in Section \ref{sec_notation}). At a certain step in the proof, referred to in this work as the \textit{geometric bound}, one must bound the volume of the intersection of two polytopal balls with different centers. In the non-archimedean setting, if we use the analogue of the exact polytope used in \cite{brumley_matz}, then the volume of intersection of the associated polytopal balls is too big to be able to complete the proof using the same proof strategy as in \cite{brumley_matz}. However, by choosing the ``right'' polytope, the volume of the intersection of the resulting polytopal balls is small enough to be able to complete the proof. Thus part of the difficulty and novelty in this work is finding the right polytope to get the requisite geometric bound, and, for related reasons, the method of proof for the geometric bound differs greatly from the proof in \cite{brumley_matz}. This subtlety is discussed further in Section \ref{compare_brumley_matz}.

Our method of proof of the geometric bound roughly goes as follows: we first classify relative positions of triples of points in the building (Section \ref{sec_classify_triples}), we use this classification to set up a ``polytopal parametrization'' of points in the intersection of two polytopal balls (Sections \ref{sec_coordinatize} and \ref{sec_polytope_parametrization}), we count the number of points which receive the same coordinates under this parametrization and find that it is an exponential in these coordinates (Section \ref{sec_coordinatize}), and finally we use the degenerate Brion's formula to sum up this exponential function over the lattice points in the parametrizing polytope (Section \ref{sec_pf_e_m_lambda_bound}).

When we classify relative positions of triples of points, we first classify pairs of sectors in the building based at the same vertex whose relative position in the local spherical building is ``nearly opposite'', i.e. one step away from the longest element in the Coxeter group. Such sectors can either remain nearly opposite in the spherical building at infinity or else become opposite past some definite point. We provide both a geometric and algebraic explanation of this phenomenon. This is closely related to the relationship between the Bruhat decomposition of elements in $\textnormal{PGL}(3, \mathcal{O})$ over the residue field $\mathbb{F}_q$ vs. over the underlying local field $F$. This is discussed in Sections \ref{sec_classify_nearly_opposite} and \ref{sec_algebraic_nearly_opposite}.
\end{remark}

\subsection*{Acknowledgements}
\subfile{./acknowledgements}

%% file: sections/acknowledgements.tex
I thank my advisor, Ralf Spatzier, for introducing me to quantum (unique) ergodicity and for his incredible support and encouragement throughout every stage of this project. I would also like to thank Farrell Brumley, Rahul Dalal, Stephen DeBacker, and Simon Marshall for listening to my ideas and providing helpful feedback. This work was partially supported by the National Science Foundation Grant DMS-2003712, RTG Grant DMS-1840234, and Deutsche Forschungsgemeinschaft (DFG, German Research Foundation) Grant SFB-TRR 358/1 2023 - 491392403.

%% file: sections/prelims_redo.tex
\sloppy

\subsection{Representation theory} \label{sec_rep_background}
A more elaborated discussion of the representation theoretic background used in this work may be found in Appendix A of \cite{my_thesis}.

Let $F$ be a non-archimedean local field of arbitrary characteristic, $\mathcal{O}$ its ring of integers, $\varpi$ a uniformizer of $\mathcal{O}$, and $q$ the order of the residue field. Let $G = \textnormal{PGL}(d, F)$ and $K = \textnormal{PGL}(d, \mathcal{O})$. We denote the Haar measure on $G$ by $\vol(\cdot)$ normalized so that $\vol(K) = 1$. Let $T < G$ denote the subgroup of diagonal matrices. Let $A < T$ denote the subgroup of matrices all of whose diagonal entries are powers of $\varpi$. Let $A^+ \subset A$ denote those elements for which the powers of $\varpi$ along the diagonal are weakly decreasing. To each element $a \in A$, we can associate a tuple $\lambda = (\lambda_1, \dots, \lambda_d)$ by recording the powers of $\varpi$ along the diagonal. Note that this is only well-defined up to shifting all entries by the same integer. We let $\varpi^\lambda$ denote the matrix $\text{diag}(\varpi^{\lambda_1}, \dots, \varpi^{\lambda_d})$. If $M_1, M_2 < G$ are subgroups, we say that a function on $G$ is \textit{$(M_1, M_2)$-invariant} if it is invariant under $M_1$-multiplication on the left and $M_2$-multiplication on the right. 

The \textit{spherical Hecke algebra} $H(G, K)$ is defined as all compactly supported $(K, K)$-invariant functions on $G$ with convolution product. The Cartan decomposition tells us that $G$ is the disjoint union of the double cosets $K \varpi^{\lambda} K$ with $\varpi^\lambda \in A^+$ (see, e.g., \cite{macdonald-symmetric}, p. 294). Hence each element in $H(G, K)$ is completely determined by its restriction to $A^+$. The Satake isomorphism  tells us that we have a $\mathbb{C}$-algebra isomorphism (\cite{macdonald-symmetric}, p. 296-297; see also Appendix A.2.4 in \cite{my_thesis}):
\begin{gather} \label{satake}
    H(G, K) \simeq \mathbb{C}[x_1^{\pm}, \dots, x_d^{\pm}]^{\mathfrak{S}_d}/(q^{-d(d-1)/2} x_1 \dots x_d - 1).
\end{gather}

A \textit{spherical function} on $G$ is a complex-valued continuous function $\omega$ on $G$ which is $(K, K)$-invariant, equal to 1 at the identity, and is a joint eigenfunction of convolution with every element in $H(G, K)$. By recording the eigenvalue, we obtain a $\mathbb{C}$-algebra homomorphism to $\mathbb{C}$. In fact, in this way there is a bijection between spherical functions and the set $\textnormal{Hom}_{\mathbb{C}-\text{alg.}}(H(G, K), \mathbb{C})$ (\cite{macdonald-spherical-function}, Prop. (1.2.6)). By the nullstellensatz and \eqref{satake}, we also have a bijection to points on an algebraic variety $\Omega$ which is a subvariety of the symmetric product of $d$ copies of $\mathbb{C}^\times$. By changing variables to $x_i = q^{\frac{1}{2}(d-1) + s_i}$ with $s_i \in \mathbb{C}$ (mod $2 \pi i \mathbb{Z}/\ln(q)$), we can parametrize points on $\Omega$ by tuples $(q^{s_1}, \dots, q^{s_d})$ such that $\prod_j q^{s_j} = 1$, modulo permutation. 

We say that a representation $(\rho, V)$ of $G$ is \textit{smooth} if for each $v \in V$, $\text{Stab}_G(v)$ contains an open subgroup. We say it is \textit{admissible} if $V^H$ is finite dimensional for each compact open subgroup $H$. A smooth, admissible, irreducible representation which contains a non-zero $K$-fixed vector is called a \textit{spherical representation}. For spherical representations, $V^K$ is one-dimensional (\cite{cartier}, p. 152).

Given a smooth, admissible representation, we get an induced action of $H(G, K)$. Furthermore each $f \in H(G, K)$ projects onto $V^K$ (\cite{cartier}, p. 117). Hence in the case of a spherical representation, we can associate an element in $\mathbb{C}$ to each element of $H(G, K)$ in such a way that we get a $\mathbb{C}$-algebra homomorphism. In fact, all $\mathbb{C}$-algebra homomorphisms arise in this way from spherical representations (\cite{cartier}, p. 152). 

Let $\chi$ be a character of $T$ which is invariant under $K \cap T$ (i.e. an \textit{unramified character}); this implies that the character only depends on $T/(K \cap T) \simeq A \simeq \mathbb{Z}^{d-1}$. Let $B < G$ be the subgroup of upper triangular matrices and let $\Delta$ be its modular character, which also defines an unramified character on $T$. Let $I_\chi$ be the $G$-representation obtained by induction of $\chi \cdot \Delta^{1/2}$ from $T$ to $G$. Let $\mathfrak{S}_d$ denote the symmetric group on $d$ elements. There is a natural $\mathfrak{S}_d$-action on $T$ which in turn induces an action on the characters of $T$. The representations $I_\chi$ and $I_{\chi'}$ are isomorphic if and only if $\chi' = \sigma.\chi$ for some $\sigma \in \mathfrak{S}_d$. Every spherical representation of $G$ is isomorphic to a subquotient of some $I_\chi$, with $\chi$ unique up to the $\mathfrak{S}_d$-action (\cite{cartier}, p. 152).

We shall be particularly interested in irreducible unitary representations of $G$. Such a representation contains a smooth, admissible, irreducible representation as a dense subspace (\cite{cartier}, Cor. 2.3). The collection of all irreducible, unitary representations is denoted $\hat{G}$ and carries a natural topology called the {\it Fell topology}. A key notion in defining the Fell topology is that of {\it weak containment}; this is in turn defined via {\it matrix coefficients}, namely functions on $G$ of the form $\langle v_1, \rho(g). v_2 \rangle$ for $v_1, v_2 \in \mathcal{H}$, and $\rho$ a unitary representation of $G$ on a Hilbert space $\mathcal{H}$. Roughly speaking, an irreducible unitary representation $\rho_1$ is weakly contained in another unitary representation $\rho_2$, denoted $\rho_1 \prec \rho_2$, if every matrix coefficient for $\rho_1$ can be well-approximated on compact sets by matrix coefficients for $\rho_2$ (see \cite{property_t_book}, Appendix F). 

An irreducible unitary representation containing a non-zero $K$-fixed vector is called {\it class 1}. Hence to each class 1 representation we can associate a spherical representation to which we can in turn associate an unramified torus character $\chi$. If this character is unitary then we say that the associated class 1 representation is in the {\it principal series}; in fact all spherical representations with unitary torus characters are unitarizable (\cite{macdonald-spherical-function}, Prop. 3.3.1). All other class 1 representations are called the {\it complementary series}. We let $\Omega^+$ denote the collection of all class 1 representations.

Now suppose $\Gamma < G$ is a cocompact lattice. Then $L^2(\Gamma \backslash G)$ decomposes as a countable direct sum of orthogonal irreducible unitary representations (\cite{deitmar_echterhoff}, Theorem 9.2.2). The subspace $L^2(\Gamma \backslash G /K)$ is preserved by convolution on the right with elements in $H(G, K)$, which all act as commuting normal operators. Therefore $L^2(\Gamma \backslash G / K)$ has an orthonormal basis of joint eigenfunctions of $H(G, K)$. In fact there is a bijection between these eigenspaces of $H(G, K)$ and the isotypic components of class 1 representations that show up in the decomposition of $L^2(\Gamma \backslash G)$ into irreducibles. The dimension of each eigenspace is equal to the multiplicity of the corresponding class 1 representation. Hence to each joint eigenfunction of $H(G, K)$ we can associate a point in $\Omega$ which we call the {\it spectral parameter} of the eigenfunction (see, for example, the discussion in Section 2 of \cite{lubotzky_samuels_vishne}, or Prop. A.15 in \cite{my_thesis}). When we coordinatize this point as $x_i = q^{\frac{1}{2}(d-1) + s_i}$, we call the tuple $(q^{s_1}, \dots, q^{s_d})$ the {\it Satake parameters} of the spectral parameter/spherical representation.

The regular representation $L^2(G)$ has a decomposition as a $G \times G$ representation as a direct integral over all representations of the form $\mathcal{H}_\rho \otimes \mathcal{H}_{\rho^*}$, where $(\rho, \mathcal{H}_\rho)$ is an irreducible unitary representation with dual representation $(\rho^*, \mathcal{H}_{\rho^*})$, with respect to a specific measure on the collection of all irreducible unitary representations of $G$ (\cite{folland}, Theorem 7.44). This measure is called the {\it Plancherel measure}, denoted $\mu$. Its restriction to class 1 representations is explicitly known and is supported on the principal series (\cite{macdonald-spherical-function}, Theorem 5.1.2). Any representation lying in the support of the Plancherel measure is called a {\it tempered representation}. The sublocus of $\Omega$ corresponding to the principal series is called the {\it tempered spectrum}, denoted $\Omega^+_{\textnormal{temp}}$. Being in the tempered spectrum is equivalent to being weakly contained in the regular representation. The Satake parameters of points in the tempered spectrum are $\{(q^{s_1}, \dots, q^{s_d}) : q^{s_d} \in S^1, \prod_j q^{s_j} = 1\}$. The Plancherel measure is absolutely continuous with respect to the Lebesgue measure on this collection of Satake parameters. 


\subsection{Coxeter complexes} The material for this section may be found in Chapter 2 of \cite{ronan}.

A {\it Coxeter group} $W$ is any group which has a presentation of the form $W = \langle r_i | r_i^2 = (r_i r_j)^{m_{ij}} = 1 \text{ for all } i, j \in I \rangle$, where $I$ is a finite set. The $m_{i j}$ are also allowed to be $\infty$. The $r_i$'s are called the {\it Coxeter generators}. We shall from now on assume that $W$ is irreducible, namely that it does not split as a direct product of smaller Coxeter groups in such a way that each factor has its Coxeter generators as a subset of $\{r_i\}_{i \in I}$.

We may construct a simplicial complex from $W$. Each $w \in W$ defines a top-dimensional cell of dimension $|I| - 1$. These cells are called {\it chambers}. Two chambers $w_1$ and $w_2$ are {\it $i$-adjacent} if $w_2 = w_1 r_i$. In this case the two chambers share a codimension one face in the complex. Codimension one faces are called {\it panels}. The group $W$ acts simply transitively on the chambers in such a way that the $I$-valued adjacency relations are preserved (multiplication on the left). This cell complex is called the {\it Coxeter complex} $X$.

Given two chambers in $X$ associated to elements $w_1, w_2 \in W$, we may define a {\it $W$-valued metric} $d_W(w_1, w_2) = w_1^{-1} w_2$ which is preserved by the left $W$-action. A {\it gallery} between chambers $w_1$ and $w_2$ is any sequence of adjacent chambers connecting $w_1$ and $w_2$. A gallery is {\it minimal} if it is as short as any other gallery connecting $w_1$ and $w_2$. The {\it length} of an element $w \in W$ is the length of any minimal gallery from $1$ to $w$ in $X$. In general there may be several minimal galleries between two chambers. A subset $Y$ of chambers of $X$ is called {\it combinatorially convex} if every minimal gallery connecting two chambers in $Y$ lies entirely in $Y$. The {\it combinatorial convex hull} of two chambers $w_1$ and $w_2$ is the intersection of all combinatorially convex subsets of $X$ containing $w_1$ and $w_2$. 

A {\it reflection} $r \in W$ is any conjugate of a Coxeter generator. Its {\it wall} $M_r$ is all simplices in $X$ fixed by $r$. We can define an equivalence relation on chambers by specifying two chambers to be equivalent if any gallery connecting them never crosses $M_r$. This partitions chambers into two distinct sets. Each such set is called a {\it root}. 

When the group $W$ is finite it is called a {\it spherical Coxeter group}. Such groups have a unique element of longest length. Because of this each chamber in $X$ has a unique {\it opposite chamber} obtained by multiplying on the right by this longest element. 

The most important spherical Coxeter group in this work is the symmetric group $\mathfrak{S}_d$ with Coxeter generators $(1 \ 2), (2 \ 3), \dots, (d-1 \ d)$. The element of longest length is the flip permutation $(1 \ d) (2 \ d-1) \dots$. In the case of $\mathfrak{S}_3$, the associated Coxeter complex may be visualized as a partition of the unit circle into six equal pieces.

Another important class of Coxeter groups are those which have faithful representations as isometries of some $d$-dimensional Euclidean space $E^d$ in such as way that each Coxeter generator corresponds to reflection across some affine hyperplane. Such Coxeter groups are called {\it affine} and have $d+1$ Coxeter generators. If we consider all hyperplanes obtained by the orbit under $W$ of the generating hyperplanes, we obtain a tesselation of $E^d$ by simplices. This tesselation, thought of as a simplicial complex, is isomorphic to the associated Coxeter complex. 

Because affine Coxeter groups embed into the isometry group of $E^d$, namely $\mathbb{R}^d \rtimes \textnormal{O}(d)$, we have a natural homomorphism from $W$ to $\textnormal{O}(d)$. The image of $W$ is in fact a spherical Coxeter group denoted $W_0$ and has $d$ Coxeter generators. If $w_1, \dots, w_d \in W$ are reflections whose images generate $W_0$, then there is a unique vertex in the Coxeter complex for $W$ which is fixed by all of these elements. Such vertices are called {\it special vertices}. The link of a special vertex is isomorphic to the Coxeter complex for $W_0$. 

Each panel of a chamber $\mathfrak{c}$ in an affine Coxeter complex determines a unique hyperplane which is the wall of some reflection. We may in turn associate a root to this panel by considering the root containing $\mathfrak{c}$. Suppose $p$ is a special vertex and $\mathfrak{c}$ is a chamber containing $p$. Then each panel of $\mathfrak{c}$ containing $p$ determines a root, and the intersection of all these roots is called a {\it sector}. We say that the sector is {\it based at $p$} and has {\it germ $\mathfrak{c}$}.

We shall be particularly interested in the affine Coxeter groups $\tilde{A}_d$ which have the presentation $W = \langle r_i| r_i^2 = (r_i r_j)^{m_{i,j}} = 1,  i \in \mathbb{Z}/(d+1) \mathbb{Z} \rangle$ with $m_{i, i+1} = m_{i, i-1} = 3$, and all other $m_{i, j} = 2$. It arises from reflections across the $d+1$ hyperplanes determined by the faces of a certain $d$-simplex in $\mathbb{R}^d$. The associated spherical Coxeter group is $\mathfrak{S}_{d+1}$. In the associated Coxeter complex, all vertices are special. In case $d = 2$, the affine Coxeter group arises via reflections across an equilateral triangle.

\subsection{The Bruhat-Tits building} \label{sec_buildings} Material for this section may be found in \cite{brown}. 

A {\it building} is a polysimplicial complex $\Delta$ which can be expressed as the union of subcomplexes $\Sigma$, called {\it apartments}, satisfying:
\begin{enumerate}
    \item Each apartment is a Coxeter complex.
    \item For any two top-dimensional polysimplicies (called {\it chambers}) $\mathfrak{c}_1$ and $\mathfrak{c}_2$, there is an apartment $\Sigma$ containing both of them.
    \item If $\Sigma$ and $\Sigma'$ are two apartments containing $\mathfrak{c}_1$ and $\mathfrak{c}_2$, then there is an isomorphism $\Sigma \to \Sigma'$ fixing $\mathfrak{c}_1$ and $\mathfrak{c}_2$ pointwise.
\end{enumerate}
Each building has as associated Coxeter group $W$. 

Associated to $G = \textnormal{PGL}(d, F)$ there is a building known as the {Bruhat-Tits building}, and it has associated Coxeter group $\tilde{A}_{d-1}$ (\cite{bruhat_tits}). Let $\mathcal{L}$ denote the set of free $\mathcal{O}$-modules of rank $d$ inside of $F^d$. We define a quotient set $\mathcal{B}$ of $\mathcal{L}$ by declaring two elements $L_1, L_2 \in \mathcal{L}$ equivalent if for some $c \in F^\times$ we have $c L_1 = L_2$. We furthermore put a simplicial complex structure on $\mathcal{B}$ by defining the top-dimensional simplices, known as chambers, to be those collections of elements $\{[L_1], \dots, [L_d]\}$ such that for some ordering of these elements and choice of representative $L_i \in [L_i]$ we have $L_{i_1} \supsetneq L_{i_2} \supsetneq \dots \supsetneq L_{i_d} \supsetneq \varpi L_{i_1}$. We clearly have a bijection between the vertices of $\mathcal{B}$ and the cosets $G/K$. The group $G$ acts on $\mathcal{B}$ by simplicial automorphisms.

The maximal split torus $T$ consisting of diagonal matrices defines a subcomplex $\mathcal{X}$ called the {\it standard apartment}. This subcomplex is in fact a Coxeter complex corresponding to the Coxeter group $\tilde{A}_{d-1}$. The vertices of this apartment are in bijection with $T/(T \cap K) \simeq A$. We may in turn identify $A$ with 
\begin{gather*}
    \Lambda := \mathbb{Z}^d/\mathbb{Z} \cdot (1, \dots, 1) \simeq \{(x_1, \dots, x_d) : x_i \in \mathbb{R}, \ \sum x_i = 0, \text{ and } x_i - x_j \in \mathbb{Z} \text{ for all } i, j \}
\end{gather*}
This naturally sits as a lattice inside the vector space
\begin{gather*}
    \mathfrak{a} := \{(x_1, \dots, x_d): x_i \in \mathbb{R} \text{ and } \sum x_i = 0 \}.
\end{gather*}
Geometrically we may view the apartment as a tesselation of $\mathbb{R}^{d-1}$ by a certain $(d-1)$-simplex. By taking the image of $\mathcal{X}$ under elements of $G$, we obtain all of the other apartments in $\mathcal{B}$. 

By the building axioms, given any two chambers in $\mathcal{B}$ there is an apartment containing both of them. This allows us to define a $W$-valued distance $d_W(\cdot, \cdot)$ between chambers by simply considering the $W$-valued distance bewteen the two chambers in any apartment containing both of them. Furthermore, we define a metric (in the usual sense of metric) on $\mathcal{B}$ as follows: given two points in $\mathcal{B}$, we may find an apartment containing both of them which may in turn be identified with $\mathfrak{a}$. We may then consider the Euclidean geodesic in $\mathfrak{a}$ connecting the corresponding two points and define the distance between the points to be the length of this geodesic. With respect to this metric $\mathcal{B}$ is a CAT(0) space, and in particular there is a unique geodesic connecting any two points. We denote the metric $d(\cdot, \cdot)$ and normalize it so that the shortest distance between 0 and any other lattice point in $\Lambda$ is 1. 

Given any fixed chamber $\mathfrak{c}$ and a fixed apartment $\Sigma$ containing $\mathfrak{c}$, there is a unique simplicial retraction $\rho_{\mathfrak{c}, \Sigma}$ which maps $\mathcal{B}$ onto $\Sigma$ in such a way that for any chamber $\mathfrak{d} \in \mathcal{B}$ we have $d_W(\mathfrak{c}, \mathfrak{d}) = d_W(\mathfrak{c}, \rho_{\mathfrak{c}, \Sigma}(\mathfrak{d}))$. The retraction $\rho_{\mathfrak{c}, \Sigma}$ has the property that it does not increase distances between chambers.

Given a vertex in $\mathcal{B}$, its link is also a building known as the {\it Tits building} $\Pi$ which has associated Coxeter group $\mathfrak{S}_d$. If $L$ is a representative for a vertex, then the vertices in its link have representatives of the form $L'$ with $L \supsetneq L' \supsetneq \varpi L$. Such $L'$ in turn are in bijection with non-trivial subspaces of $\mathbb{F}_q^d$. The simplices in $\Pi$ arise from flags in $\mathbb{F}_q^d$. The apartments in $\Pi$ arise from choosing any basis and considering all flags which may be built from that basis. Given any two chambers in $\Pi$, there is an apartment containing both.

The subset $A^+ \subset A$ corresponds to the elements in $\Lambda$ whose entries are weakly decreasing; we denote this subset by $\Lambda^+$. This naturally sits inside the subset of elements of $\mathfrak{a}$ whose elements are weakly decreasing which we denote by $\mathfrak{a}^+$; we shall also refer to $\mathfrak{a}^+$ as the {\it (standard) Weyl chamber}. The associated subcomplex of $\mathcal{X}$ is called the {\it standard sector}. More generally a {\it sector} in $\mathcal{B}$ is any subcomplex arising as the image of the standard sector under an element of $G$.

We may use the Cartan decomposition to define a {\it Weyl chamber-valued metric} $d_{A^+}(\cdot, \cdot)$ between vertices of the Bruhat-Tits building. Given vertices $v_1, v_2$, choose representatives $x_1, x_2 \in G$ and define $d_{A^+}(v_1, v_2) = \lambda$ where $x_1^{-1} x_2 = k_1 \varpi^{\lambda} k_2$ in the Cartan decomposition (with $k_1, k_2 \in K$ and $\varpi^{\lambda} \in A^+$). It is invariant under the $G$-action on $G/K$.

There are other natural ways to coordinatize $\mathfrak{a}$. We shall refer to the coordinates already given as {\it $\mathfrak{a}$-coordinates}. We may instead shift all of the entries by the same amount so that the last entry is 0. We refer to these as {\it partition coordinates}. Lastly, we may choose an ordered basis for $\mathfrak{a}$ consisting of the cone generators for the Weyl chamber, ordered and written in partition coordinates as: $(1, 0, \dots, 0), (1, 1, 0, \dots, 0), \dots, (1, \dots, 1, 0)$. When we write an element, expressed in partition coordinates, as a linear combination of elements this ordered basis, we denote the coefficients as the {\it cone coordinates}. In case $d = 3$, which is the focus of this paper, this involves writing an element of the form $(\lambda_1, \lambda_2, 0)$ as $(\lambda_1 - \lambda_2) (1, 0, 0) + \lambda_2 (1, 1, 0)$ so that the cone coordinates are $(\lambda_1 - \lambda_2, \lambda_2)$. In general we shall use the notation $(r, s)$ to denote cone coordinates. There is a partial ordering $\preceq$ on $\mathfrak{a}$ defined in cone coordinates by $(r_1, s_1) \preceq (r_2, s_2)$ if and only if $r_1 \leq r_2$ and $s_1 \leq s_2$. 

\subsection{Parametrization of $\Omega^+_{\textnormal{temp}}$} \label{sec_parametrization_tempered}

We can identify $\mathfrak{a}^*$ with $\mathfrak{a}$ using the $\mathfrak{a}$-coordinates on $\mathfrak{a}$. Let $\Lambda^* \subset \mathfrak{a}^*$ denote the dual lattice of $\Lambda$. In the natural coordinates on $\mathfrak{a}^*$, we have $\Lambda^* = \{(a_1, \dots, a_d) \in \mathfrak{a}^*: a_i \in \mathbb{Z} \}$. This is simply the root lattice of type $A_{d-1}$ and we may take as a basis the collection of simple roots: $\{(1, -1, 0, \dots, 0), (0, 1, -1, 0, \dots, 0), \dots, (0, \dots, 0, 1, -1)\}$. We have a natural $\mathfrak{S}_d$-action on this space by permuting coordinates; this action preserves the lattice $\Lambda^*$. If we take the group generated by $\mathfrak{S}_d$ and translations by $\Lambda^*$, we obtain the affine Coxeter group of type $\tilde{A}_{d-1}$. We may take as a fundamental domain any chamber in the associated Coxeter complex. In particular we can take the chamber $\mathcal{C}$ defined as the convex hull of: $\{(0, \dots, 0), (\frac{1}{d}, \dots, \frac{1}{d}, -\frac{d-1}{d}), (\frac{2}{d}, \dots, \frac{2}{d}, -\frac{d-2}{d}, -\frac{d-2}{d}), \dots, (\frac{(d-1)}{d}, -\frac{1}{d}, \dots, - \frac{1}{d})\}$.

This chamber $\mathcal{C}$ provides a parametrization of $\Omega^+_{\textnormal{temp}}$. Recall that $\Omega^+_{\textnormal{temp}} = \{(\alpha_1, \dots, \alpha_d) : \alpha_1 \in S^1 \textnormal{ and } \prod_j \alpha_j = 1\}/\mathfrak{S}_d$. We can write $\alpha_j = q^{\frac{2 \pi i}{\ln(q)} s_j}$. Furthermore, we can choose the $s_j$'s so that their sum is 0; hence we can take $(s_1, \dots, s_d) \in \mathfrak{a}^* \simeq \mathfrak{a}$. The choice of $(s_1, \dots, s_d)$ is unique up to adding an element in $\Lambda^*$ and up to permuting the entries, i.e. up to the action of the affine Coxeter group. Hence we see that $\mathcal{C}$ gives us precisely one representative for each element in $\Omega^+_{\textnormal{temp}}$. Notice that under this parametrization the extremal points of $\mathcal{C}$ correspond to the points whose Satake parameters are $(1, \dots, 1), (\zeta_d, \dots, \zeta_d), (\zeta_d^2, \dots, \zeta_d^2), \dots, (\zeta_d^{d-1}, \dots, \zeta_d^{d-1})$.

In Section \ref{spectral_bound_section}, we parametrize $\Omega^+_{\textnormal{temp}}$ slightly differently. We instead take our fundamental domain as the convex hull of $\{(0, 0, 0), \frac{2 \pi i}{\ln(q)}(\frac{1}{3}, \frac{1}{3}, -\frac{2}{3}), \frac{2 \pi i}{\ln(q)} (\frac{2}{3}, -\frac{1}{3}, -\frac{1}{3}) \}$ inside of $i \mathfrak{a}^*$, and we denote it by $\mathcal{S}$. This is a fundamental domain for the quotient of $i \mathfrak{a}^*$ by $\frac{2 \pi i}{2 \ln(q)} \Lambda^* \subset i \mathfrak{a}^*$ and the $\mathfrak{S}_3$-action. This way we can directly identify $(s_1, s_2, s_3) \in \mathcal{S}$ with $(q^{s_1}, q^{s_2}, q^{s_3}) \in \Omega^+_{\textnormal{temp}}$.

\subsection{Coloring eigenfunctions} \label{sec_coloring_eigenfunctions}

The valuation of the determinant of an element in $G$ is unchanged under multiplication by $K$. We thus have a map $\textnormal{val}: G/K \to \mathbb{Z}/d \mathbb{Z}$ by taking the valuation mod $d$ of the determinant of any coset representative. Let $B_j \subset G/K$ denote the preimage of $j \in \mathbb{Z}/d \mathbb{Z}$ under this map. This gives a coloring of $G/K$ with $d$ colors. This coloring is preserved under the action of $\textnormal{PSL}(d, F) \cdot K$ on $G/K$.

Let $\zeta_d$ be a primitive $d$th root of unity. Let $f_j$ with $j \in \mathbb{Z}/d \mathbb{Z}$ denote the function on $G/K$ defined by
\begin{gather*}
    f_j (g K) = (\zeta_d^j)^{\textnormal{val}(\det(g))}.
\end{gather*}
Let $A_\ell := \mathds{1}_{K \varpi^\lambda K}$ where $\lambda = (1, \dots, 1, 0, \dots, 0)$ is the partition with $\ell$ ones; thus $A_\ell \in H(G, K)$ and hence acts on functions on $G/K$. If $h$ is a function on $G/K$, then $A_\ell.h$ evaluated at a point $x \in B_j$ corresponds to the sum over the the neighbors of $x$ that lie in $B_{j + \ell}$. Hence we immediately see that the $f_j$ are eigenfunctions of all of the $A_\ell$; in fact $A_\ell f_j = b_{d, \ell} \zeta_d^{\ell \cdot j} f_j$ where $b_{d, \ell}$ is the number of $\ell$-dimensional subspaces inside of $\mathbb{F}_q^d$.

Any finite-dimensional unitary representation of $G$ must in fact be trivial on $\textnormal{PSL}(d, F)$ (see \cite{gorodnik_nevo}, Section 3.8.2, p. 31). If furthermore it is class 1, then it must be trivial on $\textnormal{PSL}(d, F) \cdot K$. We would thus have a representation of $G/\textnormal{PSL}(d, F) \cdot K \simeq \mathbb{Z}/d \mathbb{Z}$. Hence the only such irreducible representations are the $d$ one-dimensional representations of $\mathbb{Z}/d \mathbb{Z}$. The span of each functions $f_j$, viewed as functions on $G$, gives us a copy of each of these one-dimensional representations inside of $C^\infty(G)$.

Suppose $\Gamma < \textnormal{PSL}(d, F) \cdot K$. Then $\Gamma$ preserves the coloring on $G/K$, and hence the coloring descends to the quotient. Thus all the functions $f_\ell$ also descend to the quotient, and we subsequently get $d$ eigenfunctions of the $H(G, K)$ action on $L^2(\Gamma \backslash G / K)$. We shall also let $f_\ell$ denote the function descended to the quotient (in the sequel it will be clear from context whether we are working on the building or the quotient).

Suppose $\psi$ is an eigenfunction of $H(G, K)$ acting on $L^2(\Gamma \backslash G / K)$. Suppose the eigenvalues of $A_1, \dots, A_{d-1}$ acting on $\psi$ are $(\mu_1, \dots, \mu_{d-1})$. We claim then that $f_k \cdot \psi$, i.e. the pointwise product of $f_\ell$ and $\psi$ viewed as functions on the vertices of the quotient, is also an eigenfunction of $H(G, K)$ with eigenvalues $(\zeta_d^{k} \mu_1, \zeta_d^{k \cdot 2} \mu_2, \dots, \zeta_d^{k \cdot (d-1)} \mu_{d-1})$:
\begin{gather*}
A_\ell [ f_k \cdot \psi] (x) = \zeta_d^{k (\textnormal{val}(x) + \ell)} \mu_\ell \psi(x) = \zeta_d^{k ( \textnormal{val}(x) + \ell)} \mu_\ell \frac{ [f_k \cdot \psi] (x)}{\zeta_d^{k \cdot \textnormal{val}(x)}} = \zeta_d^{k \cdot \ell} \mu_\ell [f_k \cdot \psi] (x). 
\end{gather*}

Suppose $\psi$ has Satake parameters $(\alpha_1, \dots, \alpha_d)$. Then the Satake parameters of $f_k \cdot \psi$ are $(\zeta_d^k \alpha_1, \dots, \zeta_d^k \alpha_d)$. This follows immediately from the fact that the eigenvalue of $A_\ell$ is $q^{\ell \cdot (d-\ell)/2} \sigma_\ell(\alpha_1, \dots, \alpha_{d-1})$ where $\sigma_\ell$ is the $\ell$th elementary symmetric polynomial (which in particular is a degree $\ell$ homogeneous polynomial); see e.g. Prop. 2.1 in \cite{lubotzky_samuels_vishne}. 

More generally, we have an action of $\mathbb{Z}/d \mathbb{Z}$ on $\Omega$, the space of Satake parameters, by sending $(\alpha_1, \dots, \alpha_d)$ to $(\zeta_d^j \alpha_1, \dots, \zeta_d^j \alpha_d)$ with $j \in \mathbb{Z}/d \mathbb{Z}$. This action preserves the spaces $\Omega^+$ and $\Omega^+_\textnormal{temp}$.

The action of $\mathbb{Z}/d \mathbb{Z}$ on $\Omega^+_{\textnormal{temp}}$ using the coordinatization in Section \ref{sec_parametrization_tempered} simply corresponds to repeated iteration of the unique affine map which sends the extremal point of $\mathcal{C}$ associated to $(\zeta_d^k, \dots, \zeta_d^k)$ to the extremal point of $\mathcal{C}$ associated to $(\zeta_d^{k+1}, \dots, \zeta_d^{k+1})$. There is a unique fixed point of this action, namely the centroid $(\frac{d-1}{2d}, \frac{d-3}{2d}, \dots, \frac{1-d}{2d}) \in \mathcal{C}$ which corresponds to the point $(\zeta_{2d}^{d-1}, \zeta_{2d}^{d-3}, \dots, \zeta_{2d}^{1-d}) \in \Omega^+_{\textnormal{temp}}$. We call this point $\mathfrak{b}$. The eigenvalues of $A_1, \dots, A_{d-1}$ acting on an eigenfunction which has Satake parameters equal to $\mathfrak{b}$ are all equal to 0. If $d$ is prime, then $\mathfrak{b}$ is the only point in $\Omega$ with non-trivial stabilizer. More generally for each non-trivial subgroup of $\mathbb{Z}/d \mathbb{Z}$, the locus of points stabilized by this subgroup has dimension at least one and codimension at least one.

More generally, suppose $t = [\Gamma : \Gamma \cap \textnormal{PSL}(d, F) \cdot K]$. Then $\Gamma$ preserves a coloring of $G/K$ corresponding to the $d/t$ cosets of $t \mathbb{Z}/d \mathbb{Z} < \mathbb{Z}/d \mathbb{Z}$. The $d/t$ functions $f_0, f_t, f_{2t}, \dots, f_{d - t}$ are all constant on each coset and hence all descend to the quotient by $\Gamma$ and are all eigenfunctions of $H(G, K)$. We call these the {\it coloring eigenfunctions}. We have a $\mathbb{Z}/(d/t) \mathbb{Z}$-action on the eigenfunctions of $H(G, K)$ on $L^2(\Gamma \backslash G / K)$ given by multiplying by these $f_{\ell \cdot t}$. 


\subsection{Benjamini-Schramm convergence} \label{sec_defn_BS}

Suppose $(\Gamma_n)$ is a sequence of torsion-free lattices in $G$. Let $\hat{Y_n} = \Gamma_n \backslash \mathcal{B}$ denote the simplicial complex obtained by quotienting the building by $\Gamma_n$, and let $Y_n = \Gamma_n \backslash G / K$ denote the vertices of $\hat{Y}_n$. Let $\hat{\tau}_n: \mathcal{B} \to \hat{Y}_n$ denote the obvious projection. Suppose $y \in Y_n$. Let $\tilde{y} \in \mathcal{B}$ be any lift of $y$ under $\hat{\tau}_{n}$. The {\it injectivity radius of $y$}, denoted $\textnormal{InjRad}_{Y_n}(y)$, is the supremum of all $r$ such that the ball of radius $r$ in $\mathcal{B}$ centered at $\tilde{y}$ (with respect to the metric $d(\cdot, \cdot)$) maps injectively under $\hat{\tau}_n$ to $\hat{Y}_n$. The {\it injectivity radius of $Y_n$}, denoted $\textnormal{InjRad}(Y_n)$, is the supremum over all $y \in Y_n$ of $\textnormal{InjRad}_{Y_n}(y)$. We say that $(Y_n)$ {\it Benjamini-Schramm converges to $G/K$} if, for every $R > 0$, we have
\begin{gather}
\frac{\card(\{y \in Y_n: \textnormal{InjRad}_{Y_n}(y) \leq R \})}{\card(Y_n)} \to 0 \label{BS_convergence}
\end{gather}
as $n \to \infty$. The assumption that $\card(Y_n) \to \infty$ implies that $(Y_n)$ Benjamini-Schramm converges to $G/K$ \cite{gelander_levit}.

%% file: sections/proof_outline_redo.tex
\sloppy
\subsection{Notation} \label{sec_notation}
We continue the notation from Section \ref{preliminaries}. However, we now assume that $d = 3$. Let $\pi: G \to G/K$ denote the obvious projection. Let $\Gamma < G$ denote a torsion-free lattice, and $(\Gamma_n)$ denote a sequence of torsion-free lattices whose covolume goes to infinity. Let $\rho^\Gamma$ be the unitary $G$-representation corresponding to the right action on $L^2(\Gamma \backslash G)$. Let $E \subset G$ be a set with positive, finite Haar measure. Suppose $f \in L^2(\Gamma \backslash G)$. We define $\rho^{\Gamma}_E$ to be the operator on $L^2(\Gamma \backslash G)$ such that
\begin{gather}
[\rho_E^\Gamma.f](\Gamma h) := \frac{1}{\vol(E)} \int_E f(\Gamma h g) dg. \label{eqn_rho_gamma}
\end{gather}

We can identify $G/K$ with the set of vertices of $\mathcal{B}$ induced by identifying some particular vertex in $\mathcal{B}$ with the coset $1K$. If $B \subset G/K$ is a set of vertices of $\mathcal{B}$, we let $\card(B)$ denote the cardinality of this set; this is equivalent to computing the Haar measure of $\pi^{-1}(B) \subset G$. Let $Y = \Gamma \backslash G /K$ (and $Y_n = \Gamma_n \backslash G /K$, resp.) denote the vertices of the  simplicial complex $\Gamma \backslash \mathcal{B}$ (and $\Gamma_n \backslash \mathcal{B}$, resp.). Let $a \in L^\infty(Y)$ (and $a_n \in L^\infty(Y_n)$, resp.) be a test function with $||a||_\infty \leq 1$ (and $||a_n||_\infty \leq 1$, resp.). Let $D$ (and $D_n$, resp.) be a fundamental domain for the action of $\Gamma$ (and $\Gamma_n$, resp.) on $G/K$. Notice that functions on $Y$ are the same as $(\Gamma, K)$-invariant functions on $G$.

Suppose $Q \subset \mathfrak{a}$ is a polytope. Let $Q_m$ denote the $m$th dilate of $Q$, and let $Q_m^{\Lambda} := Q_m \cap \Lambda$ (when $m = 1$, we omit the subscript). We may naturally view any $Q_m^{\Lambda}$ as a subset of $A$, and hence also of $G$, in which case it makes sense to consider $K Q_m^\Lambda K \subset G$. 

Given $Q$, we may define associated {\it polytopal balls}: given a vertex $v \in G/K$, we define the {\it $Q$-shaped ball at $v$} to be:
\begin{gather*}
B_Q(v) := \{w \in G/K: d_{A^+}(v, w) \in Q^\Lambda \}.
\end{gather*}
We may also define a {\it polytopal norm} on $\mathfrak{a}^+$ induced by $Q$ as follows:
\begin{gather*}
|\lambda|_Q := \inf \{m \in \mathbb{R}_{\geq 0}: \lambda \in Q_m\}.
\end{gather*}
We shall often be interested in the ceiling of the polytopal norm which we denote
\begin{gather*}
|\lambda|_Q^{\textnormal{ceil}} := \lceil |\lambda|_Q \rceil.
\end{gather*}

\begin{figure}[h]
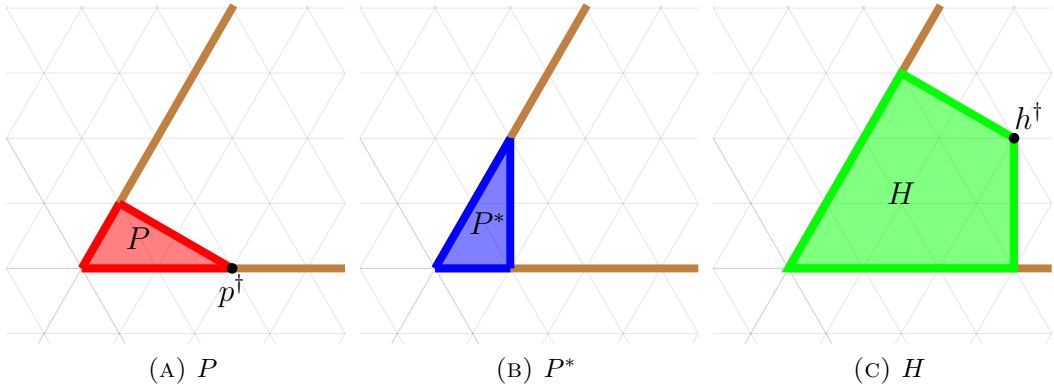

\centering
\subfloat[$P$]{\subfile{../images/simple_P_data}}
\subfloat[$P^*$]{\subfile{../images/Pstar_data}}
\subfloat[$H$]{\subfile{../images/H_data}}
\caption{Geometric realization of the polytopes $P$, $P^*$, and $H$.}
\label{table_polytopes}
\end{figure}

There are three polytopes which are of particular interest to us in this paper which we call $P$, $P^*$, and $H$; their defining inequalities in cone coordinates $(r, s)$ are given in Table \ref{table_polytopes}. The polytope $P$ has a distinguished vertex $p^\dag := (4/3, -2/3, -2/3)$.

\begin{table}[h]
\centering
\begin{tabular}{|c | c | c|}
\hline
Polytope $P$        &   Polytope $P^*$  &   Polytope $H$ \\
\hline
$r \geq 0$          &   $s \geq 0$      &   $r \geq 0$ \\
$s \geq 0$          &   $r \geq 0$      &   $s \geq 0$ \\
$r + 2s \leq 2$     &   $2r + s \leq 2$ &   $2r + s \leq 6$ \\ 
                    &                   &   $r + 2s \leq 6$ \\
\hline
\end{tabular}
\caption{Defining inequalities for $P$, $P^*$ and $H$.} \label{fig_polytope_data}
\end{table} 

We denote inner products on Hilbert spaces by $\langle \cdot, \cdot \rangle$, and we follow the convention that inner products are $\mathbb{C}$-linear in the second entry and sesquilinear in the first entry. If $U$ is an operator on a Hilbert space, we let $||U||_{\textnormal{HS}}$ denote the Hilbert-Schmidt norm of $U$, which may be defined as
\begin{gather*}
    ||U||_{\textnormal{HS}}^2 := \sum_{e_j} ||U e_j||^2,
\end{gather*}
where $\{e_j\}$ is any choice of orthonormal basis for the underlying Hilbert space. We use the notation $(\cdot, \cdot)$ for the pairing between $\mathfrak{a}$ and $\mathfrak{a}^*$. Let $\delta = (1, 0, -1)$. If we think of $\mathfrak{a}$ as the collection of trace zero diagonal matrices in $\mathfrak{sl}(3, \mathbb{R})$ and we identify $\mathfrak{a} \simeq \mathfrak{a}^*$ using the $\mathfrak{a}$-coordinates on $\mathfrak{a}$ (and the standard dot product in these coordinates), then $\delta$ corresponds to half the sum of positive roots.

Let $\mu$ denote the Plancherel measure normalized so that $\mu(\Omega^+_{\textnormal{temp}}) = 1$. We let $\Xi$ denote the exceptional locus, whose exact definition we defer until Section \ref{sec_tempered_pgl_3}.

Suppose $f$ and $g$ are two functions depending on several parameters. We use the notation $f \lesssim g$ to mean that the inequality is true up to a positive multiplicative constant not depending on any parameter. We use the notation $f \lesssim^b g$ to mean that when all parameters other than $b$ are fixed, the inequality is true for all $b$ sufficiently large up to a positive multiplicative constant not depending on $b$, but possibly depending on the other parameters. We write $f \lesssim^{b, c} g$ to mean that, assuming all other parameters other than $b$ and $c$ are fixed, then there exist $b_0$ and $c_0$ such that for all $(b, c)$ such that $b \geq b_0$ and $c \geq c_0$, the inequality is true up to a multiplicative constant not depending on $b$ or $c$ (but possibly depending on the other parameters).

\subsection{Ancillary propositions and lemmas}
In this subsection, we collect a sequence of Propositions and Lemmas which will allow us to prove Theorems \ref{main_thm} and \ref{main_thm2}. In cases where these results have straightforward or non-technical proofs, we include them here; otherwise we defer the proof to a later section.

\subsubsection{Reduction to the case of mean-zero test function} \label{sec_mean_zero}
By replacing a test function $a_n$ with the corresponding mean-zero test function $\bar{a}_n$, we may immediately simplify the summands in \eqref{main_thm_eqn} to $| \langle \psi_j^{(n)}, \bar{a}_n \psi_j^{(n)} \rangle |^2$. Thus we shall assume that $a_n$ is a mean-zero test function with $||a_n||_\infty \leq 1$.

\subsubsection{Polytopal ball averaging operators} \label{sec_polytopal_ball_sum}
Let $E_m$ denote the $P_m$-shaped ball centered at $1 K$:
\begin{gather*}
E_m := B_{P_m}(1K) = K P_m^\Lambda K.
\end{gather*}
We shall abuse notation and let
\begin{gather*}
    x E_m := B_{P_m}(x),
\end{gather*}
with $x \in G/K$. Notice that $E_m$ is a $K$-invariant set in $G/K$ and hence 
\begin{gather*}
\tilde{E}_m := \pi^{-1}(E_m)
\end{gather*}
is a $(K, K)$-invariant set in $G$. The function $\mathds{1}_{\tilde{E}_m}$ is in $H(G, K)$ and hence acts on functions on $G$ by right convolution. Furthermore it preserves $(\Gamma, K)$-invariant functions, so it acts on functions on both $G/K$ and on $Y$. We let $U_m$ denote the operator corresponding to right convolution with $\mathds{1}_{\tilde{E}_m}$. A straightforward calculation (see, e.g., Appendix B.7.1 of \cite{my_thesis}) shows that when $U_m$ is applied to a function on $G/K$ (or equivalently, a $(1, K)$-invariant function on $G$) we obtain:
\begin{gather*}
[U_m . f](g K) = \sum_{hK \in g \tilde{E}_m^{-1}} f(hK).
\end{gather*}
Also note that
\begin{gather*}
g \tilde{E}_m = \pi^{-1}(B_{P_m}(g K)).
\end{gather*}

The operator $U_m$ has kernel function $\mathcal{K}_m: G \times G \to \mathbb{C}$ defined as
\begin{gather}
\mathcal{K}_m(g, h) = \begin{cases} 
      1 & g^{-1} h \in \tilde{E}_m^{-1} \\
      0 & \textnormal{otherwise}.
   \end{cases} \label{kernel_u_m}
\end{gather}
Given a partition $\lambda = (\lambda_1, \lambda_2, 0)$, define $\lambda^\vee = (\lambda_1, \lambda_1 - \lambda_2, 0)$; in cone coordinates we have $(r, s)^\vee = (s, r)$. It follows from the fact that $(K \varpi^\lambda K)^{-1} = K \varpi^{\lambda^\vee} K$ that the polytope $P^*$ satisfies $\pi(\tilde{E}_m^{-1}) = B_{P^*_m}(1 K) = K(P^*_m)^\Lambda K$. Hence $U_m$ acting on $L^2(G/K)$ corresponds to summing up over the $P^*_m$-shaped ball centered at each vertex. From \eqref{kernel_u_m}, it is clear that the formula for kernel function of the adjoint $U_m^*$ is given by
\begin{gather}
\mathcal{K}^*_m(g, h) := \mathcal{K}_m (h, g) = \begin{cases}
    1 & g^{-1} h \in \tilde{E}_m \\
    0 & \textnormal{otherwise}.
    \end{cases} \label{kernel_u_m_star}
\end{gather}

Recall that $a$ denotes a test function on $Y$. We identify $a$ with the operator corresponding to multiplication by $a$. Suppose $M$ is a positive integer. Consider the following operator on $L^2(Y)$:
\begin{gather}
A_M := \frac{1}{M} \sum_{m = 1}^M \frac{1}{\card(E_m)} U_m^* \circ a \circ U_m. \label{A_M}
\end{gather}
When we use $a_n$ instead of $a$, we shall refer to the corresponding operator as $A_M^n$. Note that we may also consider $A_M$ as an operator on functions on $G$.

If $\psi_\nu$ is an eigenfunction of the $H(G, K)$-action on $L^2(Y)$ with spectral parameter $\nu \in \Omega^+$, then $U_m \psi_\nu = h_m(\nu) \psi_\nu$ for some complex number $h_m(\nu)$. We then get that
\begin{gather}
\Big | \langle \psi_\nu, A_M \psi_\nu \rangle \Big|^2 = \Big( \frac{1}{M} \sum_{m = 1}^M \frac{|h_m(\nu)|^2}{\card(E_m)} \Big)^2 \Big| \langle \psi_\nu, a \psi_\nu \rangle \Big|^2. \label{hms_rearrange}
\end{gather}

\subsubsection{Spectral bound}
\begin{proposition} [Spectral Bound; proof in Section \ref{pf_spectral_bound}] \label{spectral_bound}
Let $\Theta \subset \Omega^+_{\textnormal{temp}}$ be a subset of the tempered spectrum whose closure does not intersect the exceptional locus $\Xi$. Then for all $\nu \in \Theta$,
\begin{gather*}
\frac{1}{M} \sum_{m = 1}^M \frac{|h_m(\nu)|^2}{\card(E_m)} \gtrsim^M 1,
\end{gather*}
with the implied constant only depending on $\Theta$ and not on $\nu$ or $M$.
\end{proposition}

The proof of this proposition requires analysis of the spherical functions associated to each $\nu$. The strategy is as follows: we first rearrange the relevant expression so that it essentially becomes the sum of an exponential function on the lattice points in a polytope. We then apply Brion's formula, which provides an alternate expression for sums of exponential functions over lattice points in polytopes (discussed in Section \ref{sec_brion}); in particular this allows us to identify the dominating term. Finally we bound this dominating term asymptotically using a ``linear independence of characters'' argument. See Section \ref{spectral_bound_section}.

\subsubsection{Benjamini-Schramm convergence implies Plancherel convergence}
Recall the definition of $N(\Theta, Y_n)$ in \eqref{n_theta}. In Section \ref{bs_implies_plancherel_section}, we prove Theorem \ref{thm_bs_conv}. This in particular implies that for any $\Theta$ as in the statement of Theorem \ref{main_thm}, we have that Benjamini-Schramm convergence implies that
\begin{gather}
\lim_{n \to \infty} \frac{N(\Theta, Y_n)}{\card(Y_n)} = \mu(\Theta). \label{eqn_bs_plancherel}
\end{gather}

The work of Deitmar \cite{deitmar} shows a strong relationship between Benjamini-Schramm convergence and convergence to the Plancherel measure for locally compact topological groups. Similar results in the case of semisimple algebraic groups over local fields of characteristic zero appears in the work of Gelander-Levit \cite{gelander_levit}. These results involve first obtaining a form of convergence not for indicator functions of subsets of $\Omega^+$, but for a particular class of continuous functions on $\Omega^+$. By appealing to the Sauvageot density principle, one may upgrade convergence for this class of functions to the type of convergence claimed in \eqref{eqn_bs_plancherel} at least in characteristic zero. However, recently a gap in the proof of the Sauvageot density principle has been found as noted in \cite{nelson_venkatesh}. Thus in proving Theorem \ref{thm_bs_conv}, we first prove a weakened version of the Sauvageot density principle; this has the added benefit of allowing us to strengthen our main result to also include positive characteristic.

\subsubsection{Reduction to bounding the Hilbert-Schmidt norm}
The expression $|\langle \psi_j^{(n)}, a_n \psi_j^{(n)} \rangle|^2$ is often referred to as the {\it quantum variance} (of $a_n$ with respect to the wave function $\psi_j^{(n)}$).
\begin{proposition} [Quantum Variance Bounded by Hilbert-Schmidt Norm] \label{prop_reduction_to_HS}
Suppose $\Theta \subset \Omega^+_{\textnormal{temp}}$ is satisfies $\mu(\partial \Theta) = 0$, $\mu(\Theta) > 0$, and the closure of $\Theta$ does not intersect $\Xi$. Suppose $(Y_n)$ Benjamini-Schramm converge to $G/K$. Then
\begin{gather*}
\frac{1}{N(\Theta, Y_n)} \sum_{\psi_j^{(n)}: \nu_j^{(n)} \in \Theta} \bigg| \langle \psi_j^{(n)}, a_n \psi_j^{(n)} \rangle \bigg|^2 \lesssim^{n, M} \frac{1}{\textnormal{card}(Y_n)} ||A_M^n||^2_{\textnormal{HS}},
\end{gather*}
where the implied constant depends on $\Theta$, but does not depend on $n, M$ as long as both $n$ and $M$ are sufficiently large.
\end{proposition}
\begin{proof}
By \eqref{hms_rearrange}, Proposition \ref{spectral_bound}, Theorem \ref{thm_bs_conv}, and the inequality $|\langle \psi_j^{(n)}, A_M^n \psi_j^{(n)} \rangle|^2 \leq ||A_M^n \psi_j^{(n)}||^2$, we have
\begin{align*}
 \frac{1}{N(\Theta, Y_n)} \sum_{\psi_j^{(n)}: \nu_j^{(n)} \in \Theta} \big| \langle \psi_j^{(n)}, a_n \psi_j^{(n)} \rangle \big|^2 & \lesssim^{n, M} \frac{1}{\card(Y_n)} \sum_{\psi_j^{(n)}: \nu_j^{(n)} \in \Theta} \big| \langle \psi_j^{(n)}, A_M^n \psi_j^{(n)} \rangle \big|^2 \\
& \leq \frac{1}{\card(Y_n)} ||A_M^n||_{\textnormal{HS}}^2.
\end{align*}
\end{proof}

\subsubsection{Bounding the Hilbert-Schmidt norm using the kernel function on $G/K$}
Recall that when an operator has an associated kernel function, then its Hilbert-Schmidt norm may be computed by computing the $L^2$-norm of the kernel function. An operator acting on $L^2(Y)$ would then naturally have its associated kernel function as a function on $Y \times Y$. However, in many ways it's easier to work with functions on $G/K \times G/K$ because this space is homogeneous. If we have a function on $G/K \times G/K$ which is invariant under the diagonal $\Gamma$-action and is only supported near the diagonal, then by summing up over $\Gamma$, we may define an operator on $L^2(Y)$. The following lemma allows us to relate the Hilbert-Schmidt norm of this descended operator to the $L^2$-norm of the kernel function on $D \times G/K$ where $D$ is a fundamental domain for the $\Gamma$-action on $G/K$. We also obtain an error term in terms of number of points in $Y$ with small injectivity radius.

\begin{lemma}[Lifting the Kernel to $G/K$; proof in Section \ref{pf_kernel_lemma}] \label{lemma_kernel}
Let $\mathcal{K}:G/K \times G/K \to \mathbb{C}$ be a function which is invariant under the diagonal $\Gamma$-action. Suppose $R \geq 0$ is such that $\mathcal{K}(z, w) = 0$ whenever $d(z, w) \geq R$. Let $\bar{\mathcal{K}}^{\textnormal{Op}}_Y$ denote the operator on $L^2(Y)$ defined by this kernel. Then there exist $C > 0$, independent of $R$ and $\Gamma$, such that
\begin{align}
||\bar{\mathcal{K}}^{\textnormal{Op}}_Y||_\textnormal{HS}^2 & \leq \sum_{z \in D} \sum_{w \in G/K} |\mathcal{K}(z, w)|^2 \label{main_term} \\
& + C \frac{R^2 q^{4 R}}{\textnormal{InjRad}(Y)^{2}} ||\mathcal{K}||_\infty^2 \card(\{y \in Y:\textnormal{InjRad}_Y(y) \leq R\}). \label{error_term}
\end{align}
\end{lemma}

Note that the term of the form $R^2 q^{4 R} = (R q^{2 R})^2$ comes from the ``volume'' of a ball of radius $R$ in $\mathcal{B}$ (with respect to $d(\cdot, \cdot)$), and the term $\textnormal{InjRad}(Y)^{2}$ in the denominator comes from the area of a Euclidean ball of radius $\textnormal{InjRad}(Y)$ in $\mathbb{R}^2$. See Section \ref{kernel_section}.

\subsubsection{Explicit formula for the kernel function of $A_M$}
We ultimately wish to apply Lemma \ref{lemma_kernel} to the appropriate function $\mathcal{K}: G/K \times G/K \to \mathbb{C}$ so that $A_M = \bar{\mathcal{K}}_Y^{\text{Op}}$. 
\begin{proposition}[Formula for the Kernel Function] \label{prop_explicit_kernel}
Let $\mathcal{L}_M : G \times G \to \mathbb{C}$ be defined by:
\begin{gather}
\mathcal{L}_M(x, y) := \frac{1}{M} \sum_{m = 1}^M \frac{1}{\card(E_m)} \int_{ x \tilde{E}_m \cap y \tilde{E}_m} a(z) dz.
\end{gather}
Then $\mathcal{L}_M$ is invariant under the right $(K \times K)$-action and the left diagonal $\Gamma$-action and hence defines a function $\mathcal{L}_M': G/K \times G/K \to \mathbb{C}$ as in Lemma \ref{lemma_kernel}. The associated operator on $L^2(Y)$ is exactly $A_M$. 
\end{proposition}
\begin{proof}
Suppose $f$ is a function on $G$. Then by \eqref{kernel_u_m} and \eqref{kernel_u_m_star},
\begin{align}
[(U_m^* \circ a \circ U_m).f](x) &= \int_{x \tilde{E}_m} \int_{z \tilde{E}_m^{-1}} a(z) f(y) dy dz. \label{naive_kernel}
\end{align}
We now wish to change the order of integration. Suppose $y$ is fixed. We then must consider those $z$ such that $y \in z \tilde{E}_m^{-1}$ and $z \in x \tilde{E}_m$. But $y \in z \tilde{E}_m^{-1} \iff z \in y \tilde{E}_m$. Hence we need $z \in x \tilde{E}_m \cap y \tilde{E}_m$. Therefore we may rewrite (\ref{naive_kernel}) as
\begin{gather*}
\int_G \Big( \int_{x \tilde{E}_m \cap y \tilde{E}_m} a(z) dz \Big) \ f(y) dy,
\end{gather*}
from which it is clear that the kernel function for $U_m^* \circ a \circ U_m$ is exactly 
\begin{gather*}
\mathcal{H}_m(x, y) := \int_{x \tilde{E}_m \cap y \tilde{E}_m} a(z) dz. \label{h_m}
\end{gather*}
Thus by (\ref{A_M}), it is clear that the kernel function for $A_M$ is exactly
\begin{gather*}
\mathcal{L}_M (x, y) := \frac{1}{M} \sum_{m = 1}^M \frac{\mathcal{H}_m(x, y)}{\card(E_m)}.
\end{gather*}

If we replace $(x, y)$ by $(x k_1, y k_2)$ in \eqref{h_m}, with $k_1, k_2 \in K$, then we instead integrate over $x k_1 \tilde{E}_m \cap y k_2 \tilde{E}_m$, but because $\tilde{E}_m$ is $(K, 1)$-invariant, this is the same as $x \tilde{E}_m \cap y \tilde{E}_m$. Similarly, if we replace $(x, y)$ in \eqref{h_m} by $(\gamma x, \gamma y)$ with $\gamma \in \Gamma$, we now integrate $a(z)$ over $z$ in $\gamma x \tilde{E}_m \cap \gamma y \tilde{E}_m = \gamma (x \tilde{E}_m \cap y \tilde{E}_m)$, but this is the same as integrating $a(\gamma u) = a(u)$ over $u$ in $x \tilde{E}_m \cap y \tilde{E}_m$. Hence $\mathcal{L}_M$ is invariant under the right $(K \times K)$-action and the left diagonal $\Gamma$-action.
\end{proof}
\noindent If we replace $a$ with $a_n$ in the formula for $\mathcal{L}_M$, we instead write the kernel function as $\mathcal{L}_M^n$, and similarly we write $(\mathcal{L}_{M}^n)'$ for the associated function on $G/K \times G/K$. 

\subsubsection{Determining when two polytopal balls intersect}
In addition to $P$, another polytope that is of interest is the polytope $H$ which is defined in Table \ref{table_polytopes}. The following proposition explains its importance.
\begin{proposition}[When Polytopal Balls Intersect] \label{polytopal_ball_intersect}
Suppose $x, y \in G/K$. If the polytopal balls $x E_m$ and $y E_m$ intersect then $d_{A^+}(x, y) \in H_m^\Lambda$. Phrased another way, $x E_m \cap y E_m \neq \emptyset$ implies that $|d_{A^+}(x, y)|_H^{\textnormal{ceil}} \leq m$.
\end{proposition}

\begin{proof}
Proposition 4.4.4 of \cite{bruhat_tits} tells us that if we have three points $x, y, z$ such $d_{A^+}(x, y) = \lambda_1$, $d_{A^+}(y, z) = \lambda_2$, and $d_{A^+}(x, z) = \lambda_3$, then $\lambda_3$ lies in the convex hull of the Weyl group orbit of $\lambda_1 + \lambda_2$ inside of $\mathfrak{a}$. Notice that the polytope $P$ is the convex hull of the Weyl group orbit of $p^{\dag} = (4/3, -2/3, -2/3)$ and $P^*$ is the convex hull of the Weyl group orbit of $(p^*)^\dag := (2/3, 2/3, -4/3)$. Suppose $z \in x E_m \cap y E_m$. Then since $d_{A^+}(x, z) = \lambda_1 \in P_m$ and $d_{A^+}(z, y) = \lambda_2 \in (P^*)_m$, we get that $d_{A^+}(x, y) = \lambda_3$ is in the convex hull of the Weyl group orbit of $\lambda_1 + \lambda_2$ which in turn is in the convex hull of the Weyl group orbit of $p^\dag + (p^*)^\dag = (2, 0, -2)$. This is exactly the polytope $H$.
\end{proof}

The polytope $H_M$ is contained in the restriction to $\mathfrak{a}^+$ of the metric ball of radius $(2 \sqrt{3})M$ when we normalize so that adjacent vertices in the associated Coxeter complex are distance 1 apart.

\begin{corollary}[Kernel Function is Supported Near the Diagonal] \label{cor_kernel_support}
We have $\mathcal{L}_M(z, w) = 0$ if $d(zK, wK) > (2 \sqrt{3}) M$.
\end{corollary}

\subsubsection{Bounding the error term from lifting the kernel to $G/K$} 
We wish to use Lemma \ref{lemma_kernel} (Lifting the Kernel Function to $G/K$) and Proposition \ref{prop_explicit_kernel} (Explicit Formula for the Kernel Function) to reduce to analyzing a function on $G \times G$. However, there is an error term showing up in \eqref{error_term} which is handled by the following lemma.
\begin{lemma}[Bounding the Error Term from Lifting the Kernel Function] \label{lemma_error_term}
Suppose $(Y_n)$ Benjamini-Schramm converges to $G/K$. Then, for every positive integer $M$ and $\varepsilon > 0$, there exists an $n_{M, \varepsilon}$ such that for all $n \geq n_{M, \varepsilon}$, we have
\begin{gather}
C \frac{(2 \sqrt{3} M)^2 q^{4 (2 \sqrt{3})M}}{\textnormal{InjRad}(Y_n)^2} ||\mathcal{L}_M^n||_\infty^2 \frac{\card(\{y \in Y_n : \textnormal{InjRad}_{Y_n}(y) \leq (2 \sqrt{3})M \})}{\card(Y_n)} \leq \varepsilon. \label{error_term_bound}
\end{gather}
\end{lemma}
\begin{proof}
As mentioned in Remark \ref{fewer_hypotheses}, we know that there exists a positive lower bound on $\textnormal{InjRad}(Y_n)$ for all $n$ (i.e. uniform discreteness is automatic). Furthermore, by the explicit formula for $\mathcal{L}_M^n$ in Proposition \ref{prop_explicit_kernel}, it is clear that $||\mathcal{L}_M^n||_\infty \leq ||a_n||_\infty$ and we have $||a_n||_\infty \leq 1$ by assumption. Furthermore, by the definition of Benjamini-Schramm convergence \eqref{BS_convergence}, we know that for any fixed $M$,
\begin{gather*}
\lim_{n \to \infty} \frac{\card(\{y \in Y_n : \textnormal{InjRad}_{Y_n}(y) \leq (2 \sqrt{3}) M \})}{\card(Y_n)} = 0.
\end{gather*}
\noindent By combining these observations we conclude that such an $n_{M, \varepsilon}$ exists.
\end{proof}

\subsubsection{Changing variables in the integral of the kernel function} \label{sec_changing_variables}
We wish to rewrite the expression on the right hand side of \eqref{main_term}. Rather than sum over pairs $(z, w) \in D \times G/K$, we instead group elements based on their relative positions, namely $d_{A^+}(z, w)$. For any given point $z \in G/K$ and $\lambda \in A^+$, there are exactly $N_\lambda$ many $w$'s such that $d_{A^+}(z, w) = \lambda$; $N_\lambda$ is defined in \eqref{n_lambda}. The shape of any given $x E_m \cap y E_m$ with $x, y \in G/K$ only depends on $m$ and $\lambda = d_{A^+}(x, y)$. These intersections are exactly translations of the set $E_m^\lambda = E_m \cap \varpi^{\lambda} E_m$. Translating $E_m^\lambda$ to be based at all the different vertices of $Y$ and in all possible ``orientations'', and then integrating the test function $a$ over these translated sets amounts to convolving the indicator function of the lift of $E_m^\lambda$ to $G$ (namely $\tilde{E}_m^\lambda$ defined in \eqref{e_tilde_m_lambda}) with the lift of $a$ to a $K$-invariant function on $\Gamma \backslash G$ (recall the definition of $\rho^\Gamma_{E}$ in \eqref{eqn_rho_gamma}). Unraveling all of these observations ultimately results in the following proposition.

\begin{proposition} [Changing Variables in the Kernel Integral; proof in Section \ref{pf_changing_variables}] \label{prop_change_variables}
We have
\begin{gather*}
\sum_{z \in D} \sum_{w \in G/K} |\mathcal{L}'_M(z, w)|^2 =
 \frac{1}{M^2} \sum_{\lambda \in A^+} N_\lambda \int_{\Gamma \backslash G} \bigg| \sum_{m = 1}^M \frac{\card(E_m^{\lambda})}{\card(E_m)}[\rho^{\Gamma}_{\tilde{E}_m^{\lambda}}.a](\Gamma g) \bigg|^2 dg,
\end{gather*}
where 
\begin{align}
N_\lambda : & = \vol(K \varpi^\lambda K), \label{n_lambda} \\
E_m^\lambda : &= E_m \cap \varpi^\lambda E_m, \nonumber \\
\tilde{E}_m^\lambda : & = \{g \in G: d_{A^+}(1K, gK) \in P_m^\Lambda \textnormal{ and } d_{A^+}(\varpi^\lambda K, gK) \in P_m^\Lambda \} \label{e_tilde_m_lambda} \\
& = \pi^{-1} (E_m^\lambda). \nonumber
\end{align}
\end{proposition}

\subsubsection{Changing the order of integration} \label{sec_change_order_of_integration}
A straightforward application of the Minkowski integral inequality allows us to change the order of integration in the integral for the $L^2$-norm of the kernel function.

\begin{proposition}[Changing the Order of Integration in the Kernel Integral] \label{prop_minkowski}
We have 
\begin{gather}
\sum_{\lambda \in A^+} N_\lambda \int_{\Gamma \backslash G} \bigg| \sum_{m = 1}^M \frac{\card(E_m^{\lambda})}{\card(E_m)}[\rho^{\Gamma}_{\tilde{E}_m^{\lambda}}.a](\Gamma g) \bigg|^2 dg \leq \sum_{\lambda \in H_M^\Lambda} N_\lambda \Big( \sum_{m = |\lambda|_H^\ceil}^M \frac{\card(E_m^\lambda)}{\card(E_m)} ||\rho^{\Gamma}_{\tilde{E}_m^\lambda}.a||_{L^2(\Gamma \backslash G)} \Big)^2. \label{change_order}
\end{gather}
\end{proposition}

\begin{proof}
Because of Proposition \ref{polytopal_ball_intersect} (When Polytopal Balls Intersect) $\card(E_m^\lambda) \neq 0$ only if $m \geq |\lambda|_H$, so we may replace the sum on the right hand side of \eqref{change_order} with a sum starting from $m = |\lambda|_H^\ceil$. Furthermore we have $m \leq M$; thus in order to have $|\lambda|_H \leq M$, we must also have $\lambda \in H_M^\Lambda$. Therefore we can replace the sum over $A^+$ with simply the sum over $H_M^\Lambda$.

If we apply the Minkowski integral inequality to the right side of (\ref{change_order}) we obtain
\begin{align*}
\int_{\Gamma \backslash G} \Big| \sum_{m = |\lambda|_H^\ceil}^M \frac{\card(E_\lambda^m)}{\card(E_m)} [\rho_{\tilde{E}_m^\lambda}^\Gamma . a](\Gamma g) \Big|^2 dg & \leq \Bigg( \sum_{m = |\lambda|_H^\ceil}^M \frac{\card(E_m^\lambda)}{\card(E_m)} \Big| \int_{\Gamma \backslash G} \big|[\rho_{\tilde{E}_m^\lambda}^\Gamma. a] (\Gamma g) \big|^2 dg \Big|^{1/2} \Bigg)^2 \\
& = \Bigg( \sum_{m = |\lambda|_H^\ceil}^M \frac{\card(E_m^\lambda)}{\card(E_m)} \big| \big| \rho_{\tilde{E}_m^\lambda}^\Gamma . a \big| \big|_{L^2(\Gamma \backslash G)} \Bigg)^2.
\end{align*}
\end{proof}

\subsubsection{A Nevo-style ergodic theorem for $G$} \label{sec_short_nevo}
On the right hand side of \eqref{change_order} we have an expression containing $||\rho^\Gamma_{\tilde{E}_m^\lambda}.a||$. The following proposition shows that we can bound this just in terms of the volume of $\tilde{E}_m^\lambda$ and the $L^2$-norm of $a$. This is quite remarkable as many drastically different sets in $G$ have the same volume. A result of this form for semisimple Lie groups is due to Nevo \cite{nevo}. One of the ingredients in Nevo's proof is the Kunze-Stein phenomenon for semisimple Lie groups \cite{kunze_stein, cowling_ks}. The Kunze-Stein phenomenon was later shown by Veca \cite{veca} to also hold for simply connected simple algebraic groups over non-archimedean local fields (such as $\textnormal{SL}(d, F)$). In Section \ref{sec_nevo}, we explain how this in turn implies that $\textnormal{PGL}(d, F)$ also has the Kunze-Stein phenomenon. We may then retrace Nevo's proof \cite{nevo} to obtain the following proposition. See also \cite{gorodnik_nevo, gorodnik_nevo2}.

\begin{proposition}[Nevo-Style Ergodic Theorem; proof in Section \ref{sec_pf_nevo}] \label{nevo}
Suppose $a \in L^2(\Gamma \backslash G)$ is orthogonal to all finite-dimensional subrepresentations of $L^2(\Gamma \backslash G)$; in particular, if $a$ is $(1, K)$-invariant, this is equivalent to assuming that $a$ is orthogonal to all coloring eigenfunctions. Then, there exist constants $\theta > 0$ and $C > 0$, not depending on $a$ or $\Gamma$, such that for any $E \subset G$ with finite positive Haar measure, 
\begin{gather}
||\rho^\Gamma_E.a||_{L^2(\Gamma \backslash G)} \leq \frac{C}{\vol(E)^{\theta}} \ ||a||_{L^2(\Gamma \backslash G)}. \nonumber
\end{gather}
\end{proposition}

\begin{corollary}[Applying the Nevo-Style Ergodic Theorem] \label{cor_nevo}
There exists a $\theta > 0$ such that
\begin{gather}
\sum_{\lambda \in H_M^\Lambda} N_\lambda \Big( \sum_{m = |\lambda|_H^\ceil }^M \frac{\card(E_m^\lambda)}{\card(E_m)} ||\rho^{\Gamma}_{\tilde{E}_m^\lambda}.a||_{L^2(\Gamma \backslash G)} \Big)^2 \lesssim ||a||_2^2 \sum_{\lambda \in H_M^\Lambda} N_\lambda \Big( \sum_{m = |\lambda|_H^\ceil}^M \frac{\card(E_m^\lambda)^{1-\theta}}{\card(E_m)} \Big)^2, \nonumber
\end{gather}
where the implied constant does not depend on $\Gamma$ or $M$.
\end{corollary}

\begin{proof}
This follows immediately from Proposition \ref{nevo} (Nevo-Style Ergodic Theorem) using the fact that $\vol(\tilde{E}_m^\lambda) = \card(E_m^\lambda)$. 
\end{proof}

\subsubsection{Bounding the size of intersections of polytopal balls} \label{sec_pf_outline_geometric_bound}

Corollary \ref{cor_nevo} (Applying the Nevo-Style Ergodic Theorem) allows us to now focus on bounding the size of $N_\lambda$ and the cardinality of the sets $E_m$ and $\tilde{E}_m^\lambda$. There is an explicit formula for $N_\lambda$ due to MacDonald \cite{macdonald-symmetric} from which the following is an easy consequence.

\begin{proposition}[Upper Bound on $N_\lambda$; proof in Section \ref{sec_pf_n_lambda}] \label{prop_vol_n_lambda}
We have
\begin{gather}
N_\lambda  = \vol(K \varpi^\lambda K) \lesssim (q^2)^{( \delta, \lambda )}, \nonumber
\end{gather}
where the implied constant does not depend on $\lambda$.
\end{proposition}

Recall that the polytope $P$ has a distinguished vertex $p^\dag = (4/3, -2/3, -2/3)$ and that $\delta = (1, 0, -1)$ is half the sum of positive roots. Note that $( \delta, p^\dag ) > 0$. Using the explicit formula for $N_\lambda$ together with Brion's formula (Theorem \ref{thm_brion}), one may show the following:
\begin{proposition}[Lower Bound on $\card(E_m)$; proof in Section \ref{sec_pf_n_lambda}] \label{prop_vol_e_m}
We have
\begin{gather}
\card(E_m) \gtrsim (q^2)^{( \delta, m \cdot p^\dag )}, \nonumber
\end{gather}
where the implied constant does not depend on $m$.
\end{proposition}

The following proposition is the most difficult in the entire paper. It provides an upper bound on the size of $E_m^\lambda$. 

\begin{proposition}[Upper Bound on $\card(E_m^\lambda)$; proof in Section \ref{sec_pf_e_m_lambda_bound}] \label{prop_vol_e_lambda_m}
We have that
\begin{gather}
\card(E_m^\lambda) \lesssim (q^2)^{( \delta, m \cdot p^\dag - \frac{\lambda}{2} )}, \nonumber
\end{gather}
where the implied constant does not depend on $m$ or $\lambda$.
\end{proposition}

We may now combine the bounds on $N_\lambda$ (Proposition \ref{prop_vol_n_lambda}), the size of $E_m$ (Proposition \ref{prop_vol_e_m}), and the size of $E_m^\lambda$ (Proposition \ref{prop_vol_e_lambda_m}) to obtain the following.

\begin{corollary}[Combining Bounds on $N_\lambda$, $\card(E_m)$, and $\card(E_m^\lambda)$] \label{cor_reduction_to_brion}
We have
\begin{gather}
\sum_{ \lambda \in H_M^\Lambda} N_\lambda \Big( \sum_{m = |\lambda|_H^\ceil}^M \frac{\card(E_m^\lambda)^{1-\theta}}{\card(E_m)} \Big)^2 \lesssim \sum_{\lambda \in H_M^\Lambda} (q^2)^{\theta ( \delta, \lambda - 2 |\lambda|_H \cdot p^\dag )}, \label{eqn_pre_final_brion}
\end{gather}
where the implied constant does not depend on $M$.
\end{corollary}

\begin{proof}
By Propositions \ref{prop_vol_e_m} (Lower Bound on $\card(E_m)$) and \ref{prop_vol_e_lambda_m} (Upper Bound on $\card(E_m^\lambda)$), we get that
\begin{align*}
\frac{\card(E_m^\lambda)^{1- \theta}}{\card(E_m)} & \lesssim \frac{(q^2)^{( \delta, m \cdot p^\dag - \frac{\lambda}{2} ) (1-\theta)}}{(q^2)^{(\delta, m \cdot p^\dag)}} = (q^2)^{-\theta (\delta, m \cdot p^\dag)} (q^2)^{-(1-\theta)(\delta,\frac{\lambda}{2})}.
\end{align*}
So then
\begin{gather}
\Big( \sum_{m = |\lambda|_H^\ceil}^M \frac{\vol(\tilde{E}_m^\lambda)^{1-\theta}}{\card(E_m)} \Big)^2 \lesssim (q^2)^{-(1-\theta) (\delta, \lambda)} \Big( \sum_{m = |\lambda|_H^\ceil}^M (q^2)^{-\theta (\delta, m \cdot p^\dag)} \Big)^2. \label{take_out_of_sum}
\end{gather}

We have 
\begin{align}
\sum_{m = |\lambda|_H^\ceil}^M (q^2)^{-\theta (\delta, m \cdot p^\dag)} &= \frac{(q^2)^{- |\lambda|_H^\ceil \theta ( \delta, p^\dag )} - (q^2)^{- (M+1) \theta ( \delta, p^\dag )} }{1-(q^2)^{-\theta ( \delta, p^\dag )}} \nonumber \lesssim (q^2)^{- |\lambda|_H \theta ( \delta, p^\dag )}, \label{simplify_sum}
\end{align}
because $( \delta, p^\dag ) > 0$ and $|\lambda|_H^\ceil \geq |\lambda|_H$. 

Combining Proposition \ref{prop_vol_n_lambda} (Upper Bound on $N_\lambda$) with \eqref{take_out_of_sum} and \eqref{simplify_sum}, we obtain
\begin{align*}
\sum_{ \lambda \in H_M^\Lambda} N_\lambda \Big( \sum_{m = |\lambda|_H^\ceil}^M \frac{\card(E_m^\lambda)^{1-\theta}}{\card(E_m)} \Big)^2 & \lesssim \sum_{\lambda \in H_M^\Lambda} (q^2)^{( \delta, \lambda )} \cdot (q^2)^{-(1-\theta) ( \delta, \lambda )} \cdot (q^2)^{-2 |\lambda|_H \theta ( \delta, p^\dag )} \\
    &= \sum_{\lambda \in H_M^\Lambda} (q^2)^{\theta ( \delta, \lambda - 2 |\lambda|_H \cdot p^\dag )}.
\end{align*}
\end{proof}

\subsubsection{Bounding the sum over $H_M^\Lambda$}
We now wish to apply Brion's formula yet again. However, the expression on the right hand side of \eqref{eqn_pre_final_brion} is not exactly an exponential function whose power is a linear functional in the coordinates of $\Lambda$, as would be needed for Brion's formula. However, we may derive an explicit formula for $|\lambda|_H$. From this we observe that we may partition $\mathfrak{a}^+$ into two sub-polytopes such that on each one $|\lambda_H|$ is a linear functional. This allows us to split the sum on the right hand side of in \eqref{eqn_pre_final_brion} into two pieces and apply (degenerate) Brion's formula on each piece to arrive at the following proposition.

\begin{proposition}[Bounding the Sum over $H_M^\Lambda$; proof in Section \ref{sec_final_brion}] \label{prop_final_brion}
We have
\begin{gather}
\sum_{\lambda \in H_M^\Lambda} (q^2)^{\theta ( \delta, \lambda - 2 |\lambda|_H \cdot p^\dag )} \lesssim M, \nonumber
\end{gather}
where the implied constant does not depend on $M$. 
\end{proposition}

\subsection{Proof of Theorem \ref{main_thm}} \label{pf_main_thm}
\begin{proof}[Proof of Theorem \ref{main_thm} (Quantum Ergodicity in the BS Limit for $\textnormal{PGL}(3, F)$)]
We now combine all of the preceding propositions and lemmas. Recall from Theorem \ref{main_thm} that we must show that
\begin{gather}
\frac{1}{N(\Theta, Y_n)} \sum_{\psi_j^{(n)}:\nu_j^{(n)} \in \Theta} \Big| \langle \psi_j^{(n)}, a_n \psi_j^{(n)} \rangle \Big|^2 \to 0 \label{want_to_show}
\end{gather}
as $n \to \infty$, where $a_n$ is a mean-zero test function which is orthogonal to all coloring eigenfunctions and such that $||a_n||_\infty \leq 1$ (see Section \ref{sec_mean_zero}). By Proposition \ref{prop_reduction_to_HS} (Quantum Variance Bounded by Hilbert-Schmidt Norm), we know that to show \eqref{want_to_show} it suffices to show that along some sequence $(n, M_n)$ with $M_n$ eventually always larger than some specific $M_0$, we have
\begin{gather*}
\frac{1}{\card(Y_n)} ||A_{M_n}^n||_{\textnormal{HS}} \to 0
\end{gather*}
as $n \to \infty$.

By Lemma \ref{lemma_kernel} (Lifting the Kernel to $G/K$), Proposition \ref{prop_explicit_kernel} (Formula for the Kernel Function), and Corollary \ref{cor_kernel_support} (Kernel Function is Supported Near the Diagonal) it suffices to show that for some such sequence $(n, M_n)$ as above, we have both
\begin{gather}
\frac{1}{\card(Y_n)} \sum_{z \in D_n} \sum_{w \in G/K} |(\mathcal{L}_{M_n}^n)'(z, w)|^2 \to 0, \label{main_term_2} \\
C \frac{(2 \sqrt{3} M_n)^2 q^{4 (2 \sqrt{3}) M_n}}{\textnormal{InjRad}(Y_n)^2} ||(\mathcal{L}_{M_n}^n)'(z, w)||_\infty^2 \frac{\card(\{y \in Y_n: \textnormal{InjRad}_{Y_n}(y) \leq 2 \sqrt{3} M_n \})}{\card(Y_n)} \to 0. \label{error_term_2}
\end{gather}

We handle each piece separately. We first handle \eqref{main_term_2}. Combining Propositions \ref{prop_change_variables} (Changing Variables in the Kernel Integral) and Proposition \ref{prop_minkowski} (Changing the Order of Integration in the Kernel Integral) with Corollary \ref{cor_nevo} (Applying the Nevo-Style Ergodic Theorem), we get
\begin{align*}
\frac{1}{\card(Y_n)} \sum_{z \in D_n} \sum_{w \in G/K} & |(\mathcal{L}_{M_n}^n)'(z, w)|^2 \\
& = \frac{1}{\card(Y_n)}  \frac{1}{M^2} \sum_{\lambda \in A^+} N_\lambda \int_{\Gamma \backslash G} \bigg| \sum_{m = 1}^M \frac{\card(E_m^{\lambda})}{\card(E_m)}[\rho^{\Gamma}_{\tilde{E}_m^{\lambda}}.a](\Gamma g) \bigg|^2 dg \\
    & \leq \frac{1}{\card(Y_n)} \frac{1}{M^2} \sum_{\lambda \in H_M^\Lambda} N_\lambda \Big( \sum_{m = |\lambda|_H^\ceil }^M \frac{\card(E_m^\lambda)}{\card(E_m)} ||\rho^{\Gamma}_{\tilde{E}_m^\lambda}.a||_{L^2(\Gamma \backslash G)} \Big)^2 \\
    & \lesssim \frac{||a_n||_2^2}{\card(Y_n)} \frac{1}{M^2} \sum_{\lambda \in H_M^\Lambda} N_\lambda \Big( \sum_{m = |\lambda|_H^\ceil}^M \frac{\card(E_m^\lambda)^{1-\theta}}{\card(E_m)} \Big)^2.
\end{align*}

Notice that
\begin{gather}
\frac{||a_n||_2^2}{\card(Y_n)} \leq ||a_n||_\infty^2. \nonumber
\end{gather}
Hence by Corollary \ref{cor_reduction_to_brion} (Combining Bounds on $N_\lambda$, $\card(E_m)$, and $\card(E_m^\lambda)$) and Proposition \ref{prop_final_brion} (Bounding the Sum over $H_M^\Lambda$) we have
\begin{align*}
\frac{||a_n||_2^2}{\card(Y_n)} \frac{1}{M^2} \sum_{\lambda \in H_M^\Lambda} N_\lambda \Big( \sum_{m = |\lambda|_H^\ceil}^N \frac{\card(E_m^\lambda)^{1-\theta}}{\card(E_m)} \Big)^2 & \lesssim \frac{1}{M^2} \sum_{\lambda \in H_M^\Lambda} (q^2)^{\theta ( \delta, \lambda - 2 |\lambda|_H^\ceil \cdot p^\dag )} \\
    & \lesssim \frac{1}{M}.
\end{align*}
Therefore \eqref{main_term_2} is true as long as $M \to \infty$. 

Lastly we may bound \eqref{error_term_2} using Lemma \ref{lemma_error_term} (Bounding the Error Term from Lifting the Kernel Function): let $M^{(k)}$ be some sequence going to infinity. By Lemma \ref{lemma_error_term}, for every $k$ we can find an $n_{M^{(k)}}$ such that for every $n \geq n_{M^{(k)}}$, the expression in \eqref{error_term_bound} holds with $ \varepsilon = \frac{1}{k}$. Let $n_k$ be some increasing sequence such that $n_k \to \infty$ and $n_k \geq n_{M^{(k)}}$. Let $\lambda(n) = \sup \{k : n \geq n_k\}$. Notice then that $\lambda(n) \to \infty$ as $n \to \infty$. Let $M_n = M^{(\lambda(n))}$. Then, since $n \geq n_{\lambda(n)} \geq n_{M^{(\lambda(n))}}$, we have that
\begin{gather*}
\frac{C_1 q^{C_2 R}}{\textnormal{InjRad}(Y_n)^2} ||(\mathcal{L}_{M_n}^n)'(z, w)||_\infty^2 \frac{\vol(\{y \in Y_n: \textnormal{InjRad}_{Y_n}(y) \leq 2 \sqrt{3} M_n \})}{\vol(Y_n)} \leq \frac{1}{\lambda(n)}.
\end{gather*}
Hence for this choice of $M_n$, we have this quantity going to zero as $n \to \infty$. 
\end{proof}

\subsection{Coloring eigenfunctions and quantum unique ergodicity} \label{coloring_and_que}

Suppose $\Gamma < \textnormal{PSL}(d, F) \cdot K$ is a torsion-free lattice. Let $Y = \Gamma \backslash G / K$. Let $\psi_s$ be an eigenfunction whose Satake parameters are $s$, and suppose $s$ has trivial stabilizer under the $\mathbb{Z}/d \mathbb{Z}$-action described in Section \ref{sec_coloring_eigenfunctions}. Let $a \in L^\infty(Y)$ be a test function. We can write $a = a^\perp + \sum_{j \in \mathbb{Z}/d \mathbb{Z}} \beta_j f_j$ where $a^\perp$ is orthogonal to all of the coloring eigenfunctions $f_j$. We then have
\begin{align}
\langle \psi_s, a \psi_s \rangle - \frac{1}{\card(Y)} \sum_{x \in Y} a(x) &= \langle \psi_s, a^\perp \psi_s \rangle + \sum_{j \in \mathbb{Z}/d \mathbb{Z}} \beta_j \langle \psi_s, f_j \psi_s \rangle - \beta_0 \nonumber \\
    &= \langle \psi_s, a^\perp \psi_s \rangle + \sum_{j \in \mathbb{Z}/d \mathbb{Z}, j \neq 0} \beta_j \langle \psi_s, \psi_{j.s} \rangle \nonumber \\
    &= \langle \psi_s, a^\perp \psi_s \rangle. \label{coloring_cancellation}
\end{align}
This argument uses the fact that $s \neq j.s$ for $j \neq 0$. Since eigenfunctions with different spectral paramters are orthogonal, we get then that $\langle \psi_s, \psi_{j.s} \rangle = 0$. This shows that if we take $\Theta$ in the statement of Theorem \ref{main_thm} to avoid those points in $\Omega^+_{\textnormal{temp}}$ with non-trivial stabilizers under the $\mathbb{Z}/d \mathbb{Z}$-action, then the result remains true even if we do not assume that our test functions are orthogonal to the coloring eigenfunctions, as the expression which is to be bounded has the same value as the expression we would obtain if we projected our test function onto the orthogonal complement of the coloring eigenfunctions.

Let $B_k$ be the subset of vertices of $Y$ which have color $k \in \mathbb{Z}/d \mathbb{Z}$. Let $L_k$ denote the space of functions supported on $B_k$. The projection operator onto $L_k$ can be expressed as $\textnormal{proj}_{L_k} = \frac{1}{d} \sum_{j = 0}^{d-1} \zeta_d^{-k \cdot j} f_j$ where we view $f_j$ as a multiplication operator. 

Suppose now that $d$ is prime and that $\mathfrak{b}$ arises as a spectral parameter for the $H(G, K)$ action on $L^2(Y)$; then $\mathfrak{b}$ is the only point with non-trivial stabilizer. We can choose a basis of eigenfunctions of $H(G, K)$ for spectral parameters not equal to $\mathfrak{b}$ such that if $\psi$ is in the basis, then $f_j \psi$ is also in the basis. We can then group this set into subsets of size $d$ such that all of the elements in one subset differ by multiplication by some $f_j$. Suppose we choose one representative from each set. Then the image of this collection of eigenfunctions under $\textnormal{proj}_{L_k}$ must give a linearly independent set in $L_k$. Thus the dimension of the span of the projection to each different $\textnormal{proj}_{L_k}$ must be equal, and must be strictly less than $\frac{\card(Y)}{d}$ since we know that $\mathfrak{b}$ occurs as a spectral parameter.

Let $\tilde{L}_k$ be the span of these projections. Then the functions which are orthogonal to all of the $\tilde{L}_k$'s is precisely the eigenspace corresponding to spectral parameter $\mathfrak{b}$. Let $\tilde{L}_k^\perp$ denote the orthogonal complement to $\tilde{L}_k$ within $L_k$. This space has dimension at least one and consists of eigenfunctions with spectral parameter $\mathfrak{b}$. Furthermore these functions are entirely supported on $L_k$. The dimensions of $\tilde{L_k}^\perp$ are the same for all $k$; call this dimension $\ell$. 

Let $\{\phi_k^m\}$ with $1 \leq m \leq \ell$ be an orthonormal basis for $\tilde{L_k}^\perp$. As we vary over $k$ we obtain an orthonormal basis for the eigenspace associated to $\mathfrak{b}$. However, each of these functions is entirely supported on a single $B_j$. The $\phi_k^m$ are also eigenfunctions of the operator corresponding to multiplication by $f_j$ for all $f_j$. In particular, if we take $f_j$ and $j \neq 0$, we have
\begin{gather}
    \Big|\langle \phi_k^m, f_j \phi_k^m \rangle - \frac{1}{\card(Y)} \sum_{x \in Y} f_j(x) \Big|^2 = 1. \label{que_failure}
\end{gather}

We could formulate a naive version of the ``quantum unique ergodicity'' conjecture in the Benjamini-Schramm limit as follows: suppose $\{Y_n\}$ Benjamini-Schramm converges to the building. Let $\{\psi_j^{(n)}\}$ be an orthonormal basis of eigenfunctions with spectral parameters in $\Omega^+_{\textnormal{temp}}$ acting on $L^2(Y_n)$. Then is it true that
\begin{gather}
    \lim_{n \to \infty} \sup_{j} \sup_{a_n \in L^\infty(Y_n)} \frac{1}{||a_n||_\infty} \Big| \langle \psi_j^{(n)}, a_n \psi_j^{(n)} \rangle  - \frac{1}{\card(Y_n)} \sum_{x \in Y_n} a_n(x) \Big| = 0? \label{que}
\end{gather}
This is analogous to Question 1.3 in \cite{le_masson_sahlsten}. By taking a tower of covers, it is clear that this cannot be true as stated. We may strengthen the hypothesis to assume that we only consider ``new'' eigenfunctions in $L^2(Y_n)$ which are orthogonal to ``old'' eigenfunctions arising from any space covered by $Y_n$.  However, \eqref{que_failure} shows that if we can find a sequence of lattices in $\textnormal{PSL}(d, F) \cdot K$ for which $Y_n$ Benjamini-Schramm converges to the building such that every $Y_n$ has $\mathfrak{b}$ as a spectral paramter, and this spectral parameter arises from a ``new'' eigenfunction, then we have found a counterexample to \eqref{que}.

On the other hand, let $\{\chi_k^m\}$ be defined by
\begin{gather*}
    \chi_k^m := \frac{1}{\sqrt{d}} \Big( \sum_{\ell \in \mathbb{Z}/d \mathbb{Z}} \zeta_d^{k \cdot \ell} \phi_\ell^m \Big).
\end{gather*}
Then $\chi_k^m$ is an orthonormal basis for the eigenspace of $\mathfrak{b}$, but it has the property that the $L^2$-mass on each $B_j$ is equal. This in turn implies that $\langle \chi_k^m, f_j \chi_k^m \rangle = 0$ for all $k, m, j$. This allows us to repeat the argument in \eqref{coloring_cancellation}, namely we could get rid of the assumption that our test functions are orthogonal to the coloring eigenfuctions in the statement of Theorem \ref{main_thm}. Whereas the eigenfunctions $\{\phi_k^m\}$ show that the most naive formulation of quantum unique ergodicity is likely not true in this context, there is still hope that quantum unique ergodicity could be true if we take $\{\chi_k^m\}$ as part of our basis of eigenfunctions. This is somewhat analogous to the difference between taking an arbitrary orthonormal basis of eigenfunctions of the Laplacian on an arithmetic compact hyperbolic surface, versus taking an orthonormal basis of joint eigenfunctions of the Laplacian and the Hecke algebra (or at least one Hecke operator) as is usually done in arithmetic quantum unique ergodicity.

\subsection{The polytope $P$ vs. the polytope used in Brumley-Matz} \label{compare_brumley_matz}
We use the polytope $P$ to define our polytopal ball averaging operators. Brumley-Matz \cite{brumley_matz} define analogous polytopal ball averaging operators but using a polytope which, somewhat conincidentally, is essentially identical to the polytope that we call $H$. The properties that $P$ has which are used critically at various steps in the proof of Theorem \ref{main_thm} are:
\begin{enumerate}
\item The vertices of $P$ are lattice points in $\Lambda$. This is necessary to be able to apply Brion's formula such as is done in the Proposition \ref{spectral_bound} (Spectral Bound) and Proposition \ref{prop_vol_e_lambda_m} (Upper Bound on $\card(E_m^\lambda)$).
\item There is a unique vertex $p^\dag$ maximizing the dot product with $\delta$. This in particular implies that the volume of our polytopal balls $E_m$ grows like $(q^2)^{m ( \delta, p^\dag )}$ as opposed to some polynomial in $m$ times such an expression. This is used in the spectral bound and the geometric bound. This is in contrast to metric balls which do not have this property. See also \eqref{eqn_vol_e_m} and the discussion in Section 2.5 of Brumley-Matz \cite{brumley_matz}.
\item The vertex $p^\dag$ is on the boundary of $\mathfrak{a}^+$. This is the property that distinguishes our polytope from the one used by Brumley-Matz. This property is used critically in the proof of Proposition \ref{prop_vol_e_lambda_m} (Upper Bound on $\card(E_m^\lambda)$). If we instead use the polytope $H$ to define our polytopal balls $E_m$, then the analogue of Proposition \ref{prop_vol_e_lambda_m} would work out to be:
\begin{gather*}
    \card(E_m^\lambda) \lesssim \Upsilon(\lambda) (q^2)^{(\delta, m \cdot p^\dag - \frac{\lambda}{2})},
\end{gather*}
where $\Upsilon(\lambda)$ is a piecewise degree one polynomial in $\lambda$. In fact one may derive a lower bound of the same form, and hence this inequality is sharp. This is discussed in further detail in Chapter X.5 of \cite{my_thesis}. This additional polynomial factor $\Upsilon(\lambda)$ in fact prevents one from completing the proof. It is for this reason that we instead use the polytope $P$, which allows one to get rid of this polynomial factor. Partly for this reason, our method of proof for Proposition \ref{prop_vol_e_lambda_m} is very different than the method used to prove the analogous step in \cite{brumley_matz} (i.e. Proposition 5.8 therein).
\item The polytope $P$ is convex and invariant under the action of the Weyl group. This property was used crucially in the proof of Proposition \ref{polytopal_ball_intersect}.
\end{enumerate}

%% file: images/simple_P_data.tex
\begin{tikzpicture}[scale=1]
\draw[color = black, opacity = 0.1] (-1.000000, 0.000000) -- (3.500000, 0.000000);
\draw[color = black, opacity = 0.1] (-1.000000, 0.866025) -- (3.500000, 0.866025);
\draw[color = black, opacity = 0.1] (-1.000000, 1.732051) -- (3.500000, 1.732051);
\draw[color = black, opacity = 0.1] (-1.000000, 2.598076) -- (3.500000, 2.598076);
\draw[color = black, opacity = 0.1] (-1.000000, 3.464102) -- (3.500000, 3.464102);
\draw[color = black, opacity = 0.1] (-1.000000, 0.000000) -- (3.500000, 0.000000);
\draw[color = black, opacity = 0.1] (-1.000000, -0.866025) -- (3.500000, -0.866025);
\draw[color = black, opacity = 0.1] (-0.577350, -1.000000) -- (2.020726, 3.500000);
\draw[color = black, opacity = 0.1] (-1.000000, 0.000000) -- (1.020726, 3.500000);
\draw[color = black, opacity = 0.1] (-1.000000, 1.732051) -- (0.020726, 3.500000);
\draw[color = black, opacity = 0.1] (-1.000000, 3.464102) -- (-0.979274, 3.500000);
\draw[color = black, opacity = 0.1] (-0.577350, -1.000000) -- (2.020726, 3.500000);
\draw[color = black, opacity = 0.1] (0.422650, -1.000000) -- (3.020726, 3.500000);
\draw[color = black, opacity = 0.1] (1.422650, -1.000000) -- (3.500000, 2.598076);
\draw[color = black, opacity = 0.1] (2.422650, -1.000000) -- (3.500000, 0.866025);
\draw[color = black, opacity = 0.1] (3.422650, -1.000000) -- (3.500000, -0.866025);
\draw[color = black, opacity = 0.1] (0.577350, -1.000000) -- (-1.000000, 1.732051);
\draw[color = black, opacity = 0.1] (1.577350, -1.000000) -- (-1.000000, 3.464102);
\draw[color = black, opacity = 0.1] (2.577350, -1.000000) -- (-0.020726, 3.500000);
\draw[color = black, opacity = 0.1] (3.500000, -0.866025) -- (0.979274, 3.500000);
\draw[color = black, opacity = 0.1] (3.500000, 0.866025) -- (1.979274, 3.500000);
\draw[color = black, opacity = 0.1] (3.500000, 2.598076) -- (2.979274, 3.500000);
\draw[color = black, opacity = 0.1] (0.577350, -1.000000) -- (-1.000000, 1.732051);
\draw[color = black, opacity = 0.1] (-0.422650, -1.000000) -- (-1.000000, 0.000000);
\path[draw, color = brown, line width = 3] (0, 0) -- (3.500000, 0.000000);
\path[draw, color = brown, line width = 3] (0, 0) -- (2.020726, 3.500000);
\fill[color=red, opacity = 0.5] (0, 0) -- (2.000000, 0.000000) -- (0.500000, 0.866025) -- cycle;
\draw[color=red, line width=3] (0, 0) -- (2.000000, 0.000000);
\draw[color=red, line width=3] (0, 0) -- (0.500000, 0.866025);
\draw[color=red, line width = 3] (2.000000, 0.000000) -- (0.500000, 0.866025);

\fill[color = black] (2, 0) circle (2pt);
\node[yshift = -3mm] at (2, 0) {$p^\dag$};

\node at (0.75, 0.4) {$P$};
\end{tikzpicture}

%% file: images/Pstar_data.tex
\begin{tikzpicture}[scale=1]
\draw[color = black, opacity = 0.1] (-1.000000, 0.000000) -- (3.500000, 0.000000);
\draw[color = black, opacity = 0.1] (-1.000000, 0.866025) -- (3.500000, 0.866025);
\draw[color = black, opacity = 0.1] (-1.000000, 1.732051) -- (3.500000, 1.732051);
\draw[color = black, opacity = 0.1] (-1.000000, 2.598076) -- (3.500000, 2.598076);
\draw[color = black, opacity = 0.1] (-1.000000, 3.464102) -- (3.500000, 3.464102);
\draw[color = black, opacity = 0.1] (-1.000000, 0.000000) -- (3.500000, 0.000000);
\draw[color = black, opacity = 0.1] (-1.000000, -0.866025) -- (3.500000, -0.866025);
\draw[color = black, opacity = 0.1] (-0.577350, -1.000000) -- (2.020726, 3.500000);
\draw[color = black, opacity = 0.1] (-1.000000, 0.000000) -- (1.020726, 3.500000);
\draw[color = black, opacity = 0.1] (-1.000000, 1.732051) -- (0.020726, 3.500000);
\draw[color = black, opacity = 0.1] (-1.000000, 3.464102) -- (-0.979274, 3.500000);
\draw[color = black, opacity = 0.1] (-0.577350, -1.000000) -- (2.020726, 3.500000);
\draw[color = black, opacity = 0.1] (0.422650, -1.000000) -- (3.020726, 3.500000);
\draw[color = black, opacity = 0.1] (1.422650, -1.000000) -- (3.500000, 2.598076);
\draw[color = black, opacity = 0.1] (2.422650, -1.000000) -- (3.500000, 0.866025);
\draw[color = black, opacity = 0.1] (3.422650, -1.000000) -- (3.500000, -0.866025);
\draw[color = black, opacity = 0.1] (0.577350, -1.000000) -- (-1.000000, 1.732051);
\draw[color = black, opacity = 0.1] (1.577350, -1.000000) -- (-1.000000, 3.464102);
\draw[color = black, opacity = 0.1] (2.577350, -1.000000) -- (-0.020726, 3.500000);
\draw[color = black, opacity = 0.1] (3.500000, -0.866025) -- (0.979274, 3.500000);
\draw[color = black, opacity = 0.1] (3.500000, 0.866025) -- (1.979274, 3.500000);
\draw[color = black, opacity = 0.1] (3.500000, 2.598076) -- (2.979274, 3.500000);
\draw[color = black, opacity = 0.1] (0.577350, -1.000000) -- (-1.000000, 1.732051);
\draw[color = black, opacity = 0.1] (-0.422650, -1.000000) -- (-1.000000, 0.000000);
\path[draw, color = brown, line width = 3] (0, 0) -- (3.500000, 0.000000);
\path[draw, color = brown, line width = 3] (0, 0) -- (2.020726, 3.500000);
\fill[color=blue, opacity = 0.5] (0, 0) -- (1.000000, 1.732051) -- (1.000000, 0.000000) -- cycle;
\draw[color=blue, line width=3] (0, 0) -- (1.000000, 1.732051);
\draw[color=blue, line width=3] (0, 0) -- (1.000000, 0.000000);
\draw[color=blue, line width = 3] (1.000000, 1.732051) -- (1.000000, 0.000000);

\node at (0.7, 0.6) {$P^*$};
\end{tikzpicture}

%% file: images/H_data.tex
\begin{tikzpicture}[scale=1]
\draw[color = black, opacity = 0.1] (-1.000000, 0.000000) -- (3.500000, 0.000000);
\draw[color = black, opacity = 0.1] (-1.000000, 0.866025) -- (3.500000, 0.866025);
\draw[color = black, opacity = 0.1] (-1.000000, 1.732051) -- (3.500000, 1.732051);
\draw[color = black, opacity = 0.1] (-1.000000, 2.598076) -- (3.500000, 2.598076);
\draw[color = black, opacity = 0.1] (-1.000000, 3.464102) -- (3.500000, 3.464102);
\draw[color = black, opacity = 0.1] (-1.000000, 0.000000) -- (3.500000, 0.000000);
\draw[color = black, opacity = 0.1] (-1.000000, -0.866025) -- (3.500000, -0.866025);
\draw[color = black, opacity = 0.1] (-0.577350, -1.000000) -- (2.020726, 3.500000);
\draw[color = black, opacity = 0.1] (-1.000000, 0.000000) -- (1.020726, 3.500000);
\draw[color = black, opacity = 0.1] (-1.000000, 1.732051) -- (0.020726, 3.500000);
\draw[color = black, opacity = 0.1] (-1.000000, 3.464102) -- (-0.979274, 3.500000);
\draw[color = black, opacity = 0.1] (-0.577350, -1.000000) -- (2.020726, 3.500000);
\draw[color = black, opacity = 0.1] (0.422650, -1.000000) -- (3.020726, 3.500000);
\draw[color = black, opacity = 0.1] (1.422650, -1.000000) -- (3.500000, 2.598076);
\draw[color = black, opacity = 0.1] (2.422650, -1.000000) -- (3.500000, 0.866025);
\draw[color = black, opacity = 0.1] (3.422650, -1.000000) -- (3.500000, -0.866025);
\draw[color = black, opacity = 0.1] (0.577350, -1.000000) -- (-1.000000, 1.732051);
\draw[color = black, opacity = 0.1] (1.577350, -1.000000) -- (-1.000000, 3.464102);
\draw[color = black, opacity = 0.1] (2.577350, -1.000000) -- (-0.020726, 3.500000);
\draw[color = black, opacity = 0.1] (3.500000, -0.866025) -- (0.979274, 3.500000);
\draw[color = black, opacity = 0.1] (3.500000, 0.866025) -- (1.979274, 3.500000);
\draw[color = black, opacity = 0.1] (3.500000, 2.598076) -- (2.979274, 3.500000);
\draw[color = black, opacity = 0.1] (0.577350, -1.000000) -- (-1.000000, 1.732051);
\draw[color = black, opacity = 0.1] (-0.422650, -1.000000) -- (-1.000000, 0.000000);
\path[draw, color = brown, line width = 3] (0, 0) -- (3.500000, 0.000000);
\path[draw, color = brown, line width = 3] (0, 0) -- (2.020726, 3.500000);
\fill[color = green, opacity = 0.5] (0, 0) -- (3.000000, 0.000000) -- (3.000000, 1.732051) -- (1.500000, 2.598076) -- cycle;
\draw[color = green, line width = 3] (0, 0) -- (3.000000, 0.000000) -- (3.000000, 1.732051) -- (1.500000, 2.598076) -- cycle;

\fill[color = black] (3.000000, 1.732051) circle (2pt);

\node[yshift = 3mm, xshift = 2mm] at (3.000000, 1.732051) {$h^\dag$};

\node[] at (1.5, 1) {$H$};

\end{tikzpicture}

%% file: sections/polytopes_redo.tex
\subsection{Convex polytopes} \label{sec_polytope_type}
Suppose $\{\beta_i(x)\}_{i \in I}$ is a finite set of linear functionals on $\mathbb{R}^d$, and $\{b_i\}_{i \in I}$ is a finite set of real numbers; we call the elements $b_i$ {\it cutoffs}. The set of solutions to the system of inequalities $\{\beta_i(x) \leq b_i\}_{i \in I}$ is called a {\it (convex) polytope}, denoted $Q$. The interior of $Q$ is those points satisfying $\{\beta_i(x) < b_i\}_{i \in I}$. The set of all points in $Q$ along which some subset of the defining inequalities becomes equalities is called a {\it face} $\mathfrak{f}$. We can associate to $\mathfrak{f}$ the collection of functionals $\text{Func}(\mathfrak{f})$ for which the associated inequality becomes an equality.

Let $\mathfrak{f}$ be a face of $Q$. Let $\textnormal{Aff}(\mathfrak{f})$ denote the affine span of points in $\mathfrak{f}$. The {\it dimension} of $\mathfrak{f}$ is defined as the dimension of $\textnormal{Aff}(\mathfrak{f})$. The dimension of $Q$ is the dimension of the affine space spanned by all points in $Q$. We shall always assume that we are working with $d$-dimensional polytopes (i.e. the dimension is equal to the dimension of the underlying Euclidean space; all polytopes have representations for which this is the case up to an appropriate definition for isomorphism of polytopes). 
\begin{proposition} \label{prop_affine_span}
Let $B$ be the set of simultaneous solutions to $\beta_j(x) = b_j$ for all $\beta_j \in \textnormal{Func}(\mathfrak{f})$. Then $B = \textnormal{Aff}(\mathfrak{f})$.
\end{proposition}
\begin{proof}
First note that $\mathfrak{f} = B \cap Q$. Hence $\textnormal{Aff}(\mathfrak{f}) \subset B$. Suppose $B$ is bigger than $\textnormal{Aff}(\mathfrak{f})$. Let $x \in \mathfrak{f}$ and suppose we begin moving in some direction in $B$ which no longer keeps us in $\textnormal{Aff}(f)$. Then we must no longer be in $Q$ as well, and hence we must violate some inequality $\beta_k(x_1, \dots, x_n) \leq b_k$ which was valid at $x$. Consider all sets of the form $\mathfrak{f} \cap \{\beta_k(x_1, \dots, x_n) = b_k\}$ for $g_k \notin \textnormal{Func}(\mathfrak{f})$. The union of these must cover all of $\mathfrak{f}$ by the above reasoning. If none of these sets is all of $\mathfrak{f}$, then each one is codimension at least one inside $\textnormal{Aff}(\mathfrak{f})$. However, it is not possible to cover a $k$-dimensional polytope by a union of finitely many $(k-1)$-dimensional polytopes. If instead one of these sets is all of $\mathfrak{f}$, then it violates the fact that $g_k \notin \textnormal{Func}(\mathfrak{f})$. Hence in both cases we get a contradiction and can conclude that $B = \textnormal{Aff}(\mathfrak{f})$.
\end{proof}

By associating to each face $\mathfrak{f}$ the collection of functionals $\text{Func}(\mathfrak{f})$, we obtain a subset of the power set of the functionals. We call this the {\it type of $Q$}. We shall later on work with families of polytopes for which the functionals remain constant but the cutoffs change. In such a setting it makes sense to talk about when two polytopes in the family have the same type.

Suppose $\mathcal{Q}$ is a family of polytopes all of the same type. If $Q \in \mathcal{Q}$, and $\mathfrak{f}$ is a face of $Q$, then there is an analogous face $\mathfrak{f}'$ for each $Q' \in \mathcal{Q}$. Let $v$ be a vertex of $Q$ and $e$ an edge associated to $v$. 
\begin{proposition}
The direction of the ray based at $v$ in the direction of $e$ is completely determined by the type of $Q$.
\end{proposition}
\begin{proof}
By Proposition \ref{prop_affine_span}, we can pick out from $\textnormal{Func}(e)$ a subset of $d-1$ linearly independent functionals, call them $\beta_1, \dots, \beta_{d-1}$, such that $\textnormal{Aff}(e)$ is the set of solutions to $\beta_i(x) = b_i$ for $1 \leq i \leq d-1$. We then can realize $\text{Aff}(e)$ as the solutions to $M x = b$ for $M$ an $(d-1) \times d$ matrix of rank $d-1$. The solutions to this may be described as any particular solution plus vectors in the kernel of $M$. However we can take the same matrix $M$ for all polytopes of the same type.

All that is left to do now is orient this line to get a ray. Choose some $\beta_j \in \textnormal{Func}(v)$ which is linearly independent from the $\beta_1, \dots, \beta_{d-1}$ above. We orient the ray so it points in the direction so that motion in this direction decreases the value of the functional $\beta_j$.
\end{proof}

Suppose we translate $Q$ so that $v$ is moved to the origin. Let $\textnormal{Rays}_Q(v)$ be the collection of rays associated to $v$ now thought of as based at the origin. If we take the conical span of $\textnormal{Rays}_Q(v)$ we get a polyhedral cone which we shall denote by $\textnormal{Cone}_Q(v)$. The main purpose of the above discussion is the following conclusion: $\textnormal{Cone}_Q(v)$ only depends on the type of $Q$ and not on $Q$ itself.

A cone is called {\it simplicial} if it has exactly $d$ generators. A polytope is called simplicial if the cone at every vertex is simplicial. 

\subsection{Brion's formula} \label{brion_section}
Now suppose $\Lambda$ is a lattice inside $\mathbb{R}^d$. Suppose $Q$ is a polytope such that all vertices are in $\Lambda$; such polytopes are called {\it lattice polytopes}. Let $Q^\Lambda := Q \cap \Lambda$. Consider $C := \textnormal{Cone}_Q(v)$ for some vertex $v$. Then to each ray $r \in \textnormal{Rays}_Q(v)$, we can associate the smallest non-zero element in $\Lambda$ on the ray. We will call this the {\it coprime generator of the ray}. Let $\textnormal{Coprime}(C)$ denote the union of all coprime generators for all rays of $C$. If $C$ is simplicial, we let $\textnormal{PP}(C)$ be the parallelpiped formed by taking $[0, 1) \cdot r_1 + \dots [0, 1) \cdot r_d$ for coprime generators $r_k$.

Suppose $C$ is a polyhedral cone based at the origin such that all of its extremal rays intersect $\Lambda$. Let $C^*$ be the polar cone of $C$, that is $C^* = \{x \in (\mathbb{R}^n)^*: (x, y) \leq 0 \textnormal{ for all } y \in C\}$. Let $\alpha$ be a (possibly $\mathbb{C}$-valued) functional. Define
\begin{equation}
\varsigma(C; \alpha) := \sum_{\gamma \in C \cap \Lambda} q^{( \alpha, \gamma )}. \nonumber
\end{equation}
This quantity converges if $\textnormal{Re}(\alpha)$ is in the interior of $C^*$. 

Suppose now that $C$ is simplicial. Then using power series it is easy to see that 
\begin{equation}
\varsigma(C; \alpha) = \Bigg(\sum_{\beta \in \textnormal{PP}(C) \cap \Lambda} q^{( \alpha, \beta )}\Bigg) \Bigg(\prod_{r \in \textnormal{Coprime}(C)} \frac{1}{1 - q^{( \alpha, r )}} \Bigg). \label{eqn_cone_generator}
\end{equation}
If $C$ is not simplicial, then there exists a partition of $C$ into simplicial cones such that each simplicial cone has its ray generators among the extremal rays of $C$. This implies that in all cases
\begin{equation}
\varsigma(C; \alpha) = \frac{R(\alpha)}{\prod_{r \in \textnormal{Coprime}(C)} 1 - q^{( \alpha, r )}}, \nonumber
\end{equation}
where $R(\alpha)$ is a Laurent polynomial in $q^{\alpha_i}$ (where $\alpha = (\alpha_1, \dots, \alpha_n)$). Hence the whole expression is a rational function in these variables. 

\begin{theorem}[Brion's formula \cite{brion}; see also \cite{barvinok} Theorem 4.5] \label{thm_brion}
Suppose $\Lambda$ is a lattice and $Q$ is a polytope with vertices $\mathcal{V}$, and $\mathcal{V} \subset \Lambda$. Let $\alpha$ be a functional for which no ray associated to any vertex is in its kernel (i.e. $\alpha$ is not orthogonal to any face of $Q$). Then,
\begin{equation}
\sum_{\gamma \in Q^\Lambda} q^{( \alpha, \gamma )} = \sum_{v \in \mathcal{V}} \varsigma(\textnormal{Cone}_Q(v); \alpha) \cdot q^{( \alpha, v )}. \label{brion_formula}
\end{equation}
\end{theorem}

Notice that each term in the sum is a product of two terms. The first term depends only on the structure of the cone at each vertex (and on $\Lambda$ and $\alpha$) and the second term is purely exponential. In the one-dimensional case, Brion's formula reduces to the geometric series formula. 

One important takeaway of Brion's formula is that the left hand side of (\ref{brion_formula}) is dominated by its largest term. More specifically, if we fix $Q$ and then dilate it by integer multiples, then the cone and coprime ray generators at each vertex stay the same; all such polytopes have the same type. The nondegeneracy assumption on $\alpha$ implies that there is a unique vertex $v^\dag$ maximizing the dot product with $\alpha$ (if we assume that $\alpha$ is an $\mathbb{R}$-valued functional). Then the dominating term on the left hand side of (\ref{brion_formula}) is $q^{(\alpha, n \cdot v^\dag)}$. On the right hand side of (\ref{brion_formula}) we always get something of the form $\sum_{v} C(v) q^{(\alpha, n \cdot v)}$ where $C(v)$ does not depend on $n$. Hence as $n \to \infty$, the right hand side of (\ref{brion_formula}) grows like $q^{(\alpha, n \cdot v^\dag)}$ as well.

\subsection{Degenerate Brion's formula}
We shall ultimately also be interested in the case when $\alpha$ is orthogonal to some face of $Q$. We call this the {\it degenerate case}. We would like a formula which allows us to understand the growth rate of the right hand side of \eqref{brion_formula} as some parameter in a family of polytopes goes to infinity. In such a case we will have potentially an entire face $\mathfrak{f}$ of $Q$ such that on $\mathfrak{f} \cap \Lambda$ the dot product with $\alpha$ is maximized. Call this maximal dot product value $M$. 
Then heuristically we expect that the sum of our exponential function over the lattice points in $Q$ is roughly of size $\vol(\mathfrak{f}) \cdot q^M$, where of course we are taking the $\textnormal{dim}(\mathfrak{f})$-dimensional volume of $\mathfrak{f}$. Then would then imply that as we dilate $Q$ by a factor of $n$, we should expect the growth of the sum of our exponential function over lattice points to be of the form $n^{\text{dim}(\mathfrak{f})} \cdot q^{n \cdot M}$.

An illustrative example is the polytope $Q = [0, 1] \times [0, 1]$, the lattice $\Lambda = \mathbb{Z}^2$, and the functional $\alpha = (0, a)$ for some value of $a > 0$. Then $Q$ can be defined by $\{x_1 \geq 0; x_2 \geq 0; x_1 \leq 1; x_2 \leq 1\}$. The dot product with $\alpha$ is maximized at both $(0, 1)$ and $(1, 1)$. Let $Q_{m, n}$ be the polytope defined by $\{x_1 \geq 0; x_2 \geq 0; x_1 \leq m; x_2 \leq n\}$. All $Q_{m, n}$ have the same type if $m, n > 0$. It is easy to calculate that
\begin{equation}
\sum_{\gamma \in Q_{m, n}^\Lambda} q^{( \alpha, \gamma )} = (m+1) \cdot \frac{1-q^{(n+1)\cdot a}}{1-q^a}. \nonumber
\end{equation}
This grows like $n \cdot q^{n \cdot a}$ if we take $m = n \to \infty$, matching our above heuristic: the dimension of the face in $Q_{m,n}$ along which the dot product with $\alpha$ is maximized is 1, and the value of that dot product is $n \cdot a$.

Suppose $\alpha$ is a functional. For each $v \in Q$ define the {\it degeneracy of $v$ with respect to $\alpha$} as the dimension of the span of the rays which are orthogonal to $\alpha$. We say that a vertex is {\it good} if its degeneracy is zero; otherwise it is {\it bad}. We say that a ray at a bad vertex is bad if it is orthogonal to $\alpha$.

We shall now derive a degenerate form of Brion's formula. It is possible that some form of this already exists in the literature but we have been unable to find it. We also feel that this result may be of interest in its own right.

\begin{theorem} \label{prop_degenerate_brion}
Suppose $\Lambda$ is a lattice and $\mathcal{Q}(b_1, \dots, b_\ell)$ is a family of lattice polytopes of the same type defined by $\beta_i(x) \leq b_i$. Let $\alpha$ be some functional. For all $Q = \mathcal{Q}(b_1, \dots, b_\ell)$ in this family
\begin{equation}
\sum_{\gamma \in Q^\Lambda} q^{( \alpha, \gamma )} = \sum_{v \in \mathcal{V}} R_v(b_1, \dots, b_\ell) q^{( \alpha, v )}, \label{eqn_degenerate_brion}
\end{equation}
where each $R_v$ is a polynomial in the $b_j$ whose degree is bounded by greatest degeneracy of any vertex in $V$ (which is at least the largest dimension of any face orthogonal to $\alpha$) and which is degree 0 if $v$ has degeneracy 0.
\end{theorem}

\begin{proof}
Let $\tau$ be any functional which is not orthogonal to any edge in $Q$. Theorem \ref{thm_brion} allows us to express for any $\varepsilon$ small enough
\begin{equation}
\sum_{\gamma \in Q^\Lambda} q^{( \alpha + \varepsilon \cdot \tau, \gamma )} = \sum_{v} \varsigma(\textnormal{Cone}_Q(v); \alpha + \varepsilon \cdot \tau) q^{( \alpha + \varepsilon \cdot \tau, v )}. \label{degen_brion_eqn}
\end{equation}

We shall further separate out the sum on the left hand side into the sum over good and bad vertices. The terms coming from the good vertices define holomorphic functions in $\varepsilon$ near $\varepsilon = 0$. When $\varepsilon = 0$, the term corresponding to a good vertex $v$ is exactly
\[
\varsigma(\textnormal{Cone}_Q(v); \alpha) q^{( \alpha, v )}.
\]
The coefficient of $q^{( \alpha, v )}$ only depends on the type of $Q$, and hence is a non-zero constant when viewed as a function of the $b_j$'s.

The left hand sum in \eqref{degen_brion_eqn} is also clearly a holomorphic function in $\varepsilon$ as it is a finite sum of exponential functions. Consequently, the sum over the bad vertices also has a holomorphic extension to $\varepsilon = 0$.

Recall that $\varsigma(\textnormal{Cone}_Q(v); \alpha + \varepsilon \cdot \tau)$ corresponds to the (analytic continuation) of the sum of $q^{(\alpha + \varepsilon \cdot \tau, \lambda)}$ over $\lambda \in \Lambda \cap \textnormal{Cone}_Q(v)$. We may partition this cone into simplicial cones, then partition the intersection of two of these cones into simplicial cones (of smaller dimension), and so on. In this way, by using inclusion-exclusion, we can express the sum of $q^{(\alpha + \varepsilon \cdot \tau, \lambda)}$ over the lattice points of $\textnormal{Cone}_Q(v)$ as a weighted sum of $q^{(\alpha + \varepsilon \cdot \tau, \lambda)}$ over lattice points in simplicial subcones, each of whose rays are contained among the rays of $\textnormal{Cone}_Q(v)$. 

Suppose now $v$ is a bad vertex and we have partitioned $\textnormal{Cone}_Q(v)$ into simplicial cones $C$ as described above. We then have
\begin{gather*}
    \varsigma(\textnormal{Cone}_Q(v); \alpha + \varepsilon \cdot \tau) = \sum_{(v, C) \textnormal{ with $C$ simplicial}} B_{(C, v)} \cdot \frac{\textnormal{PP}(C; \alpha + \varepsilon \cdot \tau)}{\prod\limits_{\textnormal{rays $r$ of C}} 1 - q^{(\alpha + \varepsilon \cdot \tau, r)}} 
\end{gather*}
for some constants $B_{(C, v)}$. 

 We now analyze the following sum:
\begin{gather}
\sum_{(v, C) \textnormal{ s.t. $v$ bad, $C$ simplicial}} B_{(C, v)} \cdot \frac{\textnormal{PP}(C; \alpha + \varepsilon \cdot \tau)}{\prod\limits_{\textnormal{rays $r$ of C}} 1 - q^{( \alpha + \varepsilon \cdot \tau, r )}} q^{( \alpha + \varepsilon \cdot \tau, v )}. \label{bad_cone_vertex}
\end{gather}
This is simply a rewriting of the right hand side of \eqref{degen_brion_eqn} with the terms correspond to good vertices removed. 

Let $v^{\textnormal{bad}}$ be the vertex whose associated cone $C^{\textnormal{bad}}$ contains the greatest number of bad rays. Let $\textnormal{BR}$ denote these corresponding bad rays. Let $k$ be the cardinality of $\textnormal{BR}$. Notice then that $k$ is at most the maximal degeneracy of any vertex. We shall multiply (\ref{bad_cone_vertex}) by a ``fancy one'' given by the expression
\begin{gather}
\frac{\prod\limits_{t \in \textnormal{BR}} \big(1 - q^{( \alpha + \varepsilon \cdot \tau, t )} \big)}{\prod\limits_{t \in \textnormal{BR}} \big(1 - q^{( \alpha + \varepsilon \cdot \tau, t )} \big)}, \nonumber
\end{gather}
which has a removable singularity at $\varepsilon = 0$. Let 
\begin{gather}
\nu(v, C; \alpha, \varepsilon, \tau) := B_{(C, v)} \cdot \textnormal{PP}(C; \alpha + \varepsilon \cdot \tau) q^{( \alpha + \varepsilon \cdot \tau, v )}. \nonumber
\end{gather}
We thus obtain
\begin{gather}
\frac{\sum\limits_{(v, C) \textnormal{ s.t. $v$ bad}} \Big( \nu(v, C; \alpha, \varepsilon, \tau)  \prod\limits_{t \in \textnormal{BR}} \big( 1 - q^{( \alpha + \varepsilon \cdot \tau, t )} \big) \prod\limits_{r \in \textnormal{Rays}(C)} \big(1 - q^{( \alpha + \varepsilon \cdot \tau, r )} \big)^{-1} \Big)}{\prod\limits_{t \in \textnormal{BR}} \big(1 - q^{( \alpha + \varepsilon \cdot \tau, t )}\big)}. \label{pre_lhopital}
\end{gather}

We wish to take a limit as $\varepsilon \to 0$. We shall pair every term in the product over $\textnormal{Rays}(C)$ in the numerator of (\ref{pre_lhopital}) coming from a bad ray with some term in the product over $\textnormal{BR}$; we can do this because the number of bad rays at each $v$ is at most $k$. These give expressions of the form $\frac{1 - q^{a \varepsilon}}{1 - q^{b \varepsilon}}$ with $a$ and $b$ non-zero (each term has a different $a$ and $b$), and these functions extend holomorphically to $\varepsilon = 0$ (it is a removable singularity). All the remaining terms in the numerator are also holomorphic at $\varepsilon = 0$. 

We now wish to apply L'Hopital's rule. We know that the limit at $\varepsilon \to 0$ must exist because, 
 as discussed above, this function must extend holomorphically to $\varepsilon = 0$. The denominator vanishes to order $k$. Hence the numerator must also vanish to order $k$. We now apply L'Hopital $k$ times.

When we differentiate the denominator $k$ times and set $\varepsilon = 0$ we get
\begin{gather}
\prod_{t \in \textnormal{BR}} \Big(-( \tau, t ) \Big) \ln(q)^k. \nonumber
\end{gather}
This expression only depends on the type.

Now let's rewrite the numerator in the form 
\begin{gather}
\sum_{(v, C) \textnormal{ s.t. $v$ bad}} F(C; \alpha, \varepsilon, \tau) q^{( \alpha + \varepsilon \cdot \tau, v ) }. \nonumber
\end{gather}
The functions $F$ depends only on $\alpha$, $\varepsilon$, $\tau$, and the cone at $v$, but not on the specific coordinates of $v$ (that is, for a family of polytopes with the same fixed type, $F$ stays the same). The $k$th derivative of this expression is
\begin{gather}
\sum_{ \ell = 0}^k \frac{\partial^{k - \ell}F(C; \alpha, \varepsilon, \tau)}{\partial \varepsilon^{k - \ell}} q^{( \alpha + \varepsilon \cdot \tau, v )} ( \tau, v )^{\ell} \ln(q)^\ell \nonumber
\end{gather}
and when we set $\varepsilon = 0$, we get
\begin{gather}
\Bigg(\sum_{\ell = 0}^k \frac{\partial^{k - \ell}F(C; \alpha, \varepsilon, \tau)}{\partial \varepsilon^{k - \ell}}\Bigg\rvert_{\varepsilon = 0} ( \tau, v )^\ell \ln(q)^\ell \Bigg) q^{( \alpha, v )}. \nonumber
\end{gather}
The coordinates of $v$ are in turn linear functionals in the $b_i$. Hence these terms are degree (at most) $k$ polynomials in the $b_i$.
\end{proof}

\begin{remark}
Notice that in case our functional $\alpha$ is identically zero, then the left hand side of \eqref{eqn_degenerate_brion} is simply equal to the number of lattice points inside of $\mathcal{Q}(b_1, \dots, b_\ell)$. If we set all of our $b_i$ equal to $m$ and let $m$ vary, then we are simply taking dilates of a fixed polytope, call it $Q$. Then Theorem \ref{prop_degenerate_brion} says that the number of lattice points inside of dilates of $Q$ grows as a polynomial in $m$ whose degree is equal to the dimension of $Q$, i.e. we recover the Ehrhart polynomial of $Q$. Thus Theorem \ref{prop_degenerate_brion} may be seen as a result interpolating between Brion's formula and the Ehrhart polynomial.
\end{remark}

%% file: sections/spectral_bound.tex
\subsection{The tempered spectrum and exceptional locus} \label{sec_tempered_pgl_3}
Recall the parametrization $\mathcal{S} \subset i \mathfrak{a}^*$ of $\Omega^+_{\textnormal{temp}}$ from Section \ref{sec_parametrization_tempered}. Given $s = (s_1, s_2, s_3) \in \mathcal{S}$, we let $q^s := (q^{s_1}, q^{s_2}, q^{s_3}) \in \Omega^+_{\textnormal{temp}}$. We have an $\mathfrak{S}_3$-action on $\Omega^+_{\textnormal{temp}}$ and on $i \mathfrak{a}^*$ corresponding to permuting coordinates; if, e.g., $\sigma = (1 \ 2 \ 3)$, then we set $\sigma.(s_1, s_2, s_3) = (s_3, s_1, s_2)$, i.e. $\sigma$ permutes coordinates rather than indices.

Recall that $p^\dag = (4/3, -2/3, -2/3)$ and $\delta = (1, 0, -1)$. Let $L := \{1, (2 \ 3)\}$ which is the stabilizer of $p^\dag$. Notice that $(p^\dag, -\sigma_1.s + \sigma_2.s) = (-\sigma_1^{-1}.p^\dag + \sigma_2^{-1}.p^\dag, s)$. If $\sigma_1^{-1} L \neq \sigma_2^{-1} L$, then $\sigma_1^{-1}.p^\dag + \sigma_2^{-1}.p^\dag = 2 \cdot \sigma_3.\delta$ for some $\sigma_3 \in \mathfrak{S}_3$. Notice then that $q^{(\delta, \sigma.s)} = 1$ for some $\sigma$ if and only if $q^{(p^\dag/2, -\sigma_1.s + \sigma_2.s)} = 1$ for some $\sigma_1^{-1} L \neq \sigma_2^{-1} L$. Similarly, $q^{(\delta, \sigma.s)} = -1$ for some $\sigma$ if and only if $q^{(p^\dag/2, -\sigma_1.s + \sigma_2.s)} = -1$ for some $\sigma_1^{-1} L \neq \sigma_2^{-1} L$.

The exceptional locus $\Xi \subset \Omega_{\textnormal{temp}}^+$ is composed of two pieces.
\begin{align*}
\Xi_1 &:= \{q^s \in \Omega^+_{\textnormal{temp}} : \sigma.q^s = q^s \textnormal{ for some $\sigma \neq 1$}\}  = \{q^s \in \Omega^+_{\textnormal{temp}}: q^{(\delta, \sigma.s)} = 1 \textnormal{ for some } \sigma\}\\
&\phantom{:}= \{q^s \in \Omega^+_{\textnormal{temp}} : q^{(p^\dag/2, -\sigma_1.s + \sigma_2.s)} = 1 \textnormal{ for some } \sigma_1^{-1} L \neq \sigma_2^{-1} L\} \\
\Xi_2 & := \{q^s \in \Omega^+_{\textnormal{temp}}: q^{(\delta, \sigma.s)} = -1 \textnormal{ for some } \sigma\}\\
& \phantom{:} = \{q^s \in \Omega^+_{\textnormal{temp}} : q^{(p^\dag, -\sigma_1.s + \sigma_2.s)} = -1 \textnormal{ for some } \sigma_1^{-1} L \neq \sigma_2^{-1} L\}.
\end{align*}
Notice that $\Xi_1 \cup \Xi_2$ is exactly those points $q^{is}$ such that $q^{(p^\dag, -\sigma_1.s + \sigma_2.s)} = 1$ for some $\sigma_1^{-1} L \neq \sigma_2^{-1} L$. Figure \ref{fig_tempered_spectrum} shows a visualization for the tempered spectrum. The red hexagon represents a fundamental domain for the lattice $\frac{2 \pi i}{\ln(q)} \Lambda^* \subset i\mathfrak{a}^*$ as described in Section \ref{sec_parametrization_tempered}; the lattice $\frac{2 \pi i}{\ln(q)} \Lambda^*$ is represented by the black dots. The pink region represents $\mathcal{S}$, whose points are in bijection with $\Omega^+_{\textnormal{temp}}$. The solid brown lines denote $\Xi_1$ and the solid blue lines represent $\Xi_2$. The dotted brown lines represent the loci $(\delta, \cdot) = \frac{\pi i}{\ln(q)} (2 \mathbb{Z})$, and the dotted blue lines represent the loci $(\delta, \cdot) = \frac{\pi i}{\ln(q)} (1 + 2 \mathbb{Z})$. The green point denotes the point $\mathfrak{b}$, i.e. the fixed point of the $\mathbb{Z}/3 \mathbb{Z}$-action on $\Omega^+_{\textnormal{temp}}$ described in Section \ref{sec_coloring_eigenfunctions}.

\begin{figure}[!h]
\centering
\begin{minipage}{.5\textwidth}
  \centering
  \subfile{../images/exceptional}
  \caption{}
  \label{fig_tempered_spectrum}
\end{minipage}%
\begin{minipage}{.5\textwidth}
  \centering
  \subfile{../images/P_data}
  \caption{}
  \label{fig_p_data}
\end{minipage}
\end{figure}

\subsection{Explicit formula for $|h_m(s)|^2$}
Suppose $Y = \Gamma \backslash G / K$. Let $\psi_s$ be an eigenfunction of the $H(G, K)$-action on $L^2(Y)$ with spectral parameter $q^{s} \in \Omega^+_{\textnormal{temp}}$. We may think of $\psi_s$ as a $(\Gamma, K)$-invariant function on $G$. Let $\omega_s$ be the spherical function with Satake parameter $q^{s}$ (see Section \ref{sec_rep_background}). For any function $\mathfrak{h} \in H(G, K)$, we know that $\psi_s * \mathfrak{h} = \hat{\mathfrak{h}}(s) \psi_s$ with $\hat{\mathfrak{h}}(s) := (\omega_s * \mathfrak{h})(1)$. By \cite{macdonald-symmetric}, p. 299 we know that:
\begin{gather}
g_\lambda(s) := (\omega_s * \mathds{1}_{K \varpi^{\lambda} K})(1) = \frac{1}{\nu_\lambda(q^{-1})} \sum_{\sigma \in \mathfrak{S}_3} \bigg( c \big(-\sigma.s \big) \ q^{( \lambda, \delta - \sigma.s )} \bigg)
\end{gather}
where $\varpi^\lambda \in A^+$ and the $c$-function $c(s)$ is defined as
\begin{gather}
c(s) = \prod_{j < k} \frac{q^{s_j} - q^{-1} q^{s_k}}{q^{s_j} - q^{s_k}}, \label{eqn_c_fn}
\end{gather}
and
\begin{equation}
\begin{matrix*}[l]
\nu_{\lambda}(q^{-1}) =  1 & \textnormal{if } \lambda_1 > \lambda_2 > \lambda_3, \\
\nu_{\lambda}(q^{-1}) =  1 + q^{-1} & \textnormal{if } \lambda_1 = \lambda_2 > \lambda_3 \textnormal{ or } \lambda_1 > \lambda_2 = \lambda_3, \label{eqn_nu_lambda_exact}\\
\nu_{\lambda}(q^{-1}) =  (1+q^{-1})(1+q^{-1} + q^{-2}) & \textnormal{if } \lambda_1 = \lambda_2 = \lambda_3.
\end{matrix*}
\end{equation}

We are in particular interested in computing $h_m(s)$, which is defined by $\psi_s * \mathds{1}_{E_m} = h_m(s) \psi_s$. Hence $h_m(s) = (\omega_s * \mathds{1}_{E_m})(1)$. On the other hand
$\mathds{1}_{E_m} = \sum_{\lambda \in P_m^\Lambda} \mathds{1}_{K \varpi^{\lambda} K}$. We therefore arrive at the formula
\begin{gather}
h_m(s) = \sum_{\sigma \in \mathfrak{S}_3} \sum_{\lambda \in P_m^\Lambda} \frac{1}{\nu_\lambda(q^{-1})} \bigg( c \big(-\sigma.s \big) \ q^{( \lambda, \delta - \sigma.s )} \bigg) . \nonumber
\end{gather}

If $s \in \mathcal{S}$ and $\sigma \in \mathfrak{S}_3$ are fixed and we ignore the factor $\frac{1}{\nu_\lambda(q^{-1})}$, then the inner sum reduces to a sum of an exponential function over $P_m^\Lambda$ and hence may be handled by Brion's formula. However, $\nu_\lambda(q^{-1})^{-1}$ is not constant in $\lambda$: it is $\alpha_1 := 1$ on the interior of the Weyl chamber $\mathfrak{a}^+$, it is $\alpha_2 := (1 + q^{-1})^{-1}$ on the intersection of $\Lambda$ with the interior of each extremal ray of $\mathfrak{a}^+$, and it is $\alpha_3 := ((1+q^{-1})(1+q^{-1} + q^{-2}))^{-1}$ at the base vertex of $\mathfrak{a}^+$ namely $(0, 0, 0) \in \mathfrak{a}$. Therefore, using inclusion-exclusion and the labels as in Figure \ref{fig_p_data}, we can write
\begin{align*}
 h_m(s) &= \sum_{\sigma \in \mathfrak{S}_3} c(-\sigma.s) \Bigg(\alpha_1 \sum_{\lambda \in P_m^\Lambda}   q^{( \lambda, \delta-\sigma.s )} + (\alpha_3 - 2 \alpha_2 + \alpha_1) \\
 & + (\alpha_2 - \alpha_1)  \Big(\sum_{\lambda \in (e_{1, 2})_m^\Lambda}  q^{( \lambda, \delta-\sigma.s )}  + \sum_{\lambda \in (e_{1, 3})_m^\Lambda}   q^{( \lambda, \delta-\sigma.s )} \Big) \Bigg). \\
\end{align*}

\noindent We can use Brion's formula (Theorem \ref{thm_brion}) to compute each term:
\begin{align*}
& h_m(s) = \sum_{\sigma \in \mathfrak{S}_3} c(-\sigma.s) \Bigg( \alpha_1 \sum_{j \in \{1, 2, 3\}} \varsigma(\textnormal{Cone}_P(p_j); \delta - \sigma.s) q^{m ( p_j, \delta - \sigma.s )} + (\alpha_3 - 2 \alpha_2 + \alpha_1) + \\
& (\alpha_2 - \alpha_1) \Big(\sum_{j \in \{1, 2\}} \varsigma(\textnormal{Cone}_{e_{1,2}}(p_j); \delta - \sigma.s) q^{m ( p_j, \delta - \sigma.s )} + \sum_{j \in \{1, 3\}} \varsigma(\textnormal{Cone}_{e_{1,3}}(p_j); \delta - \sigma.s) q^{m ( p_j, \delta - \sigma.s )} \Big) \Bigg).
\end{align*}

\subsection{Proof of Propositions \ref{prop_vol_n_lambda} and \ref{prop_vol_e_m}} \label{sec_pf_n_lambda}
Now let's compute $\card(E_m)$. From \cite{macdonald-symmetric} p. 298 we have that:
\begin{gather} \label{eqn_n_lambda}
N_\lambda := \textnormal{vol}(K \varpi^\lambda K) = q^{2 ( \lambda, \delta )} \frac{\nu_3(q^{-1})}{\nu_\lambda(q^{-1})}.
\end{gather}
Here $\nu_3(q^{-1}) = \prod_{i = 1}^3 \frac{1-q^{-i}}{1-q^{-1}}$; notice that $\nu_3(q^{-1})$ is non-zero and does not depend on $\lambda$. Using the same style of analysis as before, we get that
\begin{align}
& \card(E_m) = \sum_{\lambda \in P_m^\Lambda} N_\lambda \nonumber = \nu_3(q^{-1}) \Bigg( \alpha_1 \sum_{j \in \{1, 2, 3\}} \varsigma(\textnormal{Cone}_P(p_j); 2 \delta) q^{2 m ( p_j, \delta )} + (\alpha_3 - 2 \alpha_2 + \alpha_1) \nonumber \\
& + (\alpha_2 - \alpha_1) \Big(\sum_{j \in \{1, 2\}} \varsigma(\textnormal{Cone}_{e_{1, 2}}(p_j); 2 \delta) q^{2 m ( p_j, \delta )} + \sum_{j \in \{1, 3\}} \varsigma(\textnormal{Cone}_{e_{1, 3}}(p_j); 2 \delta) q^{2 m ( p_j, \delta )} \Big) \Bigg) \label{eqn_vol_e_m}
\end{align}

\begin{proof}[Proof of Proposition \ref{prop_vol_n_lambda} (Upper Bound on $N_\lambda$)]
From \eqref{eqn_n_lambda} and \eqref{eqn_nu_lambda_exact}, it is clear that
$N_\lambda \leq \nu_3(q^{-1}) q^{2 ( \delta, \lambda )}$.
\end{proof}

\begin{proof}[Proof of Proposition \ref{prop_vol_e_m} (Lower Bound on $\card(E_m)$)]
Recall that $\card(E_m)$ is computed by summing up $N_\lambda$ over $P_m^\Lambda$. Furthermore $N_\lambda$ is positive for every $\lambda$. We always have $m \cdot p^\dag \in P_m^\Lambda$. If $m \geq 1$, then $N_{m \cdot p^\dag} = \frac{\nu_3(q^{-1})}{1+q^{-1}} \ q^{2 ( \delta, m \cdot p^\dag )} \leq \card(E_m)$.
If $m = 0$ we get $N_{m \cdot p^\dag} = 1 = q^{2 ( \delta, m \cdot p^\dag )} = \card(E_m)$.
\end{proof}

\subsection{Proof of Proposition \ref{spectral_bound}} \label{pf_spectral_bound}
\begin{proof}[Proof of Proposition \ref{spectral_bound} (Spectral Bound)]
We wish to find a lower bound for $\frac{1}{M} \sum_{m = 1}^M \frac{|h_m(s)|^2}{\card(E_m)}$. We first seek to bound $\frac{|h_m(s)|^2}{\card(E_m)}$. In Figure \ref{fig_p_data} one may find a labelled copy of the polytope $P$. Notice that $p^\dag$ is the unique vertex of $P$ maximizing $|q^{( \lambda, \delta - \sigma.s )}|$ (recall that $q^{(\lambda, s)}$ is on the unit circle), and that $p^\dag = p_3$ is part of $P$ and $e_{1, 3}$. Furthermore, observe that the $c$-function \eqref{eqn_c_fn} is holomorphic away from $\Xi_1$. Hence because we are assuming that the closure of $\Theta$ (which is a compact set) does not intersect $\Xi_1$, we conclude that the norm of the $c$-function has a universal upper bound on $\Theta$. Also note that all of the terms of the form $\varsigma(C; \delta - w.s)$, where $C$ is some cone, are independent of $m$ and also universally upper bounded in magnitude for all $q^s \in \Theta$ (in fact for all $q^s \in \Omega^+_{\textnormal{temp}}$). Therefore, we obtain
\begin{gather}
h_m(s) = \Big(\sum_{\sigma \in \mathfrak{S}_3} \kappa(\sigma.s) q^{m ( p^\dag, -\sigma.s )} \Big) q^{m ( \delta, p^\dag )} + R_1(m, s) q^{m (( p^\dag, \delta ) - \eta_1)}, \label{asymptotic}
\end{gather}
where
\begin{gather}
\kappa(s) = c(-s) \Big[ \varsigma(\textnormal{Cone}_P(p^\dag); \delta - s) - \frac{\varsigma(\textnormal{Cone}_{e_{1, 3}}(p^\dag); \delta - s)}{q+1}\Big], \label{kappa}
\end{gather}
and $R_1(m, s)$ is some function which is universally bounded for all $m$ and all $q^s \in \Theta$, and $\eta_1$ is some positive constant.

Since we also have that $\kappa(\sigma.s)$ is bounded for all $\sigma \in \mathfrak{S}_3$ and $q^s \in \Theta$, we get that
\begin{gather}
|h_m(s)|^2 = \big|\sum_{\sigma \in \mathfrak{S}_3} \kappa(\sigma.s) q^{m ( p^\dag, -\sigma.s )} \big|^2 q^{2 m ( \delta, p^\dag )} + R_2(m, s) q^{2m (( p^\dag, \delta ) - \eta_2)}, \nonumber
\end{gather}
where $R_2(m, s)$ is some function which is universally bounded for all $m$ and all $q^s \in \Theta$, and $\eta_2$ is some positive constant.

We now analyze $\card(E_m)$. Again, because $p^\dag$ is the unique vertex in $P$ maximizing $q^{2 ( \lambda, \delta )}$, we see that
\begin{gather}
\card(E_m) = (D + R_3(m)) q^{2m ( \delta, p^\dag )}, \nonumber
\end{gather}
where
\begin{gather}
D = \nu_3(q^{-1}) \Bigg(\varsigma(\textnormal{Cone}_P(p^\dag); 2 \delta) + (\frac{1}{1 + q^{-1}} - 1) \varsigma(\textnormal{Cone}_{e_{1, 3}}(p^\dag); 2 \delta) \Bigg), \label{eqn_spectral_bound_pf_d}
\end{gather}
and $R_3(m) \to 0$ as $m \to \infty$. A straightforward calculation shows that $D$ is positive.

Therefore, using the identity
\begin{gather}
\frac{C_1}{C_2 + R} = \frac{C_1}{C_2} - \frac{C_1 R}{C_2 (C_2 + R)}, \nonumber
\end{gather}
we obtain
\begin{align}
\frac{|h_m(s)|^2}{\card(E_m)} &= \frac{|\sum_{\sigma \in \mathfrak{S}_3} \kappa(\sigma.s) q^{m ( p^\dag, -\sigma.s )} |^2 q^{2 m ( \delta, p^\dag )} + R_2(m, s) q^{2m (( p^\dag, \delta ) - \eta_2)}}{(D + R_3(m)) q^{2m ( \delta, p^\dag )}} \nonumber \\
    &= \frac{|\sum_{\sigma \in \mathfrak{S}_3} \kappa(\sigma.s) q^{m ( p^\dag, -\sigma.s )} |^2}{D} \nonumber \\
    & - \frac{|\sum_{\sigma \in \mathfrak{S}_3} \kappa(\sigma.s) q^{m ( p^\dag, -\sigma.s )} |^2 R_3(m)}{D (D + R_3(m))} + \frac{R_2(m, s)}{(D + R_3(m)) q^{2m \cdot \eta_2}}. \label{eqn_remainder_spectral}
\end{align}
Because $R_3(m) \to 0$ as $m \to \infty$, it is clear that given $\varepsilon > 0$, there exists an $M_0$ such that for all $m \geq M_0$ and for all $q^s \in \Theta$, the first term in \eqref{eqn_remainder_spectral} is bounded in absolute value by $\varepsilon$. Furthermore, since $\eta_2 > 0$, it is also clear that there exists an $M_1$ such that for all $m \geq M_1$ and for all $s \in \Theta$, the second term in \eqref{eqn_remainder_spectral} is bounded in absolute value by $\varepsilon$.

Let $R_4(m, s)$ equal \eqref{eqn_remainder_spectral}. We define $R_5(m)$:
\begin{gather}
R_5(m) := \sup_{\{s : q^s \in \Theta \}} |R_4(m, s)|. \nonumber
\end{gather}
It is clear that $R_5(m) \to 0$ as $m \to \infty$. We therefore have that for all $q^s \in \Theta$,
\begin{gather}
\Big|\frac{1}{M} \sum_{m = 1}^M R_4(m, s) \Big| \leq \frac{1}{M} \sum_{m = 1}^M R_5(m) \to 0 \nonumber
\end{gather}
as $m \to \infty$.

In the sequel we shall prove the following lemma:
\begin{lemma} \label{lemma_spectral_bound_pf}
There exist an $M_2 \in \mathbb{N}$ and $C > 0$ such that for all $M \geq M_2$ and for all $q^s \in \Theta$, 
\begin{gather}
\frac{1}{M} \sum_{m = 1}^M \Big|\sum_{\sigma \in \mathfrak{S}_3} \kappa(\sigma.s) q^{m ( p^\dag, -\sigma.s )} \Big|^2 \geq C. \nonumber
\end{gather}
\end{lemma}

Given this lemma we can finish the proof of Proposition \ref{spectral_bound}. Let $M_3$ be such that for all $M \geq M_3$, we have
\begin{gather}
\frac{1}{M} \sum_{m = 1}^M R_5(m) \leq \frac{C}{2 D}, \nonumber
\end{gather}
with $C$ given by Lemma \ref{lemma_spectral_bound_pf} and $D$ given by \eqref{eqn_spectral_bound_pf_d}. We then have that if $M \geq M_2$ and $M \geq M_3$, then
\begin{align*}
\frac{1}{M} \sum_{m = 1}^M \frac{|h_m(s)|^2}{\card(E_m)} &= \frac{1}{M} \sum_{m = 1}^M \frac{|\sum_{\sigma \in \mathfrak{S}_3} \kappa(\sigma.s) q^{m ( p^\dag, -\sigma.s )} |^2}{D} + \frac{1}{M} \sum_{m = 1}^M R_4(m, s) \\
    & \geq \frac{C}{D} - \frac{C}{2D} = \frac{C}{2 D}.
\end{align*}
\end{proof}

\begin{proof}[Proof of Lemma \ref{lemma_spectral_bound_pf}]
Recall that $p^\dag = (4/3, -2/3, -2/3)$, which is fixed by the permutation $(23)$ (and the trivial permutation). Let $L = \{1, (23)\}$. Notice that if $\sigma_1 \sigma_2^{-1} \in L$ (i.e. $\sigma_1$ and $\sigma_2$ send the same element to 1), then $q^{( p^\dag, -\sigma_1.s + \sigma_2.s )} = 1$. Notice also that $\overline{\kappa(s)} = \kappa(-s)$. Therefore we shall expand and group terms as follows
\begin{align}
\frac{1}{M} \sum_{m = 1}^M \Big|\sum_{\sigma \in \mathfrak{S}_3} \kappa(\sigma.s) q^{m ( p^\dag, -\sigma.s )}\Big|^2 &= \frac{1}{M} \sum_{m = 1}^M \sum_{\sigma_1, \sigma_2 \in \mathfrak{S}_3} \kappa(\sigma_1.s) \kappa(-\sigma_2.s) q^{m ( p^\dag, -\sigma_1.s + \sigma_2.s )} \nonumber \\
    &= \sum_{\textnormal{right cosets $L. \sigma$}} |\kappa(\sigma.s) + \kappa\big((23)\sigma.s\big)|^2 + \textnormal{ r.t.}, \label{eqn_grouped_spectral_bound}
\end{align}
where r.t. stands for remainder terms. We wish to show that the sum in \eqref{eqn_grouped_spectral_bound} is non-zero. Evidently this amounts to showing that at least one $\kappa(\sigma.s)+\kappa\big((23)\sigma.s\big) \neq 0$ for every $s$.

Recall that $s = (s_1, s_2, s_3)$ with $s_1 + s_2 + s_3 = 0$ and each $s_i$ purely imaginary and only defined up to $2 \pi i \mathbb{Z}/\ln(q)$. For the polytope $P$, the coprime generators for the rays generating the cone at $p^\dag$ are $\{(-\frac{2}{3}, \frac{1}{3}, \frac{1}{3}), (-1, 1, 0)\}$. For the polytope $e_{1, 3}$, the coprime generator for the ray generating the cone at $p^\dag$ is $(-\frac{2}{3}, \frac{1}{3}, \frac{1}{3})$. Therefore, using \eqref{eqn_cone_generator}, we obtain
\begin{align}
\varsigma(\textnormal{Cone}_P(p^\dag); \delta - s) &= \frac{1}{1-q^{( \delta - s, (-\frac{2}{3}, \frac{1}{3}, \frac{1}{3}))}} \cdot \frac{1}{1-q^{( \delta - s, (-1, 1, 0) )}} 
    = \frac{1}{1-q^{-1+s_1}} \cdot \frac{1}{1-q^{-1+s_1-s_2}} \label{cone1} \\
\varsigma(\textnormal{Cone}_{e_{1,3}}(p^\dag); \delta - s) &= \frac{1}{1 - q^{( \delta - s, (-\frac{2}{3}, \frac{1}{3}, \frac{1}{3}) )}} = \frac{1}{1-q^{-1+s_1}}. \label{cone2}
\end{align}

Using (\ref{eqn_c_fn}), (\ref{kappa}), (\ref{cone1}), and (\ref{cone2}), we can perform a change a variables $x_1 = q^{-s_1}$, $x_2 = q^{-s_2}$, $x_3 = q^{-s_3}$, $t = q^{-1}$ and write each term as a rational function in these variables. We also have the condition $x_1 x_2 x_3 = 1$ coming from $q^{-s_2-s_2-s_3} = 1$. Using a computer algebra system (such as Mathematica \cite{Mathematica}), we find, for example, that under this change of variables $\kappa(s) + \kappa\big((23).s\big)$ becomes
\begin{equation}
\frac{x_1 (x_1^2 - t^3 x_2 x_3)}{(-t+x_1)(x_1 - x_2) (x_1 - x_3)}. \nonumber
\end{equation}
If this were zero, then $x_1 (x_1^2 - t^3 x_2 x_3) = 0$. Then $x_1^3 = t^3 x_1 x_2 x_3 = t^3$. However, $x_1$ is on the unit circle, and $t^3 = \frac{1}{q^3} < 1$. Hence this expression is never zero on the locus of interest to us.

Now we focus on the remainder terms which are all of the form:
\begin{gather}
\frac{1}{M} \sum_{m = 1}^M \Big(\kappa(\sigma_1.s) + \kappa\big((23)\sigma_1.s\big)\Big) \overline{\Big(\kappa(\sigma_2.s) + \kappa\big((23)\sigma_2.s\big)\Big)} q^{m ( p^\dag, -\sigma_1.s + \sigma_2.s )}, \nonumber
\end{gather}
where $\sigma_1, \sigma_2$ satisfy $L .\sigma_1 \neq L. \sigma_2$. Since the closure of $\Theta$ does not intersect $\Xi_1$, and hence the $c$-function is bounded on $\Theta$, we get the following bound for all $q^s \in \Theta$:
\begin{gather}
\Big|\Big(\kappa(\sigma_1.s) + \kappa\big((23)\sigma_1.s\big)\Big) \overline{\Big(\kappa(\sigma_2.s) + \kappa\big((23)\sigma_2.s\big)\Big)}\Big| \leq F \nonumber
\end{gather}
for all choices of $\sigma_1$ and $\sigma_2$, and some $F$. Furthermore, since the closure of $\Theta$ also does not intersect $\Xi_2$, we get that $|q^{( p^\dag, -\sigma_1.s + \sigma_2.s )} - 1| \geq \eta > 0$
for all $q^s \in \Theta$, all choices of $\sigma_1$ and $\sigma_2$, and some $\eta > 0$. Therefore, using geometric series, we get that each remainder term (of which there are finitely many) can be uniformly bounded by an expression of the form 
\begin{gather}
F \frac{q^{M ( p^\dag, -\sigma_1.s + \sigma_2.s )} - q^{( p^\dag, -\sigma_1.s + \sigma_2.s )}}{M (q^{( p^\dag, -\sigma_1.s + \sigma_2.s )} - 1)} \leq \frac{2 F}{M \eta}. \nonumber
\end{gather}
This bound goes to zero as $M \to \infty$. Hence Lemma \ref{lemma_spectral_bound_pf} is proven. 
\end{proof}

%% file: images/exceptional.tex
\begin{tikzpicture}[scale=0.6]
\draw[color = red, line width = 1.5] (3.000000, 1.732051) -- (0.000000, 3.464102) -- (-3.000000, 1.732051) -- (-3.000000, -1.732051) -- (-0.000000, -3.464102) -- (3.000000, -1.732051) -- cycle;
\fill[color = pink, opacity = 0.5] (0, 0) -- (3.000000, 1.732051) -- (3.000000, -1.732051) -- cycle;
\draw[color = brown, line width = 2] (0, 0) -- (3.000000, 1.732051) -- (3, -1.732051) -- cycle;

\draw[dotted, color = brown, line width = 1] (0, 6.928203) -- (0, -6.928203);
\draw[dotted, color = brown, line width = 1] (-3.000000, -6.928203) -- (-3.000000, 6.928203);
\draw[dotted, color = brown, line width = 1] (3, -6.928203) -- (3, 6.928203);
\draw[dotted, color = brown, line width = 1] (6, -6.928203) -- (6, 6.928203);
\draw[dotted, color = brown, line width = 1] (-6, -6.928203) -- (-6, 6.928203);

\draw[ color = blue, line width = 2] (1.5, 0.866025) -- (1.5, -0.866025) -- (3, 0) -- cycle;

\draw[dotted, color = blue, line width = 1] (-1.5, -6.928203) -- (-1.5, 6.928203);
\draw[dotted, color = blue, line width = 1] (1.5, -6.928203) -- (1.5, 6.928203);
\draw[dotted, color = blue, line width = 1] (-4.5, -6.928203) -- (-4.5, 6.928203);
\draw[dotted, color = blue, line width = 1] (4.5, -6.928203) -- (4.5, 6.928203);






\node (x11) at (1.500000, 0.000000) {};
\node (x12) at (0.750000, 1.299038) {};
\node (x2) at (2.250000, 1.299038) {};





\fill[color = green] (2, 0) circle (3pt);

\fill[color = black] (0.000000, 0.000000) circle (3pt);
\fill[color = black] (6.000000, 0.000000) circle (3pt);
\fill[color = black] (3.000000, 5.196152) circle (3pt);
\fill[color = black] (0.000000, 0.000000) circle (3pt);
\fill[color = black] (-3.000000, 5.196152) circle (3pt);
\fill[color = black] (3.000000, 5.196152) circle (3pt);
\fill[color = black] (-0.000000, 0.000000) circle (3pt);
\fill[color = black] (-3.000000, 5.196152) circle (3pt);
\fill[color = black] (-6.000000, 0.000000) circle (3pt);
\fill[color = black] (-0.000000, -0.000000) circle (3pt);
\fill[color = black] (-3.000000, -5.196152) circle (3pt);
\fill[color = black] (-6.000000, -0.000000) circle (3pt);
\fill[color = black] (0.000000, -0.000000) circle (3pt);
\fill[color = black] (-3.000000, -5.196152) circle (3pt);
\fill[color = black] (3.000000, -5.196152) circle (3pt);
\fill[color = black] (0.000000, 0.000000) circle (3pt);
\fill[color = black] (6.000000, 0.000000) circle (3pt);
\fill[color = black] (3.000000, -5.196152) circle (3pt);
\end{tikzpicture}

%% file: images/P_data.tex
\begin{tikzpicture}[scale=2]
\draw (-1.000000, 0.000000) -- (2.500000, 0.000000);
\draw (-1.000000, 0.866025) -- (2.500000, 0.866025);
\draw (-1.000000, 1.732051) -- (2.500000, 1.732051);
\draw (-1.000000, 0.000000) -- (2.500000, 0.000000);
\draw (-1.000000, -0.866025) -- (2.500000, -0.866025);
\draw (-0.577350, -1.000000) -- (1.443376, 2.500000);
\draw (-1.000000, 0.000000) -- (0.443376, 2.500000);
\draw (-1.000000, 1.732051) -- (-0.556624, 2.500000);
\draw (-0.577350, -1.000000) -- (1.443376, 2.500000);
\draw (0.422650, -1.000000) -- (2.443376, 2.500000);
\draw (1.422650, -1.000000) -- (2.500000, 0.866025);
\draw (2.422650, -1.000000) -- (2.500000, -0.866025);
\draw (0.577350, -1.000000) -- (-1.000000, 1.732051);
\draw (1.577350, -1.000000) -- (-0.443376, 2.500000);
\draw (2.500000, -0.866025) -- (0.556624, 2.500000);
\draw (2.500000, 0.866025) -- (1.556624, 2.500000);
\draw (0.577350, -1.000000) -- (-1.000000, 1.732051);
\draw (-0.422650, -1.000000) -- (-1.000000, 0.000000);
\path[draw, color = brown, line width = 3] (0, 0) -- (2.500000, 0.000000);
\path[draw, color = brown, line width = 3] (0, 0) -- (1.443376, 2.500000);
\fill[color=red, opacity = 0.5] (0, 0) -- (2.000000, 0.000000) -- (0.500000, 0.866025) -- cycle;
\draw[color=red, line width=3] (0, 0) -- (2.000000, 0.000000);
\draw[color=red, line width=3] (0, 0) -- (0.500000, 0.866025);
\draw[color=red, line width = 3] (2.000000, 0.000000) -- (0.500000, 0.866025);
\fill[color=blue] (0, 0) circle (2pt);
\fill[color=blue] (2.000000, 0.000000) circle (2pt);
\fill[color=blue] (0.500000, 0.866025) circle (2pt);
\node (p1) at (0, 0) {};
\node (p2) at (0.500000, 0.866025) {};
\node (p3) at (2.000000, 0.000000) {};
\node[color=blue] (p1_label) at (-0.375000, -0.216506) {$p_1$};
\node[color=blue] (p2_label) at (0.125000, 1.082532) {$p_2$};
\node[color=blue] (p3_label) at (2.000000, -0.433013) {$p_3 = p^\dag$};
\draw[->] (p1_label) edge (p1);
\draw[->] (p2_label) edge (p2);
\draw[->] (p3_label) edge (p3);
\node (e12) at (0.250000, 0.433013) {}; 
\node (e13) at (0.500000, 0.000000) {}; 
\node (e23) at (1.062500, 0.541266) {}; 
\node[color=red] (e12_label) at (-0.100000, 0.692820) {$e_{1, 2}$}; 
\node[color=red] (e13_label) at (0.500000, -0.433013) {$e_{1, 3}$}; 
\node[color=red] (e23_label) at (1.312500, 0.974279) {$e_{2, 3}$};
\draw[->] (e12_label) edge (e12);
\draw[->] (e13_label) edge (e13);
\draw[->] (e23_label) edge (e23);
\end{tikzpicture}

%% file: sections/BS_convergence_fixed.tex
\subsection{Benjamini-Schramm convergence implies Plancherel convergence}
Suppose $M$ is a topological group. Let $\hat{M}$ denote the collection of irreducible unitary representations of $M$. Suppose $\Gamma$ is a cocompact lattice in $M$. We then have
\begin{gather}
L^2(\Gamma \backslash M) = \bigoplus_{\rho \in \hat{M}} N_\Gamma(\rho) \rho \nonumber
\end{gather}
with each $N_\Gamma(\rho)$ finite, and only countably many of them non-zero. Let $C_c^\infty(M)$ denote the space of test functions on $M$ (see, e.g., Definition 1.2 of \cite{deitmar}); in case $M = \textnormal{PGL}(d, F)$, then $C_c^\infty(M)$ is the set of compactly supported locally constant functions on $M$. Each $f \in C_c^\infty(M)$ defines an operator $\hat{f}(\rho)$ on the underlying vector space $V_\rho$ via:
\begin{gather*}
    \hat{f}(\rho).v := \int_M f(m) \rho(m^{-1}).v \ dm,
\end{gather*}
with $v \in V_\rho$. It turns out this operator is in fact trace class (\cite{deitmar_echterhoff}, Theorem 9.3.2).

We shall define a measure on $\hat{M}$, called the spectral measure associated to $\Gamma$, by
\begin{gather}
\mu_\Gamma := \frac{1}{\vol(\Gamma \backslash M)} \sum_{\rho \in \hat{M}} N_\Gamma(\rho) \delta_\rho, \nonumber
\end{gather}
where $\delta_\rho$ is the indicator function for $\rho \in \hat{M}$. Suppose $(\Gamma_n)$ is a sequence of cocompact lattices and assume that $M$ is a type I group. We say that $(\Gamma_n)$ is a {\it Plancherel sequence} if, for every $f \in C_c^\infty(M)$,
\begin{gather}
 \int_{\hat{M}} \textnormal{Tr}[\hat{f}(\rho) ]\ d \mu_{\Gamma_n} (\rho) \to \int_{\hat{M}} \textnormal{Tr} [\hat{f}(\rho)] \ d \mu(\rho), \nonumber
\end{gather}
where $\mu$ is the Plancherel measure. We refer to this type of convergence as {\it Plancherel convergence}.

We say that $(\Gamma_n)$ {\it Benjamini-Schramm converges to $\{1\}$} if, for every compact subset $C \subset M$,
\begin{gather}
\frac{\textnormal{vol}(\{x \in \Gamma_n \backslash M: x^{-1} (\Gamma_n \setminus \{1\}) x \cap C \neq \emptyset \})}{\textnormal{vol}(\Gamma_n \backslash M)} \to 0. \nonumber
\end{gather}
Here in the numerator we are using the volume on $\Gamma_n \backslash M$ since the set in question is clearly $(\Gamma_n, 1)$-invariant.

Suppose $G = \textnormal{PGL}(d, F)$. Then any compact subset $C \subset G$ is contained in a disjoint union of finitely many $K$ double cosets (since such sets are compact and open); such double cosets are always of the form $K \varpi^\lambda K$ where $\lambda$ is some partition. Let $\tau_\Gamma: G/K \to \Gamma \backslash G / K$ and $\pi: G \to G/K$ denote the obvious projections. 

\begin{proposition} \label{prop_interpret_deitmar}
Suppose $\Gamma$ is a torsion-free lattice. Let $Y = \Gamma \backslash G /K$. Then
\begin{align*}
&\vol(\{x \in \Gamma \backslash G: x^{-1} (\Gamma \setminus \{1\}) x \cap K \varpi^\lambda K \neq \emptyset\}) \\
&= \card(\{y \in Y: \exists \ \tilde{y_1}, \tilde{y_2} \in \tau_{\Gamma}^{-1}(y) \subset G/K \textnormal{ with } \tilde{y_1} \neq \tilde{y_2} \textnormal{ and } d_{A^+}(\tilde{y_1}, \tilde{y_2}) = \lambda \}).
\end{align*}
\end{proposition}

Let $B_R$ denote the collection of vertices contained in a ball of radius $R$ centered at $1K$ in the building. Let $C_R := \pi^{-1}(B_R)$ be its preimage in $G$. Then $C_R$ is the disjoint union of finitely $K$ double cosets. Using Proposition \ref{prop_interpret_deitmar}, we thus have that
\begin{gather*}
\vol(\{x \in \Gamma \backslash G: x^{-1} (\Gamma \setminus \{1\}) x^{-1} \cap C_R \neq \emptyset\}) = \card(\{y \in Y: \textnormal{InjRad}_Y(y) \leq R \}).
\end{gather*}

Suppose $Y_n = \Gamma_n \backslash G / K$ is a sequence of compact quotients of the Bruhat-Tits building. In Section \ref{sec_defn_BS} we defined what it means for $Y_n$ to Benjamini-Schramm converge to $G/K$. The above observations reveal the relationship between these two notions of Benjamini-Schramm convergence.

\begin{proposition}[see also Proposition 2.4 in \cite{deitmar}]
A sequence of torsion-free lattices $(\Gamma_n)$ in $G$ Benjamini-Schramm converges to $\{1\}$ if and only if the sequence of spaces $(Y_n)$ Benjamini-Schramm converges to $G/K$.
\end{proposition}

We say that $(\Gamma_n)$ is {\it uniformly discrete} if there exists a neighborhood $U$ of the identity such that $x^{-1} \Gamma_n x \cap U = \{1\}$ for every $n$ and for every $x \in M$. Suppose $M$ is a semisimple algebraic group over a non-archimedean local field; then $M$ contains a compact open subgroup $K$. If $(\Gamma_n)$ is any sequence of torsion-free lattices in $M$, then $x^{-1}\Gamma_n x \cap K = \{1\}$ for all $n$ and all $x$ (since the intersection would be discrete and compact, i.e. finite, and hence would contain a finite order element). On the other hand since $K$ is open, we immediately conclude that $(\Gamma_n)$ is uniformly discrete.

In \cite{deitmar} we find the following result:
\begin{theorem} [\cite{deitmar} Theorem 2.6; see also Theorem 1.3 of \cite{gelander_levit}] \label{thm_deitmar}
Suppose $(\Gamma_n)$ is a sequence of cocompact, uniformly discrete lattices in a locally compact group $M$. Then the following are equivalent
\begin{enumerate}
\item $(\Gamma_n)$ is a Plancherel sequence.
\item $(\Gamma_n)$ Benjamini-Schramm converges to $\{1\}$.
\end{enumerate}
\end{theorem}
\noindent Our goal now is to relate the property of being a Plancherel sequence to the type of convergence in Theorem \ref{thm_bs_conv}.

\subsection{The Fell topology vs. the Euclidean topology on $\Omega^+_{\textnormal{temp}}$} \label{sec_compare_topologies}
Let $G = \textnormal{PGL}(d, F)$ and $K = \textnormal{PGL}(d, \mathcal{O})$. Let $\mu$ denote the Plancherel measure. Let $\Omega^+$ denote the class 1 representations of $G$. Given $\chi$, a (not necessarily unitary) unramified character of the maximal torus $T$, we can construct a representation of $G$, denoted $I_\chi$, via parabolic induction from $\chi$ (twisted by the modular character of the Borel subgroup). Any two such $I_\chi, I_{\chi'}$ do not share any irreducible subquotients unless $\chi$ and $\chi'$ are in the same Weyl group orbit.

Each element of $\Omega^+$ arises as a subquotient of a unique $I_\chi$. Note that the group of unramified characters of $T$ is isomorphic to the group of characters of $A$, which is in turn isomorphic to the lattice $\mathbb{Z}^d/\mathbb{Z} \cdot (1, \dots, 1)$. Let $\tilde{A}$ denote the collection of characters of $A$, and let $\hat{A}$ denote the collection of unitary characters of $A$. These space carries the natural topology of pointwise convergence of characters. These spaces are in turn homeomorphic to $(\mathbb{C}^\times)^{d-1}$ and $(S^1)^{d-1}$, respectively, with their usual Euclidean topologies.

We have an injection $f: \Omega^+ \to \tilde{A}/\mathfrak{S}_d$. The former space naturally carries the Fell topology and the latter space carries a ``Euclidean topology''. A theorem of Tadic (\cite{tadic-1983}, Theorem 5.6) tells us that with respect to these topologies $f$ is continuous. Furthermore, since $K$ is an open subgroup, $\Omega^+_{\textnormal{temp}}$ is compact with respect to the Fell topology (see \cite{macdonald-spherical-function}, p. 12). The space $\tilde{A}/\mathfrak{S}_d$ is clearly Hausdorff. Thus $f$ is a continuous bijective map from a compact space to a Hausdorff space (namely its image in $\tilde{A}/\mathfrak{S}_d$), and thus it must be a homeomorphism.

\begin{proposition} \label{prop_homeo}
The natural map $f: \Omega^+ \to \tilde{A}/\mathfrak{S}_d$ is a homeomorphism onto its image when viewing $\Omega^+$ with the Fell topology and $\hat{A}/\mathfrak{S}_d$ with the Euclidean topology. In particular this map provides a homeomorphism between $\Omega^+_{\textnormal{temp}}$ and $\hat{A}/\mathfrak{S}_d$. 
\end{proposition}

\subsection{Plancherel convergence vs. weak convergence}

In \cite{sauvageot}, one finds a result often referred to as the Sauvageot density principle. It says that when $M$ is a product of groups, each one a linear reductive group over a local field of characteristic zero, then Plancherel convergence implies that
\begin{gather}
\mu_{\Gamma_n}(E) \to \mu(E) \label{weak_conv}
\end{gather}
for every $E$ which is relatively compact and $\mu$-regular with respect to the Fell topology (see e.g. Remark 1.6 of \cite{deitmar}). However, a gap was recently found in the proof and noted in \cite{nelson_venkatesh}. If the Sauvageot density principle is true, then Theorem \ref{thm_bs_conv}, at least in characteric zero, could be deduced from Theorem \ref{thm_deitmar}, Proposition \ref{prop_homeo}, and the fact that the Plancherel measure is absolutely continuous with respect to the Lebesgue measure. However, we shall instead prove a weaker form of the Sauvageot density principle which only concerns $\Omega^+$. This has the added benefit of allowing us to prove Theorem \ref{thm_bs_conv} in both zero and positive characteristic.

\begin{proposition} \label{prop_sauvageot}
Suppose $E$ is a measurable subset of $\Omega^+$ with $\mu(\partial E) = 0$. Let $\varepsilon > 0$ be given. Then there exist $\phi, \psi \in C_c^\infty(G)$ such that
\begin{enumerate}
    \item $0 \leq \mathbf{1}_E(\rho) - \textnormal{Tr}[\hat{\phi}(\rho)] \leq \textnormal{Tr}[\hat{\psi}(\rho)]$ for all $\rho \in \hat{G}$,
    \item $\int_{\hat{G}} \textnormal{Tr}[\hat{\psi}] d \mu \leq \varepsilon$.
\end{enumerate}
\end{proposition}

\begin{proof}
The main ingredients of the proof are the Stone-Weierstrass theorem and Urysohn's lemma. First off, note that if $\phi \in H(G, K)$, then $\textnormal{Tr}[\hat{\phi}(\rho)] = 0$ if $\rho$ is not class 1, and $\textnormal{Tr}[\hat{\phi}(\rho)]$ is equal to the eigenvalue associated to convolving $\phi$ with the spherical function associated to $\rho$ in case $\rho$ is class 1.

If $\phi$ is the indicator of $K$, then $\textnormal{Tr}[\hat{\phi}(\rho)]$ is identically equal to one on $\Omega^+$. Recall that the Satake isomorphism identifies $H(G, K)$ with the coordinate ring of $\tilde{A}/\mathfrak{S}_d$ when viewed as an affine variety. The points on this variety parametrize the collection of all spherical representations. Given $\phi \in H(G, K)$, its image $\theta(\phi)$ under the Satake isomorphism is a complex-valued function on $\tilde{A}/\mathfrak{S}_d$ whose value at a point $\rho$, when viewed as a class 1 representation, is equal to the eigenvalue associated to convolving the associated spherical function with $\phi$. Since coordinate rings of affine varieties separate points, we get that the algebra of functions on $\Omega^+ \subset \tilde{A}/\mathfrak{S}_d$ arising from taking the trace of elements in $H(G, K)$ separates points in $\Omega^+$. Since $\Omega^+$ is a compact Hausdorff space, and since the image of $H(G, K)$ is an algebra of continuous functions containing constants and separating points, by the Stone-Weierstrass theorem we get that this image is dense in the space of all continuous functions on $\Omega^+$ with the topology of uniform convergence. This means that for any continuous function $g$ on $\Omega^+$, and any choice of $\delta > 0$, we can find a $\psi \in H(G, K)$ such that $|g(\rho) - \hat{\psi}(\rho)| < \delta$ for all $\rho \in \Omega^+$. 

We now use Urysohn's lemma to approximate indicator functions by continuous functions. Suppose $\delta > 0$. Recall that $\mu$ is absolutely continuous with respect to Lebesgue measure, and Lebesgue measure is a regular measure. Because $\mu(\partial E) = 0$, we can find an open set $A$ covering $\partial E$ such that $\mu(A) < \varepsilon$. Let $W = \bar{E} \cap A^c$. Then $W$ is closed, $\mu(E \setminus W) < \delta$, and $W \cap \bar{E^c} = \emptyset$. Let $U$ be an open set containing $\bar{E}$ such that $\mu(U \setminus E) < \delta$. By Urysohn's lemma, we can find continuous functions $f, g$ on $\Omega^+$ taking values in $[0, 1]$ such that $f|_{W} \equiv 1$, $f|_{\bar{E^c}} \equiv 0$, $g|_{\bar{E}} \equiv 1$, and $g|_{U^c} \equiv 0$. We shall rescale $f$ by a factor of $1-\delta$, and rescale $g$ by a factor of $1 + \delta$. We thus have $0 \leq f + \delta \leq \mathbf{1}_E \leq g - \delta \leq 1$. We may then use the Stone-Weierstrass theorem to find $\tilde{f} = \textnormal{Tr}[\hat{\alpha}]$ and $\tilde{g} = \textnormal{Tr}[\hat{\beta}]$ such that $|\tilde{f}(\rho) - f(\rho)| < \delta$ and $|\tilde{g}(\rho) - g(\rho)| < \delta$ for all $\rho$, and such that $\alpha, \beta \in H(G, K)$. We then have $\tilde{f} \leq \mathbf{1}_E \leq \tilde{g}$. 

Writing $\tilde{g} - \tilde{f} = (\tilde{g} - g) + (g - f) + (f - \tilde{f})$, we get that
\begin{gather*}
\int_{\hat{G}} \tilde{\psi} - \tilde{\phi} \ d \mu \leq \delta \cdot \mu(\Omega^+) + 2 \delta \cdot \mu(E) + (2 + 2 \delta) \cdot \delta + \delta \cdot \mu(\Omega^+).
\end{gather*}
We can choose $\delta$ to make the right-hand side of the above inequality less than $\varepsilon$. After picking such a $\delta$, we may thus take $\phi = \alpha$ and $\psi = \beta - \alpha$ in the statement of the proposition.
\end{proof}

\subsection{Proof of Theorem \ref{thm_bs_conv}} \label{pf_bs_plancherel}
\begin{proof}[Proof of Theorem \ref{thm_bs_conv} (BS Convergence Implies Plancherel Convergence)]
Let $E \subset \Omega^+_{\textnormal{temp}}$ be measurable with measure zero boundary. Let $\varepsilon > 0$ be given and let $\phi, \psi \in H(G, K)$ be as in Proposition \ref{prop_sauvageot}. We have
\begin{align*}
    \Big|\frac{N(E, Y_n)}{\card(Y_n)} - \mu(E) \Big| &=  \Big| \int_{\hat{G}} \mathds{1}_E d \mu_{\Gamma_n} - \int_{\hat{G}} \mathds{1}_E d \mu \Big| \\
    & = \Big| \int_{\hat{G}} (\mathds{1}_E - \textnormal{Tr} [\hat{\phi}]) + \textnormal{Tr} [\hat{\phi}] \ d \mu_{\Gamma_n} - \Big(\int_{\hat{G}} (\mathds{1}_E -\textnormal{Tr} [\hat{\phi}]) + \textnormal{Tr} [\hat{\phi}] \ d \mu \Big) \Big| \\
    & \leq \int_{\hat{G}} \textnormal{Tr} [\hat{\psi}] \ d \mu_{\Gamma_n} + \Big| \int_{\hat{G}} \textnormal{Tr} [\hat{\phi}] d \mu_{\Gamma_n} - \int_{\hat{G}} \textnormal{Tr} [\hat{\phi}] d \mu \Big| + \varepsilon \\ 
    & \leq \Big| \int_{\hat{G}} \textnormal{Tr} [\hat{\psi}] \ d \mu_{\Gamma_n} - \int_{\hat{G}} \textnormal{Tr} [\hat{\psi}] d \mu\Big|  + \Big| \int_{\hat{G}} \textnormal{Tr} [\hat{\phi}] d \mu_{\Gamma_n} - \int_{\hat{G}} \textnormal{Tr} [\hat{\phi}] d \mu \Big| + 2 \varepsilon
\end{align*}

We now follow an adaptation of the argument in Proposition 2.9 in \cite{deitmar}. Let $C_\phi, C_\psi \subset G$ denote the supports of $\phi, \psi$; such sets are necessarily $(K, K)$-invariant. Using the trace formula, we get that
\begin{align*}
\Big| \int_{\hat{G}} \textnormal{Tr}[\hat{\phi}] d \mu_{\Gamma_n} - \int_{\hat{G}} \textnormal{Tr}[\hat{\phi}] d \mu \Big| &= \frac{1}{\card(Y_n)} \Big| \int_{\Gamma_n \backslash G} \sum_{\gamma \in \Gamma_n \setminus \{1\}} \phi(x^{-1} \gamma x) d x \Big| \\
    & \leq \frac{||\phi||_\infty}{\card(Y_n)} \int_{\Gamma_n \backslash G} \# \{x^{-1} (\Gamma_n \setminus \{1\}) x \cap C_\phi \} \\
    & = \frac{||\phi||_\infty}{\card(Y_n)} \sum_{x \in \Gamma_n \backslash G / K} \#\{x^{-1} (\Gamma_n \setminus \{1\}) x \cap C_\phi\}.
\end{align*}

We could bound this expression as follows: we can enlarge $C_\phi$ to some set $C_R$ equal to the preimage of the collection of vertices contained in a ball of radius $R$ centered at $1K$ in the building. Given $x \in Y_n$, let $\tilde{x}$ be any choice of lift in $G / K$. Then 
\begin{align*}
\#\{x^{-1}(\Gamma_n \backslash \{1\}) x \cap C_\phi\} &\leq \#\{x^{-1}(\Gamma_n \backslash \{1\}) x \cap C_R\} \\
    & = \# \{\gamma \in \Gamma_n \setminus \{1\} : d(\tilde{x}, \gamma \cdot \tilde{x}) \leq R \} \\
    & \leq \frac{\vol_{G/K} (B_R)}{\vol_{G/K} (B_{\textnormal{InjRad}_{Y_n}(x)})}.
\end{align*}
Here $\vol_{G/K}(B_R)$ denotes the number of vertices in the ball of radius $R$ in the building centered at any vertex. Therefore
\begin{align*}
\sum_{x \in \Gamma_n \backslash G / K} \# \{x^{-1}(\Gamma_n \setminus \{1\}) x \cap C_\phi \} \leq \card(\{x : \textnormal{InjRad}_{Y_n}(x) \leq R\}) \cdot \frac{\vol_{G / K}(B_R)}{\vol_{G/ K}(B_{\textnormal{InjRad}(Y_n)})}.
\end{align*}
Putting everything together we obtain:
\begin{align*}
\Big| \frac{N(E, Y_n)}{\card(Y_n)} - \mu(E) \Big| \leq &2 \varepsilon + \frac{1}{\vol_{G/K}(B_{\textnormal{InjRad}(Y_n)})} \times \\
&\Bigg( ||\phi||_\infty \cdot \frac{\card (\{x \in Y_n: \textnormal{InjRad}_{Y_n}(x) \leq R_{\phi}\})}{\card(Y_n)} \cdot \vol_{G/ K} (B_{R_\phi}) \\
&+ ||\psi||_\infty \cdot \frac{\card (\{x \in Y_n: \textnormal{InjRad}_{Y_n}(x) \leq R_{\psi}\})}{\card(Y_n)} \cdot \vol_{G/ K} (B_{R_\psi}) \Bigg),
\end{align*}
where $R_\phi$ and $R_\psi$ are the radii of the support of $\phi$ and $\psi$ as functions on $G/K$.

In Section \ref{kernel_section} we further discuss $\vol_{G/K}(B_R)$. In particular there exist constant $C_1, C_2$ such that $C_1 q^{2 R} \leq \vol_{G/K}(B_R) \leq C_2 R q^{2 R}$ for all $R$. 
\end{proof}

%% file: sections/kernel.tex
\subsection{The volume of balls in $\mathcal{B}$} \label{sec_vol_balls}
Recall that elements in $G/K$ correspond to vertices on the building $\mathcal{B}$. Let $d(\cdot, \cdot)$ be the Euclidean metric on $\mathcal{B}$ normalized so that adjacent vertices are distance 1 apart. With respect to this metric $\mathcal{B}$ is a CAT(0) space, and hence in particular there is a unique geodesic joining any two points. The group $G$ acts on $\mathcal{B}$ by isometries with respect to this metric on $\mathcal{B}$.

There is a natural measure on $\mathcal{B}$ coming from the Lebesgue measure on each apartment; we shall denote this measure by $\vol_{d(\cdot, \cdot)}$. Under the above normalization, each chamber has measure $\frac{\sqrt{3}}{4}$ (the area of an equilateral triangle with side lengths 1). Since $G$ acts transitively on vertices, the volume of a ball of radius $R$ centered at any vertex is independent of that vertex. We let $B(x, R)$ denote the ball of radius $R$ centered at vertex $x$, and we let $\vol_{d(\cdot, \cdot)}(B_R)$ denote the volume of the ball of radius $R$ centered at any vertex.

There is another natural measure on $\mathcal{B}$ coming from the Haar measure on $G$ which is obtained by simply counting the number of vertices in a given set. We shall refer to this as $\vol_{G/K}$ (this is essentially the same as what has previously been referred to as $\card(\cdot)$). 

\begin{proposition}[Lemma 2 in Leuzinger \cite{leuzinger}]
There exist constants $C, D > 0$ such that for all $R \geq 0$,
\begin{gather}
\vol_{d(\cdot, \cdot)}(B_R) \leq C \cdot \vol_{G/K}(B_R) \leq \vol_{d(\cdot, \cdot)}(B_{R + D}). \nonumber
\end{gather}
\end{proposition}

On the other hand, Leuzinger \cite{leuzinger} also shows that

\begin{proposition}[\cite{leuzinger}, p. 487] \label{leuzinger2}
There exist constants $C_1, C_2 > 0$ such that for all $R \geq 0$,
\begin{gather}
C_1 q^{2 R} \leq \vol_{G/K}(B_R) \leq C_2 R q^{2 R}.\nonumber
\end{gather}
\end{proposition}
Combining these facts, we get the following:
\begin{proposition} \label{leuzinger}
There exists $C > 0$ such that for all $R$,
\begin{gather}
\vol_{d(\cdot, \cdot)}(B_R) \leq C R q^{2 R} . \nonumber
\end{gather}
\end{proposition}

We shall also make use of the following proposition.
\begin{proposition} \label{n_gamma}
There exist a universal constant $C > 0$ such that for every $R$, for every lattice $\Gamma$, and for every $z, w \in G/K$,
\begin{gather}
N_\Gamma(z, w; R) := \# \{ \gamma \in \Gamma : d(z, \gamma . w) \leq R \} \leq \frac{C R q^{2 R}}{\textnormal{InjRad}(\Gamma \backslash G / K)^{2}}. \label{eqn_n_gamma}
\end{gather}
\end{proposition}

\begin{proof}
This is simply the combination of (the proof of) Lemma 6.18 in \cite{7_samurai} and (the proof of) Lemma 5.1 in \cite{brumley_matz}. Though the proof contains no new ideas besides the cited lemmas, a complete proof may be found in Proposition VI.5 of \cite{my_thesis} .
\end{proof}

\subsection{Proof of Lemma \ref{lemma_kernel}} \label{pf_kernel_lemma}
\begin{proof}[Proof of Lemma \ref{lemma_kernel} (Lifting the Kernel to $G/K$)]
Let $D$ be a fundamental domain in $G/K$ for the action of $\Gamma$ and let $Y = \Gamma \backslash G /K$. Given $x \in Y$, we let $\tilde{x}$ denote its unique lift to $D$, and given $w \in D$, we let $\bar{w}$ denote its projection to $Y$.

Suppose $\mathcal{K}: G/K \times G/K \to \mathbb{C}$ is invariant under the diagonal $\Gamma$-action and satisfies $\mathcal{K}(z, w) = 0$ if $d(z, w) > R$ for some $R \geq 0$. Let $\mathcal{K}^{\textnormal{Op}}$ denote the associated operator on $L^2(G/K)$. Then $\mathcal{K}^{\text{Op}}$ preserves $(\Gamma, 1)$-invariant functions and hence descends to a well-defined operator on $L^2(Y)$ which we denote $\bar{\mathcal{K}}^{\text{Op}}$. On this space this operator is represented by the kernel
\begin{gather}
\bar{\mathcal{K}}(x, y) = \sum_{\gamma \in \Gamma} \mathcal{K}(\tilde{x}, \gamma . \tilde{y}),
\end{gather}
with $x, y \in Y$. Thus the Hilbert-Schmidt norm of $\bar{\mathcal{K}}^{\text{Op}}$ is
\begin{gather}
||\bar{\mathcal{K}}^{\text{Op}}||_\textnormal{HS}^2 = \sum_{y \in Y} \sum_{x \in Y} \Big|\sum_{\gamma \in \Gamma} \mathcal{K}(\tilde{x}, \gamma . \tilde{y})\Big|^2  = \sum_{\tilde{y} \in D} \sum_{\tilde{x} \in D} \Big|\sum_{\gamma \in \Gamma} \mathcal{K}(\tilde{x}, \gamma . \tilde{y})\Big|^2. \nonumber
\end{gather}

Now we use the assumption that $\mathcal{K}(z, w) = 0$ whenever $d(z, w) > R$ to bound this Hilbert-Schmidt norm. Let $I(R)$ be the points in $Y$ with injectivity radius greater than $R$. Let $I(R)^C$ be the remaining points in $Y$. We claim that if $y \in I(R)$, then the injectivity radius at $y$ is greater than $R$, and thus $d(z, \gamma . \tilde{y}) \leq R$ for at most one $\gamma \in \Gamma$ for a given $z \in G/K$.

Note that $y \in I(R)$ implies that $d(\tilde{y}, \gamma . \tilde{y}) > 2 R$ for $\gamma \neq 1$. Suppose $d(z, \tilde{y}) \leq R$. Suppose for the sake of contradiction that $d(z, \gamma . \tilde{y}) \leq R$ for some $\gamma \neq 1$. Then by joining geodesic segments from $\tilde{y}$ to $z$, then from $z$ to $\gamma . \tilde{y}$, we'd get that $d(\tilde{y}, \gamma . \tilde{y}) \leq 2R$, a contradiction. Hence there is at most one $\gamma$ such that $d(z, \gamma . \tilde{y}) \leq R$. This implies that for $y \in I(R)$, we have 
\begin{gather}
\Big|\sum_{\gamma \in \Gamma} \mathcal{K}(z, \gamma . \tilde{y})\Big|^2 = \sum_{\gamma \in \Gamma} |K(z, \gamma . \tilde{y})|^2 \nonumber
\end{gather}
since both sides have at most one term.

We clearly have $|\mathcal{K}(\tilde{x}, \gamma . w)| \leq ||\mathcal{K}||_\infty$, and by Cauchy-Schwarz we have:
\begin{align*}
\Big|\sum_{\gamma \in \Gamma} \mathcal{K}(\tilde{x}, \gamma . w) \cdot 1 \Big|^2 & \leq N_\Gamma(\tilde{x}, w; R) \sum_{\gamma \in \Gamma} |\mathcal{K}(\tilde{x}, \gamma . w)|^2.
\end{align*}

Therefore, combining this with Proposition \ref{n_gamma}, we get
\begin{align*}
\sum_{\bar{w} \in I(R)^C} \sum_{\tilde{x} \in D} \Big| \sum_{\gamma \in \Gamma} \mathcal{K}(\tilde{x}, \gamma . w) \Big|^2 & \leq \frac{C R q^{2 R}}{\textnormal{InjRad}(Y)^{2}} \sum_{\bar{w} \in I(R)^C} \sum_{\tilde{x} \in D} \sum_{\gamma \in \Gamma} |\mathcal{K}(\tilde{x}, \gamma . w)|^2 \\
    &= \frac{C R q^{2 R}}{\textnormal{InjRad}(Y)^{2}} \sum_{y \in I(R)^C} \sum_{\{u \in G/K: d(u, \tilde{y}) \leq R\}} |\mathcal{K}(u, \tilde{y})|^2 \\
    &\leq \frac{C R q^{2 R}}{\textnormal{InjRad}(Y)^{2}}  \vol_{G/K}(B_R) ||\mathcal{K}||_\infty^2  \card(I(R)^C).
\end{align*}

Putting everything together and using Proposition \ref{leuzinger} we get that there exists $\tilde{C} \geq 0$ such that:
\begin{align*}
||\bar{\mathcal{K}}^{\text{Op}}||_\textnormal{HS}^2 &= \sum_{\bar{w} \in I(R)} \sum_{\tilde{x} \in D} \Big|\sum_{\gamma \in \Gamma} \mathcal{K}(\tilde{x}, \gamma. w)\Big|^2 + \sum_{\bar{w} \in I(R)^c} \sum_{\tilde{x} \in D} |\sum_{\gamma \in \Gamma} \mathcal{K}(\tilde{x}, \gamma .w)|^2 \\
    & \leq \sum_{\bar{w} \in I(R)} \sum_{\tilde{x} \in D} \sum_{\gamma \in \Gamma} |\mathcal{K}(\tilde{x}, \gamma. w)|^2 + \frac{C R q^{2 R}}{\textnormal{InjRad}(Y)^{2}} \vol_{G/K}(B_R) ||\mathcal{K}||_\infty^2 \card(I(R)^C) \\
    & \leq \sum_{z \in D} \sum_{w \in G/K} |\mathcal{K}(z, w)|^2 + \frac{\tilde{C} R^2 q^{4 R}}{\textnormal{InjRad}(Y)^{2}} ||\mathcal{K}||_\infty^2 \card(\{y \in Y:\textnormal{InjRad}_Y(y) \leq R\}).
\end{align*}
\end{proof}

%% file: sections/changing_variables.tex
\subsection{Geometric interpretation of $\Gamma \backslash G / M_{\lambda}$}
Let $o$ denote the point $1K$ in $G/K$. Given $\lambda \in A^+$, let $z_\lambda := \varpi^\lambda K$. Then $d_{A^+}(o, z_\lambda) = \lambda$. Let $M_\lambda$ be the joint (pointwise) stabilizer of $o$ and $z_\lambda$, namely $M_\lambda = K \cap (\varpi^\lambda K \varpi^{-\lambda})$. Let $D$ be a fundamental domain for the $\Gamma$-action on $G/K$. Recall that we always assume that $\Gamma$ is torsionfree.

\begin{proposition}
Assuming that the Haar measure on $G$ is normalized such that $\vol(K) = 1$, we have
\begin{gather}
\vol(M_\lambda) = \frac{1}{\vol(K \varpi^\lambda K)} = \frac{1}{N_\lambda}. \nonumber
\end{gather}
\end{proposition}
\begin{proof}
We clearly have $M_\lambda < K$. Consider the right action of $K$ on right cosets of $M_\lambda$ in $K$. Given a coset $M_\lambda k$, we may associate to it the point $z_k := k^{-1} z_\lambda$. This point is clearly independent of the choice of coset representative (since $M_\lambda$ stabilizes $z_\lambda$). Furthermore, if $k' \in K$ satisfies $(k')^{-1}.z_\lambda = k^{-1}.z_\lambda$, then $k' k^{-1} \in M_\lambda$, and thus $M_\lambda k = M_\lambda k'$. Therefore we have an injection from the cosets $M_\lambda k$ to the points in the $K$-orbit of $z_\lambda$. However, by the Cartan decomposition we know that the $K$-orbit of $z_\lambda$ is exactly those points $w$ satisfying $d_{A^+}(o, w) = \lambda$. Therefore we have a surjection, and hence a bijection, between the cosets of $M_\lambda$ in $K$ and those points in the $K$ orbit of $z_\lambda$, namely $K \varpi^\lambda K$. Therefore, the number of $M_\lambda$ cosets is exactly equal to $\card(K z_\lambda) = \vol(K \varpi^\lambda K)$. This in turn is equal to the index of $M_\lambda$ in $K$.
\end{proof}

\begin{proposition} \label{prop_identify_double_coset}
Let $D_\lambda$ be defined as
\begin{gather}
    D_\lambda := \{(x, y) \in D \times G/K : d_{A^+}(x, y) = \lambda \}. \nonumber
\end{gather}
Consider the map
\begin{gather}
    \mathfrak{f}_\lambda: D_\lambda \to \Gamma \backslash G / M_\lambda \nonumber
\end{gather}
defined as follows: given $(x, y) \in D_\lambda$, choose any $g \in G$ such that $g.(o, z_\lambda) = (x, y)$; set
\begin{gather}
    \mathfrak{f}_\lambda(x, y) = \Gamma \backslash g / M_\lambda. \nonumber
\end{gather}
Then $\mathfrak{f}_\lambda$ is well-defined and defines a bijection between $D_\lambda$ and $\Gamma \backslash G / M_\lambda$. 
\end{proposition}

\begin{proof}
First we show that a $g$ as in the statement of the proposition exists. We can clearly find some $g_1$ such that $g_1.x = o$. Then $d_{{A^+}}(o, g_1.y) = \lambda$. Suppose $g_1.y = h K$. Writing $h = k_1 a k_2$ with $a \in A^+$ (the Cartan decomposition), we must have that $a = \varpi^\lambda$. Thus $k_1^{-1} h K = \varpi^\lambda K$, so we may take $g = g_1^{-1} k_1$.

Next we show that the map $\mathfrak{f}_\lambda$ is well-defined, i.e. does not depend on the choice of $g$. If $g$ and $g'$ both map $(o, z_\alpha)$ to $(x, y)$, then $g^{-1} g' \in M_\lambda$, so they generate the same double coset. 

We now show injectivity of $\mathfrak{f}_\lambda$. Consider the set $G_D := \{g \in G: g.o \in D\}$. Then clearly the collection of $g$'s which arise in the definition of the map $\mathfrak{f}_\lambda$ (for all choices of $(x, y) \in D_\lambda$) is exactly $G_D$. Suppose $g_1, g_2 \in G_D$ are in the same $(\Gamma, M_\lambda)$-double coset. Then $g_1 = \gamma g_2 m$ with $\gamma \in \Gamma$ and $m \in M_\lambda$. Therefore $g_1.o = \gamma g_2 m.o = \gamma g_2.o$. Since $g_1 \in G_D$, we know that $g_1.o \in D$. Since we also have $g_2 \in G_D$, we must have $\gamma = 1$. Therefore $g_1.o = g_2.o$. We clearly then also have that $g_1.z_\lambda = \gamma g_2 m.z_\lambda = g_2.z_\lambda$. Therefore the double coset $\Gamma g_1 M_\lambda = \Gamma g_2 M_\lambda$ has exactly one preimage, namely $(g_1.o, g_1.z_\lambda)$. Therefore $\mathfrak{f}_\lambda$ is injective.

Lastly we show that $\mathfrak{f}_\lambda$ is surjective. Let $g$ be some double coset representative. Then for a unique $\gamma \in \Gamma$, $\gamma g.o \in D$. Thus this double coset is the image under $\mathfrak{f}_\lambda$ of $(\gamma g.o, \gamma g.z_\lambda) \in D_\lambda$.
\end{proof}

\begin{proposition} \label{prop_bij_g_gamma}
Suppose $\{g_i\}$ is a complete set of double coset reprensenatives of $\Gamma \backslash G /M_\lambda$. Then the map:
\begin{gather}
\mathfrak{j}_\lambda: \Gamma \backslash G /M_\lambda \times M_\lambda \to \Gamma \backslash G, \ \ \ \ \ (g_i, m) \mapsto \Gamma g_i m, \nonumber
\end{gather}
is a bijection. Furthermore this identification is measure-preserving when we take $M_\lambda$ to have measure $1/N_\lambda$.
\end{proposition}
\begin{proof}
Because $\Gamma$ is torsionfree (by assumption), no conjugate of $\Gamma$ can intersect $M_\lambda$ non-trivially because the intersection would be a compact discrete group which necessarily has torsion. We claim that this implies that every $g \in G$ may be represented uniquely as $g = \gamma g_j m$ with $\gamma \in \Gamma$ and $m \in M_\lambda$. To see this, first note that $g_j$ clearly must be the unique element in $\{g_j\}$ generating the same $(\Gamma, M_\lambda)$-double coset as $g$. Suppose we have some equation of the following form:
\begin{gather}
\gamma_1 g_j m_1 = \gamma_2 g_j m_2 \iff g_j^{-1} (\gamma_2^{-1} \gamma_1) g_j = m_2 m_1^{-1}. \nonumber
\end{gather}
Because we know that $g_j^{-1} \Gamma g_j \cap M_\lambda = \{1\}$, we must have $\gamma_1 = \gamma_2$ and $m_1 = m_2$.

From this we can conclude that $\mathfrak{j}_\lambda$ is injective: if $\Gamma g_j m_1 = \Gamma g_k m_2$, then $\gamma_1 g_j m_1 = \gamma_2 g_k m_2$, and thus $g_j = g_k$ and $m_1 = m_2$. 

Lastly we show that $\mathfrak{j}_\lambda$ is surjective. Given a coset $\Gamma g$, we write $g = \gamma g_j m$. Therefore $\mathfrak{j}_\lambda(g_j, m) = \Gamma g$. 

Clearly this map is measure-preserving as $\Gamma \backslash G / M_\lambda$ is a discrete set.
\end{proof}

\subsection{Proof of Proposition \ref{prop_change_variables}} \label{pf_changing_variables}
\begin{proof}[Proof of Proposition \ref{prop_change_variables} (Changing Variables in the Kernel Integral)]
Suppose $x \in G/K$. We now wish to rewrite the following integral:
\begin{gather}
\sum_{x \in D} \sum_{y \in G/K} \Big|\frac{1}{M} \sum_{m = 1}^M \frac{1}{\card(E_m)} \sum_{z \in x E_m \cap y E_m} a(z) \Big|^2. \label{eqn_change_variables_start}
\end{gather}
Recall the definition of the following sets:
\begin{gather*}
E_m^\lambda := E_m \cap \varpi^\lambda E_m, \hspace{10mm}
\tilde{E}_m^\lambda := \pi^{-1} (E_m^\lambda).
\end{gather*}
Notice that both sets are setwise invariant under left multiplication by elements in $M_\lambda$. Also notice that we have a (non-canonical) identification
$\tilde{E}_m^\lambda \simeq E_m^\lambda \times K$
by choosing coset representatives $\tilde{w}$ for each $w \in E_m^\lambda$ (recall $w \in G/K$). This identification is measure-preserving because $E_m^\lambda$ is a discrete set.

Notice that if $\mathfrak{f}_\lambda(x, y) = \Gamma g M_\lambda$ as in Proposition \ref{prop_identify_double_coset}, then $gE_m^\lambda = xE_m \cap yE_m$. We therefore have (recalling that $a$ may be viewed as a $(1, K)$-invariant function on $\Gamma \backslash G$):

\begin{align*} [\rho^{\Gamma}_{\tilde{E}_m^\lambda}.a](\Gamma g) & = \frac{1}{\vol(\tilde{E}_m^\lambda)} \int_{\tilde{E}_m^\lambda} a(\Gamma g h) dh = \frac{1}{\vol(\tilde{E}_m^\lambda)} \sum_{w \in E_m^\lambda} \int_K a(\Gamma g \tilde{w} k) dk \\
    &= \frac{1}{\vol(\tilde{E}_m^\lambda)} \sum_{w \in E_m^\lambda} a(\Gamma g \tilde{w}) = \frac{1}{\vol(\tilde{E}_m^\lambda)} \sum_{z \in gE_m^\lambda} a(\Gamma \tilde{z})  = \frac{1}{\vol(\tilde{E}_m^\lambda)} \sum_{z \in xE_m \cap yE_m} a(z).
\end{align*}
Let $\{g_i\}$ be a complete set of coset representatives for $\gmodgamma / M_\lambda$. Thus our original integral from \eqref{eqn_change_variables_start} can now be written as 
\begin{align*}
\sum_{x \in D} \sum_{y \in G/K} \Big|  \sum_{m = 1}^M \frac{1}{\card(E_m)} \sum_{z \in x E_m \cap y E_m} a(z) \Big|^2
 = \sum_{\lambda \in A^+} \sum_{g_i \in \gmodgamma / M_\lambda} \Big| \sum_{m=1}^M \frac{\vol(\tilde{E}_m^\lambda)}{\card(E_m)} \rho^\Gamma_{\tilde{E}_m^\lambda} a(\Gamma g_i) \Big|^2.
\end{align*}

We now utilize Proposition \ref{prop_bij_g_gamma}:
\begin{align*}
\int_{\gmodgamma} \Big| \sum_{m = 1}^M \frac{\vol(\tilde{E}_m^\lambda)}{\card(E_m)} [\rho^\Gamma_{\tilde{E}_m^\lambda}.a](\Gamma g)\Big|^2 dg & = \sum_{g_i \in \gmodgamma / M_\lambda} \sum_{n \in M_\lambda} \Big| \sum_{m = 1}^M\frac{\vol(\tilde{E}_m^\lambda)}{\card(E_m)} [\rho^\Gamma_{\tilde{E}_m^\lambda}.a](\Gamma g_i n)\Big|^2  \\
&= \sum_{g_i \in \Gamma \backslash G / M_\lambda} \sum_{n \in M_\lambda} \Big|  \sum_{m = 1}^M\frac{\vol(\tilde{E}_m^\lambda)}{\card(E_m)} [\rho^\Gamma_{n.\tilde{E}_m^\lambda}.a](\Gamma g_i)\Big|^2  \\
&= \sum_{g_i \in \Gamma \backslash G / M_\lambda} \sum_{n \in M_\lambda} \Big|  \sum_{m = 1}^M\frac{\vol(\tilde{E}_m^\lambda)}{\card(E_m)} [\rho^\Gamma_{\tilde{E}_m^\lambda}.a](\Gamma g_i)\Big|^2  \\
&= \sum_{g_i \in \Gamma \backslash G / M_\lambda} \vol(M_\lambda) \Big|  \sum_{m = 1}^M\frac{\vol(\tilde{E}_m^\lambda)}{\card(E_m)} [\rho^\Gamma_{\tilde{E}_m^\lambda}.a](\Gamma g_i)\Big|^2. 
\end{align*}
Here we have used that $\tilde{E}_m^\lambda$ is invariant (as a set) under $M_\lambda$.
\end{proof}

%% file: sections/kunze_stein.tex
\subsection{The Kunze-Stein phenomenon}
Suppose $M$ is a locally compact topological group. Let $\mathcal{M}(M)$ denote the collection of measurable functions on $M$. We say that $M$ {\it satisfies the Kunze-Stein phenomenon} (or is a {\it KS group}, for short), if for every $1 \leq p < 2$, the map $L^p(M) \times L^2(M) \to \mathcal{M}(M)$ defined by convolution is continuous and has image contained in $L^2(M)$. This means that, for a fixed $p$, there exists a $C_p$ such that for all $f \in L^p(M)$ and $g \in L^2(M)$ we have 
\begin{gather}
||f \ast g||_2 \leq C_p ||f||_p ||g||_2. \label{eqn_ks_op_norm}
\end{gather}
Another way of expressing this is that for all $f \in L^p(M)$ we have $||\hat{f}(\lambda_M)|| \leq C_p ||f||_p$, where $\lambda_M$ is the (left) regular representation and $\hat{f}(\lambda_M)$ is the operator on $L^2(M)$ given by (left) convolution with $f$.

\begin{lemma} [see Lemma 7.1 in \cite{cowling_ks}] \label{finite}
Suppose $M$ is a KS group, and $N \lhd M$ is a compact normal subgroup. Then $M/N$ is a KS group.
\end{lemma}

\begin{lemma} \label{finite_index}
Suppose $N$ is a KS group, and $N \vartriangleleft M$ is a normal subgroup such that $M/N$ is compact. Suppose $N, M$ are unimodular and $\sigma$-finite with respect to Haar measure. Then $M$ is a KS group.
\end{lemma}

\begin{proof}
Choose a measurable section $\tilde{L} \subset M$ of $M \twoheadrightarrow M/N$. Because $M, N$ are unimodular, then, assuming we appropriately normalize our Haar measure, we have that for any measurable function $f$ on $M$ (see, e.g., Theorem B.1.4 of \cite{property_t_book}):
\begin{gather*}
    \int_M f(m) dm = \int_{\tilde{L}} \int_N f(\tilde{\ell} \cdot n) dn d \mu_{M/N}( \tilde{\ell} N).
\end{gather*}
We shall use $d \tilde{\ell}$ as shorthand for $d \mu_{M/N}(\tilde{\ell} N)$. Suppose $a \in L^p(M)$ and $b \in L^2(M)$. Suppose $p = \tilde{\alpha} \cdot \beta \in M$ with $\tilde{\alpha} \in \tilde{L}$ and $\beta \in N$. Then
\begin{align*}
a * b (\tilde{\alpha} \beta) = \int_M a(m) b (m^{-1} \tilde{\alpha} \beta) dm &= \int_{\tilde{L}} \int_N a(\tilde{\ell} n) b (n^{-1} \tilde{\ell}^{-1} \tilde{\alpha} \beta) dn d \tilde{\ell} \\
    &= \int_{\tilde{L}} \int_N a (\tilde{\ell} n) b \big((\tilde{\ell}^{-1} \tilde{\alpha}) (\tilde{\ell}^{-1} \tilde{\alpha})^{-1} n^{-1} (\tilde{\ell}^{-1} \tilde{\alpha}) \beta \big) dn d \tilde{\ell}.
\end{align*}

If $m \in M$, let $\phi_m: N \to N$ denote the automorphism obtained by conjugating by $m$, i.e. $n \mapsto m^{-1} n m$. We define the following functions on $N$:
$\hat{a}_{\tilde{\ell}}(n) := a( \tilde{\ell} n)$ and $\hat{b}_{\tilde{\ell}}^{\tilde{\alpha}}(n) := b \big( \tilde{\ell}^{-1} \tilde{\alpha} \cdot \phi_{\tilde{\ell}^{-1} \tilde{\alpha}}(n) \big)$. Notice then that
\begin{gather*}
\hat{a}_{\tilde{\ell}} * \hat{b}_{\tilde{\ell}}^{\tilde{\alpha}} \Big(\phi_{(\tilde{\ell}^{-1} \tilde{\alpha})^{-1}}(\beta)\Big) = \int_N a (\tilde{\ell} n) b \big((\tilde{\ell}^{-1} \tilde{\alpha}) (\tilde{\ell}^{-1} \tilde{\alpha})^{-1} n^{-1} (\tilde{\ell}^{-1} \tilde{\alpha}) \beta \big) dn.
\end{gather*}

We wish to compute the $L^2$-norm of $a * b$. We have
\begin{align*}
||a * b||_2^2 &= \int_{\tilde{L}} \int_{N} \Big| \int_{\tilde{L}} \hat{a}_{\tilde{\ell}} * \hat{b}_{\tilde{\ell}}^{\tilde{\alpha}} \big( \phi_{(\tilde{\ell}^{-1} \tilde{\alpha})^{-1}}(\beta) \big) d \tilde{\ell} \Big|^2 d\beta d \tilde{\alpha} \\
    &\leq \vol(M/N)^2 \int_{\tilde{L}} \int_N \int_{\tilde{L}} \big| \hat{a}_{\tilde{\ell}} * \hat{b}_{\tilde{\ell}}^{\tilde{\alpha}} \big( \phi_{(\tilde{\ell}^{-1} \tilde{\alpha})^{-1}}(\beta) \big) \big|^2 d \tilde{\ell} d \beta d \tilde{\alpha} \\
    &= \vol(M/N)^2 \int_{\tilde{L}} \int_{\tilde{L}} \int_N \big| \hat{a}_{\tilde{\ell}} * \hat{b}_{\tilde{\ell}}^{\tilde{\alpha}} \big( \phi_{(\tilde{\ell}^{-1} \tilde{\alpha})^{-1}}(\beta) \big) \big|^2 d \beta d \tilde{\ell} d \tilde{\alpha}.
\end{align*}
To move from the first line to the second we have used Cauchy-Schwarz, and to get from the second to the third we have used Fubini's theorem (which holds because $N$ and $\tilde{L}$ are $\sigma$-finite measure spaces by assumption). Note that clearly $||\hat{a}_{\tilde{\ell}}||_p$ is equal to the $L^p$-norm of $a$ restricted to $\tilde{\ell}N$, which we denote by $a_{\tilde{\ell}N}$. Because $N$ is unimodular and $M/N$ is compact, conjugation of $N$ by elements in $M$ preserves the Haar measure on $N$. This implies that for any $p \in M$ and any measurable function $\psi$ on $N$, $\int_N \psi(n) dn = \int_N \psi \big( \phi_p(n) \big) dn$. This implies that $||\hat{b}||_2$ is equal to the $L^2$-norm of $b$ restricted to $\tilde{\ell}^{-1} \tilde{\alpha} N$, which we denote by $b_{\tilde{\ell}^{-1} \tilde{\alpha} N}$. Clearly we have that
\begin{gather*}
  \int_N \big| \hat{a}_{\tilde{\ell}} * \hat{b}_{\tilde{\ell}}^{\tilde{\alpha}} \big( \phi_{(\tilde{\ell}^{-1} \tilde{\alpha})^{-1}}(\beta) \big) \big|^2 d \beta = \int_N \big| \hat{a}_{\tilde{\ell}} * \hat{b}_{\tilde{\ell}}^{\tilde{\alpha}} \big( \beta \big) \big|^2 d \beta. 
\end{gather*}
Therefore, because $N$ is a KS group, we have that the integral on the above line is bounded by $C_p^2 ||a_{\tilde{\ell}N}||_p^2 \cdot ||b_{\tilde{\ell}^{-1} \tilde{\alpha} N}||_2^2$. We therefore have
\begin{gather*}
||a * b||_2^2 \leq \vol(M/N)^2 \int_{\tilde{L}} \int_{\tilde{L}} C_p^2 \cdot ||a_{\tilde{\ell}N}||_p^2 \cdot ||b_{\tilde{\ell}^{-1} \tilde{\alpha} N}||_2^2 d \tilde{\ell} d \tilde{\alpha} = \vol(M/N)^2 C_p^2 \cdot  ||a||_p^2 \cdot ||b||_2^2.
\end{gather*}
\end{proof}

We now wish to show that $G = \textnormal{PGL}(d, F)$ is a KS group. We utilize the above lemmas and the following result.

\begin{theorem} [Veca \cite{veca}] \label{veca}
Let $F$ be a non-discrete, totally disconnected local field, and let $M$ be the group of $F$-rational points of a simply connected algebraic group defined over $F$. Then $M$ is a KS group.
\end{theorem}

\noindent In particular this applies to the group $\textnormal{SL}(d, F)$.

\begin{corollary}
The group $\textnormal{PGL}(d, F)$ is a KS group.
\end{corollary}

\begin{proof}
We have the following long exact sequence:
\begin{gather}
1 \to \mu_d(F) \to \textnormal{SL}(d, F) \to \textnormal{PGL}(d, F) \to F^\times/(F^\times)^d \to 1, \nonumber
\end{gather}
where $\mu_d(F)$ is the group of $d$th roots of unity of $F$. Since $\textnormal{SL}(d, F)$ is a KS group by Theorem \ref{veca}, and since $\mu_d(F) \lhd \textnormal{SL}(d, F)$ is finite, $\textnormal{PSL}(d, F) := \textnormal{SL}(d, F)/\mu_d(F)$ is also a KS group by Lemma \ref{finite}. We claim that $F^\times/(F^\times)^d$ is a compact group; note that if the characteristic is relatively prime to $d$ this group is in fact finite (\cite{serre}), but if the characteristic divides $d$ then this group could be infinite. Recall that after choosing a uniformizer we have an isomorphism $F^\times \simeq \mathbb{Z} \times \mathcal{O}^\times$. Under the projection $F^\times \to F^\times/(F^\times)^d$, the kernel, viewed as a subgroup of $\mathbb{Z} \times \mathcal{O}^\times$, is clearly $(d \mathbb{Z}) \times (\mathcal{O}^\times)^d$. Therefore, $F^\times/(F^\times)^d \simeq (\mathbb{Z}/d \mathbb{Z}) \times \mathcal{O}^\times/(\mathcal{O}^\times)^d$. This group is clearly compact since $\mathcal{O}^\times$ is compact. Thus by Lemma \ref{finite_index}, we get that $\textnormal{PGL}(d, F)$ is a KS group.
\end{proof}

\subsection{Proof of Proposition \ref{nevo}} \label{sec_pf_nevo}
Let $M$ be a semisimple algebraic group over a local field. Suppose $\psi \in L^1(M)$ and $(\rho, \mathcal{H})$ is a unitary representation of $M$. Recall the notation $\hat{\psi}(\rho)$ from Section \ref{bs_implies_plancherel_section}. Given an operator $T$ on a Hilbert space, the notation $||T||$ refers to its operator norm.
\begin{lemma} \label{nevo_lemma}
Suppose $(\rho, \mathcal{H})$ is a unitary representation of $M$. Let $\psi \in L^1(M)$ be such that $\int_M \psi dg = 1$ and $\psi$ is non-negative and real-valued. Then for any even positive integer $2m$,
\begin{gather}
|| \hat{\psi}(\rho) ||^{2m} \leq ||\hat{\psi}(\rho^{\otimes 2m})||. \nonumber
\end{gather}
\end{lemma}

\begin{proof}
This essentially follows from the proof of Theorem 1 in Nevo \cite{nevo}. The main idea here is Jensen's inequality (and the convexity of the function $t^{2m}$). 
\end{proof}

Given a unitary representation $(\rho, \mathcal{H})$, we define its integrability exponent $q(\rho)$ to be
\begin{gather*}
    \inf \{p : \langle v, \rho(g).v \rangle \in L^p(G) \text{ for $v \in V \subset \mathcal{H}$ a dense subset} \}.
\end{gather*}
In \cite{cowling_haagerup_howe} it is shown that if $q(\rho) \leq 2$, then $(\rho, \mathcal{H})$ is tempered. Assuming $d \geq 3$, i.e. $G$ has $F$-rank at least 2, then by Theorem 5.6 in \cite{gorodnik_nevo} we get that there exists a $p = p(G) < \infty$ such that for every unitary representation $\rho$ not containing any finite-dimensional $G$-subrepresentations, we have $q(\rho) < p$. It follows from H\"{o}lder's inequality that if $\rho$ has integrability exponent $t$, then $\rho^{\otimes n}$ has integrability exponent $t/n$.

\begin{corollary} \label{tensor_power}
Let $(\rho, \mathcal{H})$ be a unitary representation of $M$ with integrability exponent $q(\rho) < \infty$. Let $m$ be an integer such that $q(\rho) < 4m$. Then $\rho^{\otimes 2m}$ is weakly contained in the regular representation.
\end{corollary}

\begin{proof}
The integrability exponenent of $\rho^{\otimes 2 m}$ is $q(\rho)/2m \leq 2$, hence $\rho^{\otimes 2m} \prec \lambda_M$. 
\end{proof}

\begin{lemma}
Suppose $M$ is a KS group. Let $\rho$ be a unitary representation with integrability exponent $q(\rho)$. Let $m$ be an integer such that $q(\rho) < 4m$. Let $1 \leq p < 2$. Suppose $\psi \in L^1(M) \cap L^p(M)$ with $\psi$ non-negative, real-valued, and $\int_M \psi dg = 1$. Then there exists $C_p$ such that
\begin{gather}
||\hat{\psi}(\rho)|| \leq C_p ||\psi||_p^{1/2m}. \nonumber
\end{gather}
\end{lemma}

\begin{proof}
By Corollary \ref{tensor_power}, we have that $\rho^{\otimes 2m}$ is weakly contained in the regular representation. As discussed in Theorem F.4.4 in \cite{property_t_book}, if we have two unitary representaitons $\pi_1$ and $\pi_2$ with $\pi_1 \prec \pi_2$, then for every $f \in L^1(G)$ we have $||\hat{f}(\pi_1)|| \leq ||\hat{f}(\pi_2)||$. Combining this with Lemma \ref{nevo_lemma}, we have that
$||\hat{\psi}(\rho)||^m \leq ||\hat{\psi}(\rho^{\otimes 2m })|| \leq ||\hat{\psi}(\lambda_M)||$.
By the KS property, we get that there exists a $C_p$ depending only on $p$ such that for all $\psi \in L^p(G)$, 
$||\hat{\psi}(\lambda_M)|| \leq C_p^{2m} ||\psi||_p$.
\end{proof}

\begin{theorem} \label{thm_nevo_original}
Suppose $M$ is a KS group. Suppose $(\rho, \mathcal{H})$ is a unitary representation of $M$ with finite integrability exponent. Then there exist $C > 0$ depending only on $M$, and $\theta > 0$ depending only on the integrability exponent, such that for every measurable subset $E \subset M$ with finite positive measure and every $f \in \mathcal{H}$,
\begin{gather}
\Big| \Big|\frac{1}{\vol(E)} \int_E \rho(m).f dm \Big| \Big|_2 \leq C \vol(E)^{-\theta} ||f||_2. \label{eqn_nevo_pf}
\end{gather}
\end{theorem}

\begin{proof}
The function $\psi = \frac{1}{\textnormal{vol}(E)} \mathds{1}_E$ is in $L^1(M) \cap L^p(M)$ for every $p$ and is non-negative, real, and has $L^1$-norm equal to 1. Therefore, if $1 \leq p < 2$ we have
$||\hat{\psi}(\rho)|| \leq C_p ||\psi||_p^{1/2m}$. 
Notice that in fact $\hat{\psi}(\rho).f$ is exactly the expression inside the $|| \cdot ||_2$ on the left hand side of \eqref{eqn_nevo_pf}. 

We also have that $||\psi||_p = \vol(E)^{-1 + 1/p}$.
Therefore,
\begin{gather}
||\hat{\psi}(\rho)|| \leq C_p \vol(E)^{-1/(2m) + 1/(2pm)}. \nonumber
\end{gather}
We see that in fact we can choose $\theta$ to be of the form $-1/(4m) + \delta$ for any $\delta > 0$.
\end{proof}

\begin{proof}[Proof of Proposition \ref{nevo} (Nevo-Style Ergodic Theorem)]
This follows almost immediately from Theorem \ref{thm_nevo_original}. Let $\mathcal{H} \subset L^2(\Gamma \backslash G)$ be the orthogonal complement to the finite-dimensional subrepresentations. Let $a \in L^2(\Gamma \backslash G)$. Notice that if $a$ is $(1, K)$-invariant, then $a \in \mathcal{H}$ if and only if $a$ is orthogonal to the coloring eigenfunctions. This follows from the fact that the coloring eigenfunctions exactly correspond to the finite-dimensional spherical representations in $L^2(\Gamma \backslash G)$ (see Section \ref{sec_coloring_eigenfunctions}), as well as the fact that all finite-dimensional non-spherical representations in $L^2(\Gamma \backslash G)$ are orthogonal to $(1, K)$-invariant functions. 

By Theorem 5.6 in \cite{gorodnik_nevo}, we know that the $G$-action on $\mathcal{H}$ has finite integrability exponent which only depends on $G$ (and not on $\Gamma$). Therefore by Theorem \ref{thm_nevo_original}, we know that there exists a $\theta > 0$ and $C > 0$, both independent of $\Gamma$, such that
\begin{gather}
||\rho^\Gamma_{E}.a||_{L^2(\Gamma \backslash G)} \leq \frac{C}{\vol(E)^{\theta}} ||a||_{L^2(\Gamma \backslash G)}. \nonumber
\end{gather}
\end{proof}

%% file: sections/classification_redo.tex
\sloppy

In this section we classify in a certain sense relative positions of triples of vertices in $\tilde{A}_2$ buildings. This classification allows us to assign coordinates to triples of vertices $(x, y, z)$. Suppose $d_{A^+}(x, y) = (r, s)$ in cone coordinates. Suppose we wish to compute, or at least bound, the volume of the intersection of two polytopal balls centered at $x$ and $y$. Under the coordinatization described in this section, the collection of possible coordinates of points $z$ lying in the intersection of these polytopal balls may be described as the lattice points in a collection of finitely many polytopes whose defining inequalities are derived from the data of $(r, s)$ and the defining inequalities for the polytopes underlying our polytopal balls. In Section \ref{sec_geometric_bound}, we count, for each lattice point in one of these finitely many polytopes parametrizing possible coordinates, the number of points in the intersection of the polytopal balls which correspond to that lattice point. These counting functions will be (approximately) exponential functions whose powers are linear functionals in the coordinates of the lattice point. This reduces computing the volume of the intersection of two polytopal balls to an application of (the possibly degenerate case of) Brion's formula. 

We highly encourage the reader to first read Chapter III.5 in \cite{my_thesis} where the proof method for the geometric bound is explained in the case of the regular tree.

\subsection{Parallelograms}
Suppose $x, y, z \in \Lambda \subset \mathfrak{a}$. We say that $(x, y; z)$ are an \textit{additive triple} if $d_{A^+}(x, z) + d_{A^+}(z, y) = d_{A^+}(x, y)$. We define the \textit{parallelogram} of $x$ and $y$, denoted $\para(x, y)$, to be all $z$ such that $(x, y; z)$  are an additive triple. In case $x = (0, 0)$ and $y = (r, s) \in \mathfrak{a}^+$ in cone coordinates, we clearly have $\para(x, y) = \{(r', s') : (0, 0) \preceq (r', s') \preceq (r, s) \}$.

We may tesselate $\mathfrak{a}$ by equilateral triangles to obtain $\mathcal{X}$, the Coxeter complex associated to $\tilde{A}_2$. The vertices of this simplicial complex correspond to the lattice $\Lambda$. We define a {\it combinatorial path} in $\mathcal{X}$ to be any path between vertices along edges, and we define the length of such paths to be the number of edges traversed in the path. We define a {\it combinatorial geodesic} between two vertices $x$ and $y$ to be any shortest combinatorial path from $x$ to $y$. We define the {\it combinatorial distance} between vertices, $d_c(x, y)$, to be the length of any combinatorial geodesic connecting them.

\begin{proposition} \label{prop_dist_coxeter}
    We have $d_{A^+}(x, y) = (r, s)$ for some $(r, s)$  with $r+s = n$ if and only if $d_c(x, y) = n$. Furthermore, all combinatorial geodesics from $x$ to $y$ lie within $\para(x, y)$.
\end{proposition}
\begin{proof}
    We use induction on $n$. The base case is obvious as $d_{A^+}(x, y) = (1, 0)$ or $(0, 1)$ if and only if $x$ and $y$ are adjacent. For the inductive step, suppose $d_{A^+}(x, y) = (r, s)$ with $r+s = n$. Let $z$ be a neighbor of $y$ such that $z \in \para(x, y)$; such a $z$ must have $d_{A^+}(x, z) \in \{(r-1, s), (s, r-1)\}$. By induction $d_c(x, z) = n-1$ so that $d_c(x, y) \leq n$. However, by induction if $d_c(x, y) < n$, then $r+s < n$ which would be a contradiction. Thus $d_c(x, y) = n$.

    On the other hand suppose $d_c(x, y) = n$ and $d_{A^+}(x, y) = (r, s)$. Suppose $z$ is the penultimate vertex of some combinatorial geodesic from $x$ to $y$. Then $z$ is a neighbor of $y$ such that $d_{A^+}(x, z) \in \{(r-1, s), (r+1, s), (r, s-1), (r, s+1), (r+1, s-1), (r-1, s+1)\}$. Furthermore $d_c(x, z) = n-1$. By the inductive hypothesis $d_{A^+}(x, z) = (r', s')$ with $r' + s' = n-1$ which implies that in fact $(r', s') = (r, s-1)$ or $(r-1, s)$ and thus that $r+s = n$. Such $z$ must then also lie within $\para(x, y)$. By induction all combinatorial geodesics from $x$ to $z$ lie within $\para(x, z)$ (which in turn lies within $\para(x, y)$) and thus the combinatorial geodesic from $x$ to $y$ lies within $\para(x, y)$. 
\end{proof}

Since apartments in the building $\mathcal{B}$ are Coxeter complexes, it makes sense to talk about parallelograms in $\mathcal{B}$.

\begin{proposition} \label{prop_para_building}
    Suppose $x$ and $y$ are two vertices of $\mathcal{B}$. Let $\Sigma$ be any apartment containing $x$ and $y$. Then all combinatorial geodesics in $\mathcal{B}$ from $x$ to $y$ lie in $\para(x, y) \subset \Sigma$, and $\para(x, y)$ is independent of the choice of $\Sigma$.
\end{proposition}
\begin{proof}
    Suppose $d_{A^+}(x, y) = (r, s)$. Suppose $P$ is a combinatorial geodesic from $x$ to $y$. The image of $P$ under the retraction onto $\Sigma$ (based at any choice of chamber in $\Sigma$; see Section \ref{sec_buildings}) is a combinatorial path from $x$ to $y$ which must be at least as short as the shortest path from $x$ to $y$ in $\Sigma$. By Proposition \ref{prop_dist_coxeter} we conclude that the length of $P$ is $r+s$.

    Suppose $P$ does not lie entirely in $\para(x, y) \subset \Sigma$. Then there exist $w, z$ such that $w \not \in \para(x, y)$, $z \in \para(x, y)$, $w$ and $z$ are adjacent, and $w$ is closer to $x$ than $z$. Suppose $d_{A^+}(x, z) = (r', s')$. Then by Proposition \ref{prop_dist_coxeter} we must have $d_{A^+}(z, w) \in \{(r'-1, s'), (r', s'-1)\}$ and there is exactly one such choice, namely the unique such element in $\para(x, z)$. Therefore $w \in \para(x, z) \subseteq \para(x, y)$. 

    If for another choice of apartment $\Sigma'$ we obtain a different parallelogram for $x$ and $y$, then by Proposition \ref{prop_dist_coxeter} we would get a shortest path from $x$ to $y$ which does not lie in $\para(x, y) \subset \Sigma$ which is a contradiction.
\end{proof}

\subsection{Classification of nearly opposite sectors} \label{sec_classify_nearly_opposite}

Suppose $p$ is a vertex in $\mathcal{B}$. Recall that the link of $p$ is itself a spherical building (see Section \ref{sec_buildings}). Hence given chambers $\mathfrak{c}_1$ and $\mathfrak{c}_2$ containing $p$, we may associate an element $d_W(\mathfrak{c}_1, \mathfrak{c}_2) \in \mathfrak{S}_3$ using the Coxeter group-valued metric. Up to simplicial automorphisms of the Coxeter complex of $\mathfrak{S}_3$, there are four relative positions that two chambers in a given apartment in the local spherical building can be in:
\begingroup
\renewcommand\labelenumi{(\theenumi)}
\begin{enumerate}
\item $d_W (\mathfrak{c}_1, \mathfrak{c}_2) = 1$, i.e. $\mathfrak{c}_1 = \mathfrak{c}_2$.
\item $d_W (\mathfrak{c}_1, \mathfrak{c}_2) = (12)$ or $(23)$, i.e. $\mathfrak{c}_1$ and $\mathfrak{c}_2$ are {\it adjacent}.
\item $d_W(\mathfrak{c}_1, \mathfrak{c}_2) = (132)$ or $(1 2 3)$, i.e. $\mathfrak{c}_1$ and $\mathfrak{c}_2$ are {\it nearly opposite}.
\item $d_W(\mathfrak{c}_1, \mathfrak{c}_2) = (1 3)$, i.e. $\mathfrak{c}_1$ and $\mathfrak{c}_2$ are {\it opposite}.
\end{enumerate}
\endgroup

Suppose $x, y \in \mathcal{B}$. Let $\mathfrak{a}^{++}$ be those elements $(r, s) \in \mathfrak{a}^+$ for which $r, s > 0$. If $d_{A^+}(x, y) \in \mathfrak{a}^{++}$, then there is a unique chamber in $\para(x, y)$ containing $x$ which we denote $\mathfrak{c}_{x, y}$.

Now suppose $x, y, p$ are vertices of $\mathcal{B}$. We say that $(x, y; p)$ is a {\it primitive triple} if $\para(x, p) \cap \para(y, p) = \{p\}$. If we have $d_{A^+}(x, p) \in \mathfrak{a}^{++}$ and $d_{A^+}(y, p) \in \mathfrak{a}^{++}$, then any sector $\mathcal{S}$ based at $p$ containing $x$ must have $\mathfrak{c}_1 := \mathfrak{c}_{p, x}$ as its germ, and similarly for $\mathfrak{c}_2 := \mathfrak{c}_{p, y}$. We clearly have $\mathfrak{c}_{1} \cap \mathfrak{c}_2 \subset \para(x, p) \cap \para(y, p) = \{p \}$. Therefore $\mathfrak{c}_{1}$ and $\mathfrak{c}_{2}$ must be either opposite or nearly opposite (Cases (3) and (4) above). 

If instead $(x, y; p)$ is a primitive triple but at least one of $d_{A^+}(p, x)$ and $d_{A^+}(p, y)$ is not in $\mathfrak{a}^{++}$, then we may find sectors $\mathcal{S}_1$ based at $p$ containing $x$, and $\mathcal{S}_2$ based at $p$ containing $y$ such that their germs $\mathfrak{c}_1$ and $\mathfrak{c}_2$ are either opposite or nearly opposite.

We now wish to classify primitive triples. There are a couple of facts about affine buildings that we shall utilize.

\begin{proposition} [\cite{brown} p. 169] \label{add simplex} 
Suppose $\mathfrak{H}$ is a half apartment in $\mathcal{B}$ with boundary wall $\partial \mathfrak{H}$. Suppose $\mathfrak{c}$ is a chamber in $\mathcal{B}$ with a panel lying in $\partial \mathfrak{H}$. Then there exists an apartment in $\mathcal{B}$ containing $\mathfrak{H}$ and $\mathfrak{c}$.
\end{proposition}

\begin{proposition} [\cite{bennett_schwer}, therein referred to as property (EC)] \label{EC}
Suppose $\Sigma_1$ and $\Sigma_2$ are two apartments that intersect in some half apartment $\mathfrak{H}$. Let $\mathfrak{H}_i$ be the other half apartment of $\Sigma_i$. Then $\mathfrak{H}_1$ and $\mathfrak{H}_2$ together form an apartment.
\end{proposition}

\begin{proposition} [\cite{bennett_schwer}, therein referred to as property (CO)] \label{CO}
Suppose $\mathcal{S}_1$ and $\mathcal{S}_2$ are two sectors in $\mathcal{B}$ based at the same point $p$. Suppose the $\mathcal{S}_i$'s determine opposite chambers in the link of $p$ (which is a spherical building). Then $\mathcal{S}_1$ and $\mathcal{S}_2$ are contained in a unique apartment.
\end{proposition}

If $(x, y; p)$ is a primitive triple and $\mathfrak{c}_1 = \mathfrak{c}_{p, x}$ and $\mathfrak{c}_2 = \mathfrak{c}_{p, y}$ are opposite, then by Proposition \ref{CO} one can find an apartment containing all three points. It is then clear that in fact $p \in \para(x, y)$.

We thus are left with understanding the situation when $\mathfrak{c}_1$ and $\mathfrak{c}_2$ are nearly opposite. We extend the terminology and say that two sectors based at the same vertex are {\it nearly opposite} if their germs are.

\subsubsection{Strips and levels of sectors}
Recall that a {\it wall} in the Coxeter complex $\mathcal{X}$ is the line fixed by any reflection in the underlying Coxeter group. Concretely this corresponds to any union of edges in $\mathcal{X}$ which form an infinite line in the Euclidean sense. We define a {\it half-wall} to be any union of edges which form a ray in the Euclidean sense. We define a {\it combinatorial line segment} to be any union of edges which forms a line segment in the Euclidean sense. We define a {\it strip of width $k$} (with $k \geq 1$) in the Coxeter complex $\mathcal{X}$ to be the region bounded between two parallel walls whose combinatorial distance apart is $k$. We define a {\it half-strip} to be the intersection of a strip with any root in $\mathcal{X}$ which is not parallel to the walls defining the strip. Given a half-strip there is a unique chamber which has one of its edges along one of the defining walls, and the other edge along the boundary of the defining root; we call this the {\it germ of the half-strip}.

Given a sector $\mathcal{S}$, and a choice of one of its bounding half-walls $\ell$, we get an associated partition of $\mathcal{S}$ into {\it levels}, each of which is a half-strip of width one, which we shall index by $k \in \mathbb{N}$ (starting from $k = 1$). Let $\textnormal{Level}(\mathcal{S}, k)$ denote the $k$th level, and let $\mathcal{S}^{(k)}$ denote the sector obtained by removing the first $k$ levels (i.e. $\mathcal{S}^{(k)}$ is the union of levels $k+1, k+2, \dots$). See Figure \ref{fig_strip_levels}.

\begin{figure}[!h]
\centering
\begin{minipage}{.5\textwidth}
  \centering
  \subfile{../images/levels}
  \caption{}
  \label{fig_strip_levels}
\end{minipage}%
\begin{minipage}{.5\textwidth}
  \centering
  \subfile{../images/strip}
  \caption{}
  \label{fig_convex_hull}
\end{minipage}
\end{figure}

Suppose $\Sigma$ is some apartment in $\mathcal{B}$. Suppose $\mathfrak{c} \subset \Sigma$ is a chamber and $\ell \subset \Sigma$ is a combinatorial line segment whose edges form part of a wall in $\Sigma$. Suppose the endpoints of $\ell$ are the vertices $v_1$ and $v_2$. Suppose the edge of $\ell$ containing $v_1$ is also one of the edges of $\mathfrak{c}$. Let $w$ be the opposite vertex in $\mathfrak{c}$ to this edge. Then the convex hull of $\mathfrak{c}$ and $\ell$ is clearly $\para(w, v_2) \cup \mathfrak{c}$. See Figure \ref{fig_convex_hull}. 

Repeated application of this observation gives the following:

\begin{lemma} \label{lemma_convex_hull}
Suppose $\Sigma$ is an apartment in $\mathcal{B}$, and $\mathfrak{c} \subset \Sigma$ is a chamber, and $\ell \subset \Sigma$ is a half-wall whose initial edge is one of the edges of $\mathfrak{c}$. Then the combinatorial convex hull of $\mathfrak{c}$ and $\ell$ is the half-strip in $\Sigma$ of width one which has one of its bounding half-walls equal to $\ell$ and has germ equal to $\mathfrak{c}$. If instead $\ell$ is a wall, one of whose edges is one of the edges of $\mathfrak{c}$, then the convex hull of $\mathfrak{c}$ and $\ell$ is the strip of width one containing $\mathfrak{c}$ and which has $\ell$ as one of its bounding walls.
\end{lemma}

\begin{corollary} \label{cor_convex_hull}
Suppose $\mathfrak{H}$ is a half-apartment in $\mathcal{B}$. Suppose $\mathfrak{c}$ is a chamber not lying in $\mathfrak{H}$ but which shares a panel with $\partial \mathfrak{H}$. Then all apartments containing $\mathfrak{H}$ and $\mathfrak{c}$ (such apartments exist by Proposition \ref{add simplex}) must contain the strip determined by $\partial \mathfrak{H}$ and $\mathfrak{c}$ as in Lemma \ref{lemma_convex_hull}.
\end{corollary}

\subsubsection{Classification of nearly opposite sectors}
Let $\mathcal{S}_1$ and $\mathcal{S}_2$ be nearly opposite sectors with common vertex $p$. Let $\mathcal{P}$ be the local spherical building obtained from the link of $p$. Let $a_1, b_1$ be the boundary half-walls of $\mathcal{S}_1$ and $a_2, b_2$ the boundary half-walls of $\mathcal{S}_2$ with $a_1$ and $a_2$ determining opposite vertices in the link of $p$. Let $\mathfrak{s}_1$ and $\mathfrak{s}_2$ be the germs of $\mathcal{S}_1$ and $\mathcal{S}_2$ respectively. 

We claim that the union of $a_1$ and $a_2$ forms an infinite (Euclidean) geodesic (and hence there exists some apartment containing $a_1 \cup a_2$ in which this set corresponds to a wall). Let $\bar{a}_1$ and $\bar{a}_2$ be the neighbors of $p$ in $a_1$ and $a_2$ respectively. Notice that $a_1$ is a geodesic, $a_2$ is a geodesic, and the path $\bar{a}_1 \to p \to \bar{a}_2$ is a geodesic. Hence $a_1 \cup a_2$ is locally geodesic and therefore globally geodesic by the CAT(0) geometry. 

In the sequel, any time we talk about levels of the sectors $\mathcal{S}_j$, it is with respect to the half-wall $a_j$.

\begin{lemma} \label{lemma_add_level}
Suppose $\mathcal{S}_1$ and $\mathcal{S}_2$ are as above. Then there exists a half-apartment $\mathfrak{H}$ such that $\partial \mathfrak{H} = a_1 \cup a_2$, and $\mathfrak{H}$ does not otherwise intersect $\mathcal{S}_1$ or $\mathcal{S}_2$. Furthermore $\mathfrak{H}$ can be extended to a half-apartment $\mathfrak{H}_1$ by appending $\textnormal{Level}(\mathcal{S}_1, 1)$, $\textnormal{Level}(\mathcal{S}_2, 1)$, and a uniquely determined chamber $\mathfrak{t}$. 
\end{lemma}

\begin{proof}
Because the link of $p$ is a spherical building, we know that there is some apartment $\Pi$ in the local spherical building $\mathcal{P}$ containing $\mathfrak{s}_1$ and $\mathfrak{s}_2$. Let $\mathfrak{t}$ be the chamber connecting $\mathfrak{s}_1$ and $\mathfrak{s}_2$ along the minimal gallery in $\Pi$ from $\mathfrak{s}_1$ to $\mathfrak{s}_2$. Then $\mathfrak{t}$ must have $p$ as a vertex as well as the neighbors of $p$ in $b_1$ and $b_2$. Since any three vertices determine at most one chamber, $\mathfrak{t}$ is well-defined independently of the choice of $\Pi$. Notice also that $\{\mathfrak{s}_1, \mathfrak{t}, \mathfrak{s}_2\}$ form a half-apartment in $\mathcal{P}$, and therefore by \cite{ronan} Lemma 6.3, any chamber adjacent to $\mathfrak{s}_2$ which has $\mathfrak{s}_2 \cap a_2$ as one of its panels must be opposite to $\mathfrak{s}_1$. 

Let $\Sigma$ be any apartment in $\mathcal{B}$ containing $\mathcal{S}_2$. Let $\mathcal{R}$ be the sector based at $p$ in $\Sigma$ which shares the half-wall $a_2$ with $\mathcal{S}_2$. Let $\mathfrak{r}$ be the germ of $\mathcal{R}$. Because $\mathfrak{r}$ shares a panel with $\mathfrak{s}_2$, we conclude that it is opposite to $\mathfrak{s}_1$ in $\mathcal{P}$. Therefore by Proposition \ref{CO}, there exists a unique apartment $\Sigma'$ containing $\mathcal{S}_1$ and $\mathcal{R}$. Furthermore the half-apartment $\mathfrak{H}$ in $\Sigma'$ containing both $a_1 \cup a_2$ and $\mathcal{R}$ clearly does not intersect $\mathcal{S}_1$ other than along $a_1$. We claim that it also cannot intersect $\mathcal{S}_2$ other than along $a_2$. This shall follow from the following in which we show that we may append a level to $\mathfrak{H}$ to obtain a bigger half-apartment as in the statement of the lemma.

Consider the convex hull of $\mathfrak{s}_1$ and $\partial \mathfrak{H}$. By Proposition \ref{add simplex}, this adds on a strip to $\mathfrak{H}$. This strip clearly contains $\textnormal{Level}(\mathcal{S}_1, 1)$ by Corollary \ref{cor_convex_hull}. Furthermore the convex hull of $\mathfrak{s}_1$ and $\mathfrak{r}$ must contain $\mathfrak{t}$ and $\mathfrak{s}_2$ (because they are opposite, they determine a unique apartment in $\mathcal{P}$ by \cite{ronan} Chapter 6.1). Since the convex hull of $\mathfrak{s}_2$ and $a_2$ is $\textnormal{Level}(\mathcal{S}_2, 1)$, we conclude the proof of the lemma. See Figure \ref{fig_sector_configuration}.
\end{proof}

We define an {\it equilateral triangle} in $\mathcal{B}$ to be the intersection of a sector $\mathcal{S}$ with any half-apartment $\mathfrak{J}$ which has an extension to an apartment containing $\mathcal{S}$ and such that $\partial \mathfrak{J}$ is transverse (i.e. not parallel) to the half-walls of $\mathcal{S}$. Such a region lies in an apartment, and in that apartment it is a Euclidean equilateral triangle. The size of the equilateral triangle is the length of any one of its sides.

\begin{figure}[!h]
\centering
\begin{minipage}{.5\textwidth}
  \centering
  \subfile{../images/sector_configuration}
  \caption{}
  \label{fig_sector_configuration}
\end{minipage}%
\begin{minipage}{.5\textwidth}
  \centering
  \subfile{../images/apartment_perspective}
  \caption{}
  \label{fig_triples_flat}
\end{minipage}
\end{figure}

\begin{lemma} \label{lemma_classify_triples}
Suppose $\mathcal{S}_1$ and $\mathcal{S}_2$ are as before. Let $\mathcal{S}_1^{(n)}$ and $\mathcal{S}_2^{(n)}$ be the sectors obtained by removing the first $n$ levels. Let $p_1^{(n)}$ and $p_2^{(n)}$ be the base vertices of these sectors. Let $a_1^{(n)}$ and $a_2^{(n)}$ be the bouding half-walls of these sectors parallel to $a_1$ and $a_2$. Then either $\mathcal{S}_1$ and $\mathcal{S}_2$ are contained in an apartment, or there exists a $k$ such that all of the following hold: 
\begingroup
\renewcommand\labelenumi{(\theenumi)}
\begin{enumerate}
\item $\mathcal{S}_1^{(k)}$ and $\mathcal{S}_2^{(k)}$ are contained in an apartment $\Sigma$,
\item $p_1^{(k)}$ and $p_2^{(k)}$ are the endpoints of a line segment of length $k$ in $\Sigma$,
\item $\ell \cup a_1^{(k)} \cup a_2^{(k)}$ is a wall in $\Sigma$,
\item given $x \in \mathcal{S}_1^{(k)}$ and $y \in \mathcal{S}_2^{(k)}$, we have $\ell \subset \para(x, y) \subset \Sigma$,
\item and $p_1^{(k)}$, $p_2^{(k)}$, and $p$ form the vertices of an equilateral triangle of size $k$ and one of whose sides is $\ell$.
\end{enumerate}
\endgroup
\end{lemma}

\begin{figure}[!h]
\centering
\subfile{../images/perspective}
\caption{} \label{fig_triples_perspective}
\end{figure}

\begin{proof}
The strategy is to attempt to add levels to the half-apartment $\mathfrak{H}$ from Lemma \ref{lemma_add_level}. In fact Lemma \ref{lemma_add_level} already tells us that we can add the first level. Let $\mathfrak{H}_k$ denote the half apartment obtained by appending the first $k$ levels of $\mathcal{S}_1$ and $\mathcal{S}_2$ to $\mathfrak{H}$, assuming we are able to do so. We thus know that $\mathfrak{H}_1$ exists.
If we are able to keep adding levels indefinitely, then, by taking the union of these levels together with $\mathfrak{H}$, we obtain a apartment containing $\mathcal{S}_1$ and $\mathcal{S}_2$.

Suppose instead that $n$ is such that we have successfully added $n$ levels (and hence constructed $\mathfrak{H}_n$) but are not able to add the $(n+1)$st level. By Proposition \ref{CO}, we may find apartments $\Sigma_1$ and $\Sigma_2$ such that $\Sigma_j$ contains $\mathfrak{H}_n$ and $\mathcal{S}_j^{(n)}$. We clearly have $\Sigma_1 \cap \Sigma_2 \supset \mathfrak{H}_n$. Suppose some other chamber $\mathfrak{c}$ were in the intersection. Let $\mathfrak{c}'$ be the last chamber not in $\mathfrak{H}_n$ along any minimal gallery from $\mathfrak{c}$ to any chamber in $\mathfrak{H}_n$. This gallery must lie in the intersection $\Sigma_1 \cap \Sigma_2$ as the intersection is combinatorially convex. Therefore $\mathfrak{c}'$ is a chamber in the intersection which is not contained in $\mathfrak{H}_n$ but shares a panel with $\partial \mathfrak{H}_n$. By Corollary \ref{cor_convex_hull}, any apartment containing $\mathfrak{H}_n$ and $\mathfrak{c}'$ must also contain the strip parallel to $\partial \mathfrak{H}_n$ containing $\mathfrak{c}'$ and hence this strip is in $\Sigma_1 \cap \Sigma_2$. However, this strip clearly contains $\textnormal{Level}(\mathcal{S}_j, n+1)$ for $j = 1, 2$, so we would be able to form $\mathfrak{H}_{n+1}$, which is a contradiction. Therefore, $\Sigma_1 \cap \Sigma_2 = \mathfrak{H}_n$. 

Let $\mathcal{J}_1$ and $\mathcal{J}_2$ be the other half apartments of $\Sigma_1$ and $\Sigma_2$. By Proposition \ref{EC}, we may form an apartment $\Sigma = \mathcal{J}_1 \cup \mathcal{J}_2$. It is clear from the construction that $\Sigma$ satisfies Properties (1), (2), (3), and (5) (with $n = k$). Property (4) follows from the observation that in $\Sigma$, the sectors $\mathcal{S}_1^{(n)}$ and $\mathcal{S}_2^{(n)}$ have opposite orientations.
\end{proof}

\begin{remark}
    Figures \ref{fig_triples_flat} and \ref{fig_triples_perspective} illustrate the geometry of Lemma \ref{lemma_classify_triples}. Given vertices $x, y$, and $p$, we can consider $\para(x, y)$, which is the black parallelogram in Figure \ref{fig_triples_perspective}. If $\para(x, p)$ and $\para(y, p)$ only intersect at $p$, and the associated chambers $\mathfrak{c}_{p, x}$ and $\mathfrak{c}_{p, y}$ are nearly opposite, then we may find nearly opposite sectors at $p$, one of which contains $x$ (the blue sector) and the other of which contains $y$ (the red sector). Certain half-strips of these sectors may be combined with an equilateral triangle (the green triangle) to form a strip. The remaining part of these sectors may then we placed in some apartment (containing $x$ and $y$) which is the brown apartment in the figure. Figure \ref{fig_triples_flat} shows essentially the same picture as Figure \ref{fig_triples_perspective} from the ``perspective'' of the brown apartment. The blue sector is exactly $\mathcal{S}_1^{(k)}$ and the red sector is exactly $\mathcal{S}_2^{(k)}$ as in Lemma \ref{lemma_classify_triples} (here $k = 3$).
\end{remark}

\subsection{Algebraic interpretation of the classification of nearly opposite sectors} \label{sec_algebraic_nearly_opposite}
We now explain how to interpret Lemma \ref{lemma_classify_triples} as a statement about the relationship between the Bruhat decomposition over the residue field vs. over the underlying local field. In this section we assume some familiarity with the spherical building at infinity associated to an affine building. Further information about this can be found, for example, in Chapter 9.3 of \cite{ronan}.

Suppose we have two sectors, $\mathcal{S}_1$ and $\mathcal{S}_2$ both based at the point $1 K$ in $G/K$ ($G = \textnormal{PGL}(3, F)$). Without loss of generality, we can assume that $\mathcal{S}_1$ is the ``standard sector'' whose associated chamber in the spherical building at infinity is stabilized by the standard Borel $B$ consisting of upper triangular matrices. Some element $\kappa \in K$ must take $\mathcal{S}_1$ to $\mathcal{S}_2$.

We have a natural map $\theta: K \to \textnormal{PGL}(3, \mathbb{F}_q)$ given by choosing any representative in $\textnormal{GL}(3, \mathcal{O})$ and taking all of its entries modulo $\varpi \mathcal{O}$. Let $I$ be the preimage of the standard Borel subgroup $\tilde{B}$ of $\textnormal{PGL}(3, \mathbb{F}_q)$; the group $I$ is often referred to as the Iwahori subgroup. By a slight refinement of the Iwahori decomposition, we have a decomposition:
\begin{gather*}
    K = \bigsqcup_{w \in \mathfrak{S}_3} I w (B \cap K)
\end{gather*}
where we interpret $w \in \mathfrak{S}_3$ as a permutation matrix. If we write $\kappa = i \cdot  w \cdot b$ with respect to this decomposition, then we must have that $\theta(\kappa) \in \tilde{B} w \tilde{B}$. This $w$ exactly tells us the relative position of the germs of the sectors $\mathcal{S}_1$ and $\mathcal{S}_2$ in the local spherical building. By assumption these sectors are nearly opposite, and thus, without loss of generality, we can take
\begin{gather*}
    w = \begin{bmatrix}
        0 & 1 & 0 \\
        0 & 0 & 1 \\
        1 & 0 & 0
    \end{bmatrix}.
\end{gather*}

Elements in $I$ can in turn be uniquely expressed as $i = u \cdot d \cdot \ell$ where $u$ is a strictly upper triangular matrix with all entries in $\mathcal{O}$, $d$ is a diagonal matrix with diagonal entries in $\mathcal{O}^\times$, and $\ell$ is a lower triangular matrix with all off diagonal entries in $\varpi \mathcal{O}$ (see for instance Proposition 7.1 and 11.2 in \cite{casselman_gl_n}). The element $\ell$ can in turn be uniquely expressed as 
\begin{gather*}
    \ell = \begin{bmatrix}
        1 & 0 & 0 \\
        \alpha & 1 & 0 \\
        0 & 0 & 1
    \end{bmatrix}
    \begin{bmatrix}
        1 & 0 & 0 \\
        0 & 1 & 0 \\
        \beta & \gamma & 1
    \end{bmatrix}
\end{gather*}
for some $\alpha, \beta, \gamma \in \varpi \mathcal{O}$. Notice that
\begin{gather*}
    w^{-1} \begin{bmatrix}
        1 & 0 & 0 \\
        0 & 1 & 0 \\
        \beta & \gamma & 1
    \end{bmatrix} w = \begin{bmatrix}
        1 & \beta & \gamma \\
        0 & 1 & 0 \\
        0 & 0 & 1
    \end{bmatrix} \in B.
\end{gather*}
Therefore
\begin{gather*}
\kappa = u \cdot d \cdot \begin{bmatrix}
        1 & 0 & 0 \\
        \alpha & 1 & 0 \\
        0 & 0 & 1
    \end{bmatrix} w \Bigg(w^{-1} \begin{bmatrix}
        1 & 0 & 0 \\
        0 & 1 & 0 \\
        \beta & \gamma & 1
    \end{bmatrix} w \Bigg) b = b' \cdot \begin{bmatrix}
        1 & 0 & 0 \\
        \alpha & 1 & 0 \\
        0 & 0 & 1
    \end{bmatrix} w \cdot b''
\end{gather*}
where $b', b'' \in B$. 

There are now two cases. In the first case $\alpha = 0$ and thus $\kappa \in B w B$. This implies that $\mathcal{S}_1$ and $\mathcal{S}_2$ define nearly opposite chambers in the spherical building at infinity which in turn implies that $\mathcal{S}_1$ and $\mathcal{S}_2$ lie in some apartment and are nearly opposite in that apartment. 

In the second case we have $\alpha \neq 0$. We then notice that
\begin{gather*}
\begin{bmatrix}
        1 & 0 & 0 \\
        \alpha & 1 & 0 \\
        0 & 0 & 1
    \end{bmatrix} w = \begin{bmatrix}
        0 & 1 & 0 \\
        0 & \alpha & 1 \\
        1 & 0 & 0
    \end{bmatrix} = \begin{bmatrix}
        1 & \alpha^{-1} & 0 \\
        0 & 1 & 0 \\
        0 & 0 & 1
    \end{bmatrix} \begin{bmatrix}
        0 & 0 & 1 \\
        0 & 1 & 0 \\
        1 & 0 & 0
    \end{bmatrix} 
    \begin{bmatrix}
        1 & 0 & 0 \\
        0 & \alpha & 1 \\
        0 & 0 & \alpha^{-1}
    \end{bmatrix}.
\end{gather*}
This implies that $\kappa \in B \tilde{w} B$ where $\tilde{w}$ is the longest element in $\mathfrak{S}_3$. This in turn implies that $\mathcal{S}_1$ and $\mathcal{S}_2$ define opposite chambers in the spherical building at infinity, which further implies that there is a unique apartment in the spherical building at infinity containing both chambers. In terms of the affine building, this implies that there is a unique apartment in the affine building which contains subsectors of $\mathcal{S}_1$ and $\mathcal{S}_2$, and in this apartment these subsectors are oriented opposite. These are precisely the subsectors $\mathcal{S}_1^{(k)}$ and $\mathcal{S}_2^{(k)}$ as in Lemma \ref{lemma_classify_triples}. The specific value of $k$ as in the lemma corresponds to the valuation of $\alpha$.

\subsection{Classification of primitive triples} \label{sec_classify_triples}

\subsubsection{Existence of confluence points}

\begin{lemma} \label{lemma_exist_confluence}
Suppose $x, y$, and $p$ are vertices of $\mathcal{B}$. Let $\mathcal{D} = \para(x, p) \cap \para(y, p)$. Then there exist points $z \in \mathcal{D}$ which simultaneously minimize $d_c(x, \cdot)$ and $d_c(y, \cdot)$ over all points in $\mathcal{D}$. Furthermore, for any such $z$, $(x, y; z)$ forms a primitive triple.
\end{lemma}

\begin{proof}
Let $z \in \mathcal{D}$. By Proposition \ref{prop_para_building}, this implies that there exist combinatorial geodesics from $x$ to $p$ and from $y$ to $p$ which both pass through $z$, and hence
$d_c(x, p) = d_c(x, z) + d_c(z, p)$ and $d_c(y, p) = d_c(y, z) + d_c(z, p)$.

Now suppose $z \in \mathcal{D}$ is as close to $x$ with respect to $d_c(\cdot, \cdot)$ as any other point in $\mathcal{D}$. Suppose $w \in \mathcal{D}$ is some other point. We wish to show that $d_c(y, z) \leq d_c(y, w)$. By the additive triple property we have
$d_c(x, z) + d_c(z, p) = d_c(x, w) + d_c(w, p)$ and $d_c(y, z) + d_c(z, p) = d_c(y, w) + d_c(w, p)$, which imply that $d_c(x, z) - d_c(y, z) = d_c(x, w) - d_c(y, w)$ and $d_c(x, z) - d_c(x, w) = d_c(y, z) - d_c(y, w)$.
By assumption $d_c(x, z) - d_c(x, w) \geq 0$. Therefore, the same holds for $d_c(y, z) - d_c(y, w)$, and therefore $z$ is as close to $x$ and $y$ as any other point in $\mathcal{D}$. 

We now show that $(x, y; z)$ is a primitive triple. Suppose $w \in \para(x, z) \cap \para(y, z)$. Then, on the one hand $w \in \mathcal{D}$, and on the other hand by Proposition \ref{prop_para_building}, $w$ lies along some combinatorial geodesic from $x$ to $z$ and hence, unless $w = z$, $w$ is closer to $x$ than $z$ is. But $z$ is as close to $x$ as any other point in $\mathcal{D}$. Therefore $w = z$.
\end{proof}

Given $x, y$, and $p$, we call any point $z \in \para(x, p) \cap \para(y, p)$ satisfying the conditions in Lemma \ref{lemma_exist_confluence} a {\it confluence point} of $(x, y; p)$. 

\begin{remark} \label{remark_many_confluence_points}
In general $(x, y; p)$ may have several confluence points. For example suppose $d_{A^+}(x,  y) = (1, 1)$. Consider the unique edge $e \in \para(x, y)$ which passes through the interior of $\para(x, y)$. Let $p$ be any point such that $e \cup p$ forms a chamber $\mathfrak{c} \subset \mathcal{B}$. Then $\para(x, p) \cap \para(y, p) = \mathfrak{c}$. Therefore, either endpoint of $e$ may be considered a confluence point of $(x, y; p)$.
\end{remark}

\subsubsection{Directions of the bounding line segments of a parallelogram}
Suppose $d_{A^+}(x, y) = (r, s)$. The parallelogram $\para(x, y)$ viewed as a Euclidean parallelogram has as its corners $x, y$ as well as two other vertices $w, z$. One of these, say $w$, satisfies $d_{A^+}(x, w) = (r, 0)$, and the other, say $z$, satisfies $d_{A^+}(x, w) = (0, s)$. In such a case we say that the line from $x$ to $w$ is {\it in the $(1, 0)$-direction with respect to $x$}, and the edge from $x$ to $z$ is {\it in the $(0, 1)$-direction with respect to $x$}. Combinatorial line segments inside of $\para(x, y)$ which are parallel to, e.g., the combinatorial line segment from $x$ to $w$ are also said to be in the $(1, 0)$-direction with respect to $x$. 


\subsubsection{Branch lines and statement of the classification}
The following is a summary of the content of Section \ref{sec_classify_nearly_opposite} applied towards the classification of primitive triples (together with some straightforward calculations which have been suppressed).

\begin{lemma} \label{lemma_actual_classification}
Suppose $(x, y; p)$ is a primitive triple. Let $\mathcal{S}_x$ (and $\mathcal{S}_y$, resp.) be any sector based at $p$ containing $x$ (containing $y$, resp.). 
\begin{enumerate}
\item If $\mathcal{S}_x$ and $\mathcal{S}_y$ are opposite, then $p \in \para(x, y)$. We define $\ell = p$ to be the branch line of the primitive triple (we may consider a point to be a line of length 0).
\item Suppose $\mathcal{S}_x$ and $\mathcal{S}_y$ are nearly opposite.
\begin{enumerate}
\item Suppose there exists an apartment $\Sigma$ containing $\mathcal{S}_x$ and $\mathcal{S}_y$. Let $\mathcal{S}'$ be the sector based at $p$ which connects these two sectors in $\Sigma$. Let $\ell$ be all points in $\mathcal{S}' \cap \para(x, y)$ minimizing $d_c(p, \cdot)$. Then $\ell$ consists of a line segment; we define $\ell$ to be the branch line of the primitive triple.
\item Suppose there does not exist an apartment containing $\mathcal{S}_x$ and $\mathcal{S}_y$. Let $\ell$ be as in Lemma \ref{lemma_classify_triples}. Then $\ell$ is defined to be the branch line of the primitive triple.
\end{enumerate}
\end{enumerate}
\end{lemma}

%% file: images/levels.tex
\begin{tikzpicture}[scale=0.5]
\fill[color = blue, opacity = 0.3]  (0.000000, 0.000000) --  (0.500000, -0.866025) --  (10.000000, -0.866025) --  (10.000000, 0.000000) -- cycle;
\fill[color = red, opacity = 0.3]  (0.500000, -0.866025) --  (1.000000, -1.732051) --  (10.000000, -1.732051) --  (10.000000, -0.866025) -- cycle;
\fill[color = blue, opacity = 0.3]  (1.000000, -1.732051) --  (1.500000, -2.598076) --  (10.000000, -2.598076) --  (10.000000, -1.732051) -- cycle;
\fill[color = red, opacity = 0.3]  (1.500000, -2.598076) --  (2.000000, -3.464102) --  (10.000000, -3.464102) --  (10.000000, -2.598076) -- cycle;
\fill[color = blue, opacity = 0.3]  (2.000000, -3.464102) --  (2.500000, -4.330127) --  (10.000000, -4.330127) --  (10.000000, -3.464102) -- cycle;
\fill[color = red, opacity = 0.3]  (2.500000, -4.330127) --  (3.000000, -5.196152) --  (10.000000, -5.196152) --  (10.000000, -4.330127) -- cycle;
\fill[color = blue, opacity = 0.3]  (3.000000, -5.196152) --  (3.500000, -6.062178) --  (10.000000, -6.062178) --  (10.000000, -5.196152) -- cycle;
\draw[color = black, opacity = 0.1] (0.000000, 0.000000) -- (10.000000, 0.000000);
\draw[color = black, opacity = 0.1] (0.500000, 0.866025) -- (10.000000, 0.866025);
\draw[color = black, opacity = 0.1] (0.000000, 0.000000) -- (0.577350, 1.000000);
\draw[color = black, opacity = 0.1] (0.000000, 0.000000) -- (0.000000, 0.000000);
\draw[color = black, opacity = 0.1] (1.000000, 0.000000) -- (0.500000, 0.866025);
\draw[color = black, opacity = 0.1] (2.000000, 0.000000) -- (1.422650, 1.000000);
\draw[color = black, opacity = 0.1] (3.000000, 0.000000) -- (2.422650, 1.000000);
\draw[color = black, opacity = 0.1] (4.000000, 0.000000) -- (3.422650, 1.000000);
\draw[color = black, opacity = 0.1] (5.000000, 0.000000) -- (4.422650, 1.000000);
\draw[color = black, opacity = 0.1] (6.000000, 0.000000) -- (5.422650, 1.000000);
\draw[color = black, opacity = 0.1] (7.000000, 0.000000) -- (6.422650, 1.000000);
\draw[color = black, opacity = 0.1] (8.000000, 0.000000) -- (7.422650, 1.000000);
\draw[color = black, opacity = 0.1] (9.000000, 0.000000) -- (8.422650, 1.000000);
\draw[color = black, opacity = 0.1] (10.000000, 0.000000) -- (9.422650, 1.000000);
\draw[color = black, opacity = 0.1] (1.000000, 0.000000) -- (1.577350, 1.000000);
\draw[color = black, opacity = 0.1] (2.000000, 0.000000) -- (2.577350, 1.000000);
\draw[color = black, opacity = 0.1] (3.000000, 0.000000) -- (3.577350, 1.000000);
\draw[color = black, opacity = 0.1] (4.000000, 0.000000) -- (4.577350, 1.000000);
\draw[color = black, opacity = 0.1] (5.000000, 0.000000) -- (5.577350, 1.000000);
\draw[color = black, opacity = 0.1] (6.000000, 0.000000) -- (6.577350, 1.000000);
\draw[color = black, opacity = 0.1] (7.000000, 0.000000) -- (7.577350, 1.000000);
\draw[color = black, opacity = 0.1] (8.000000, 0.000000) -- (8.577350, 1.000000);
\draw[color = black, opacity = 0.1] (9.000000, 0.000000) -- (9.577350, 1.000000);
\draw[color = black, opacity = 0.1] (0.000000, 0.000000) -- (0.000000, 0.000000);
\draw[color = black, opacity = 0.1] (-0.500000, 0.866025) -- (0.500000, 0.866025);
\draw[color = black, opacity = 0.1] (0.000000, 0.000000) -- (0.577350, 1.000000);
\draw[color = black, opacity = 0.1] (-0.500000, 0.866025) -- (-0.422650, 1.000000);
\draw[color = black, opacity = 0.1] (0.000000, 0.000000) -- (-0.577350, 1.000000);
\draw[color = black, opacity = 0.1] (0.500000, 0.866025) -- (0.422650, 1.000000);
\draw[color = black, opacity = 0.1] (-1.000000, 0.000000) -- (0.000000, 0.000000);
\draw[color = black, opacity = 0.1] (-1.000000, 0.866025) -- (-0.500000, 0.866025);
\draw[color = black, opacity = 0.1] (0.000000, 0.000000) -- (0.000000, 0.000000);
\draw[color = black, opacity = 0.1] (-1.000000, 0.000000) -- (-0.500000, 0.866025);
\draw[color = black, opacity = 0.1] (0.000000, 0.000000) -- (-0.577350, 1.000000);
\draw[color = black, opacity = 0.1] (-1.000000, 0.000000) -- (-1.000000, 0.000000);
\draw[color = black, opacity = 0.1] (-1.000000, 0.000000) -- (0.000000, 0.000000);
\draw[color = black, opacity = 0.1] (-1.000000, -1.732051) -- (0.000000, 0.000000);
\draw[color = black, opacity = 0.1] (-1.000000, 0.000000) -- (-1.000000, 0.000000);
\draw[color = black, opacity = 0.1] (0.000000, 0.000000) -- (0.000000, 0.000000);
\draw[color = black, opacity = 0.1] (-1.000000, -0.866025) -- (-0.500000, -0.866025);
\draw[color = black, opacity = 0.1] (-0.500000, -0.866025) -- (-1.000000, 0.000000);
\draw[color = black, opacity = 0.1] (0.000000, 0.000000) -- (0.000000, 0.000000);
\draw[color = black, opacity = 0.1] (-1.000000, -1.732051) -- (0.000000, 0.000000);
\draw[color = black, opacity = 0.1] (3.464102, -6.000000) -- (0.000000, 0.000000);
\draw[color = black, opacity = 0.1] (-0.500000, -0.866025) -- (0.500000, -0.866025);
\draw[color = black, opacity = 0.1] (-1.000000, -1.732051) -- (1.000000, -1.732051);
\draw[color = black, opacity = 0.1] (-1.000000, -2.598076) -- (1.500000, -2.598076);
\draw[color = black, opacity = 0.1] (-1.000000, -3.464102) -- (2.000000, -3.464102);
\draw[color = black, opacity = 0.1] (-1.000000, -4.330127) -- (2.500000, -4.330127);
\draw[color = black, opacity = 0.1] (-1.000000, -5.196152) -- (3.000000, -5.196152);
\draw[color = black, opacity = 0.1] (-1.000000, -3.464102) -- (0.500000, -0.866025);
\draw[color = black, opacity = 0.1] (-1.000000, -5.196152) -- (1.000000, -1.732051);
\draw[color = black, opacity = 0.1] (-0.464102, -6.000000) -- (1.500000, -2.598076);
\draw[color = black, opacity = 0.1] (0.535898, -6.000000) -- (2.000000, -3.464102);
\draw[color = black, opacity = 0.1] (1.535898, -6.000000) -- (2.500000, -4.330127);
\draw[color = black, opacity = 0.1] (2.535898, -6.000000) -- (3.000000, -5.196152);
\draw[color = black, opacity = 0.1] (2.464102, -6.000000) -- (-0.500000, -0.866025);
\draw[color = black, opacity = 0.1] (1.464102, -6.000000) -- (-1.000000, -1.732051);
\draw[color = black, opacity = 0.1] (0.464102, -6.000000) -- (-1.000000, -3.464102);
\draw[color = black, opacity = 0.1] (-0.535898, -6.000000) -- (-1.000000, -5.196152);
\draw[color = black, opacity = 0.1] (0.000000, 0.000000) -- (10.000000, 0.000000);
\draw[color = black, opacity = 0.1] (0.000000, 0.000000) -- (0.000000, 0.000000);
\draw[color = black, opacity = 0.1] (3.464102, -6.000000) -- (0.000000, 0.000000);
\draw[color = black, opacity = 0.1] (4.464102, -6.000000) -- (1.000000, 0.000000);
\draw[color = black, opacity = 0.1] (5.464102, -6.000000) -- (2.000000, 0.000000);
\draw[color = black, opacity = 0.1] (6.464102, -6.000000) -- (3.000000, 0.000000);
\draw[color = black, opacity = 0.1] (7.464102, -6.000000) -- (4.000000, 0.000000);
\draw[color = black, opacity = 0.1] (8.464102, -6.000000) -- (5.000000, 0.000000);
\draw[color = black, opacity = 0.1] (9.464102, -6.000000) -- (6.000000, 0.000000);
\draw[color = black, opacity = 0.1] (10.000000, -5.196152) -- (7.000000, 0.000000);
\draw[color = black, opacity = 0.1] (10.000000, -3.464102) -- (8.000000, 0.000000);
\draw[color = black, opacity = 0.1] (10.000000, -1.732051) -- (9.000000, 0.000000);
\draw[color = black, opacity = 0.1] (10.000000, 0.000000) -- (10.000000, 0.000000);
\draw[color = black, opacity = 0.1] (0.500000, -0.866025) -- (10.000000, -0.866025);
\draw[color = black, opacity = 0.1] (1.000000, -1.732051) -- (10.000000, -1.732051);
\draw[color = black, opacity = 0.1] (1.500000, -2.598076) -- (10.000000, -2.598076);
\draw[color = black, opacity = 0.1] (2.000000, -3.464102) -- (10.000000, -3.464102);
\draw[color = black, opacity = 0.1] (2.500000, -4.330127) -- (10.000000, -4.330127);
\draw[color = black, opacity = 0.1] (3.000000, -5.196152) -- (10.000000, -5.196152);
\draw[color = black, opacity = 0.1] (0.500000, -0.866025) -- (1.000000, 0.000000);
\draw[color = black, opacity = 0.1] (1.000000, -1.732051) -- (2.000000, 0.000000);
\draw[color = black, opacity = 0.1] (1.500000, -2.598076) -- (3.000000, 0.000000);
\draw[color = black, opacity = 0.1] (2.000000, -3.464102) -- (4.000000, 0.000000);
\draw[color = black, opacity = 0.1] (2.500000, -4.330127) -- (5.000000, 0.000000);
\draw[color = black, opacity = 0.1] (3.000000, -5.196152) -- (6.000000, 0.000000);
\draw[color = black, opacity = 0.1] (3.535898, -6.000000) -- (7.000000, 0.000000);
\draw[color = black, opacity = 0.1] (4.535898, -6.000000) -- (8.000000, 0.000000);
\draw[color = black, opacity = 0.1] (5.535898, -6.000000) -- (9.000000, 0.000000);
\draw[color = black, opacity = 0.1] (6.535898, -6.000000) -- (10.000000, -0.000000);
\draw[color = black, opacity = 0.1] (7.535898, -6.000000) -- (10.000000, -1.732051);
\draw[color = black, opacity = 0.1] (8.535898, -6.000000) -- (10.000000, -3.464102);
\draw[color = black, opacity = 0.1] (9.535898, -6.000000) -- (10.000000, -5.196152);
\end{tikzpicture}

%% file: images/strip.tex
\begin{tikzpicture}[scale=1.0]
\fill[color = blue, opacity = 0.3]  (1.000000, 0.000000) --  (1.500000, 0.866025) --  (0.500000, 0.866025) -- cycle;

\draw[color = brown, line width = 2]  (0.500000, 0.866025) --  (5.500000, 0.866025) -- cycle;
\fill[color = red, opacity = 0.3]  (1.000000, 0.000000) --  (5.000000, 0.000000) --  (5.500000, 0.866025) --  (1.500000, 0.866025) -- cycle;
\draw[color = black, opacity = 0.1] (0.000000, 0.000000) -- (7.000000, 0.000000);
\draw[color = black, opacity = 0.1] (0.500000, 0.866025) -- (7.000000, 0.866025);
\draw[color = black, opacity = 0.1] (0.000000, 0.000000) -- (0.577350, 1.000000);
\draw[color = black, opacity = 0.1] (0.000000, 0.000000) -- (0.000000, 0.000000);
\draw[color = black, opacity = 0.1] (1.000000, 0.000000) -- (0.500000, 0.866025);
\draw[color = black, opacity = 0.1] (2.000000, 0.000000) -- (1.422650, 1.000000);
\draw[color = black, opacity = 0.1] (3.000000, 0.000000) -- (2.422650, 1.000000);
\draw[color = black, opacity = 0.1] (4.000000, 0.000000) -- (3.422650, 1.000000);
\draw[color = black, opacity = 0.1] (5.000000, 0.000000) -- (4.422650, 1.000000);
\draw[color = black, opacity = 0.1] (6.000000, 0.000000) -- (5.422650, 1.000000);
\draw[color = black, opacity = 0.1] (7.000000, 0.000000) -- (6.422650, 1.000000);
\draw[color = black, opacity = 0.1] (1.000000, 0.000000) -- (1.577350, 1.000000);
\draw[color = black, opacity = 0.1] (2.000000, 0.000000) -- (2.577350, 1.000000);
\draw[color = black, opacity = 0.1] (3.000000, 0.000000) -- (3.577350, 1.000000);
\draw[color = black, opacity = 0.1] (4.000000, 0.000000) -- (4.577350, 1.000000);
\draw[color = black, opacity = 0.1] (5.000000, 0.000000) -- (5.577350, 1.000000);
\draw[color = black, opacity = 0.1] (6.000000, 0.000000) -- (6.577350, 1.000000);
\draw[color = black, opacity = 0.1] (7.000000, 0.000000) -- (7.000000, 0.000000);
\draw[color = black, opacity = 0.1] (0.000000, 0.000000) -- (0.000000, 0.000000);
\draw[color = black, opacity = 0.1] (-0.500000, 0.866025) -- (0.500000, 0.866025);
\draw[color = black, opacity = 0.1] (0.000000, 0.000000) -- (0.577350, 1.000000);
\draw[color = black, opacity = 0.1] (-0.500000, 0.866025) -- (-0.422650, 1.000000);
\draw[color = black, opacity = 0.1] (0.000000, 0.000000) -- (-0.577350, 1.000000);
\draw[color = black, opacity = 0.1] (0.500000, 0.866025) -- (0.422650, 1.000000);
\draw[color = black, opacity = 0.1] (-1.000000, 0.000000) -- (0.000000, 0.000000);
\draw[color = black, opacity = 0.1] (-1.000000, 0.866025) -- (-0.500000, 0.866025);
\draw[color = black, opacity = 0.1] (0.000000, 0.000000) -- (0.000000, 0.000000);
\draw[color = black, opacity = 0.1] (-1.000000, 0.000000) -- (-0.500000, 0.866025);
\draw[color = black, opacity = 0.1] (0.000000, 0.000000) -- (-0.577350, 1.000000);
\draw[color = black, opacity = 0.1] (-1.000000, 0.000000) -- (-1.000000, 0.000000);
\draw[color = black, opacity = 0.1] (-1.000000, 0.000000) -- (0.000000, 0.000000);
\draw[color = black, opacity = 0.1] (-0.577350, -1.000000) -- (0.000000, 0.000000);
\draw[color = black, opacity = 0.1] (-1.000000, 0.000000) -- (-1.000000, 0.000000);
\draw[color = black, opacity = 0.1] (0.000000, 0.000000) -- (0.000000, 0.000000);
\draw[color = black, opacity = 0.1] (-1.000000, -0.866025) -- (-0.500000, -0.866025);
\draw[color = black, opacity = 0.1] (-0.500000, -0.866025) -- (-1.000000, 0.000000);
\draw[color = black, opacity = 0.1] (0.000000, 0.000000) -- (0.000000, 0.000000);
\draw[color = black, opacity = 0.1] (-0.577350, -1.000000) -- (0.000000, 0.000000);
\draw[color = black, opacity = 0.1] (0.577350, -1.000000) -- (0.000000, 0.000000);
\draw[color = black, opacity = 0.1] (-0.500000, -0.866025) -- (0.500000, -0.866025);
\draw[color = black, opacity = 0.1] (0.422650, -1.000000) -- (0.500000, -0.866025);
\draw[color = black, opacity = 0.1] (-0.422650, -1.000000) -- (-0.500000, -0.866025);
\draw[color = black, opacity = 0.1] (0.000000, 0.000000) -- (7.000000, 0.000000);
\draw[color = black, opacity = 0.1] (0.000000, 0.000000) -- (0.000000, 0.000000);
\draw[color = black, opacity = 0.1] (0.577350, -1.000000) -- (0.000000, 0.000000);
\draw[color = black, opacity = 0.1] (1.577350, -1.000000) -- (1.000000, 0.000000);
\draw[color = black, opacity = 0.1] (2.577350, -1.000000) -- (2.000000, 0.000000);
\draw[color = black, opacity = 0.1] (3.577350, -1.000000) -- (3.000000, 0.000000);
\draw[color = black, opacity = 0.1] (4.577350, -1.000000) -- (4.000000, 0.000000);
\draw[color = black, opacity = 0.1] (5.577350, -1.000000) -- (5.000000, 0.000000);
\draw[color = black, opacity = 0.1] (6.577350, -1.000000) -- (6.000000, 0.000000);
\draw[color = black, opacity = 0.1] (7.000000, 0.000000) -- (7.000000, 0.000000);
\draw[color = black, opacity = 0.1] (0.500000, -0.866025) -- (7.000000, -0.866025);
\draw[color = black, opacity = 0.1] (0.500000, -0.866025) -- (1.000000, 0.000000);
\draw[color = black, opacity = 0.1] (1.422650, -1.000000) -- (2.000000, 0.000000);
\draw[color = black, opacity = 0.1] (2.422650, -1.000000) -- (3.000000, 0.000000);
\draw[color = black, opacity = 0.1] (3.422650, -1.000000) -- (4.000000, 0.000000);
\draw[color = black, opacity = 0.1] (4.422650, -1.000000) -- (5.000000, 0.000000);
\draw[color = black, opacity = 0.1] (5.422650, -1.000000) -- (6.000000, 0.000000);
\draw[color = black, opacity = 0.1] (6.422650, -1.000000) -- (7.000000, 0.000000);

\fill[color = green] (1.000000, 0.000000) circle (2pt);
\node[yshift=-3mm] at (1, 0) {$w$};
\node[yshift=6mm] at (1, 0) {$\mathfrak{c}$};

\node[yshift = 3mm] at (2.7500000, 0.866025) {$\ell$};

\node[xshift = -3mm] at (0.500000, 0.866025) {$v_1$};
\fill[color = brown] (0.500000, 0.866025) circle (2pt);
\fill[color = brown] (5.500000, 0.866025) circle (2pt);
\node[xshift = 3mm] at (5.500000, 0.866025) {$v_2$};
\end{tikzpicture}

%% file: images/sector_configuration.tex
\begin{tikzpicture}[scale=0.75]
\fill[color = green, opacity = 0.3]  (0.000000, 0.000000) --  (0.500000, -0.866025) --  (-0.500000, -0.866025) -- cycle;

\node[yshift=-4mm] at (0, 0) {$\mathfrak{t}$};

\node[yshift = 3mm] at (0, 0) {$p$};

\node[] at (0.5, -0.3) {$\mathfrak{s}_2$};

\node[] at (-0.5, -0.3) {$\mathfrak{s}_1$};

\draw[color = red, line width = 0.5mm] (0, 0) -- (5, 0);

\draw[color = red, line width = 0.5mm] (0, 0) -- (2.886751, -5.000000);

\draw[color = blue, line width = 0.5mm] (0, 0) -- (-5, 0);

\draw[color = blue, line width = 0.2mm] (-1, 0) -- (-0.5, -0.866);

\draw[color = red, line width = 0.2mm] (1, 0) -- (0.5, -0.866);

\draw[color = green, line width = 0.2mm] (0.5, -0.866) -- (-0.5, -0.866);

\draw[color = blue, line width = 0.5mm] (0, 0) -- (-2.886751, -5.000000);

\node[yshift = 5mm] at (-1.44, -2.5) {$b_1$};

\node[yshift = 5mm] at (1.44, -2.5) {$b_2$};

\node[yshift = -2mm] at (-2.5, 0) {$a_1$};

\node[yshift = -2mm] at (2.5, 0) {$a_2$};

\node[] at (-4, -3) {$\mathcal{S}_1$};

\node[] at (4, -3) {$\mathcal{S}_2$};

\node[] at (-2, 2) {$\mathfrak{H}$};

\draw[color = red, opacity = 0.1] (0.000000, 0.000000) -- (5.000000, 0.000000);
\draw[color = red, opacity = 0.1] (0.000000, 0.000000) -- (0.000000, 0.000000);
\draw[color = red, opacity = 0.1] (2.886751, -5.000000) -- (0.000000, 0.000000);
\draw[color = red, opacity = 0.1] (3.886751, -5.000000) -- (1.000000, 0.000000);
\draw[color = red, opacity = 0.1] (4.886751, -5.000000) -- (2.000000, 0.000000);
\draw[color = red, opacity = 0.1] (5.000000, -3.464102) -- (3.000000, 0.000000);
\draw[color = red, opacity = 0.1] (5.000000, -1.732051) -- (4.000000, 0.000000);
\draw[color = red, opacity = 0.1] (5.000000, 0.000000) -- (5.000000, 0.000000);
\draw[color = red, opacity = 0.1] (0.500000, -0.866025) -- (5.000000, -0.866025);
\draw[color = red, opacity = 0.1] (1.000000, -1.732051) -- (5.000000, -1.732051);
\draw[color = red, opacity = 0.1] (1.500000, -2.598076) -- (5.000000, -2.598076);
\draw[color = red, opacity = 0.1] (2.000000, -3.464102) -- (5.000000, -3.464102);
\draw[color = red, opacity = 0.1] (2.500000, -4.330127) -- (5.000000, -4.330127);
\draw[color = red, opacity = 0.1] (0.500000, -0.866025) -- (1.000000, 0.000000);
\draw[color = red, opacity = 0.1] (1.000000, -1.732051) -- (2.000000, 0.000000);
\draw[color = red, opacity = 0.1] (1.500000, -2.598076) -- (3.000000, 0.000000);
\draw[color = red, opacity = 0.1] (2.000000, -3.464102) -- (4.000000, 0.000000);
\draw[color = red, opacity = 0.1] (2.500000, -4.330127) -- (5.000000, -0.000000);
\draw[color = red, opacity = 0.1] (3.113249, -5.000000) -- (5.000000, -1.732051);
\draw[color = red, opacity = 0.1] (4.113249, -5.000000) -- (5.000000, -3.464102);
\fill[color = red, opacity = 0.3] (0.000000, 0.000000) -- (2.886751, -5.000000) -- (5.000000, -5.000000) -- (5.000000, 0.000000) -- cycle;
\draw[color = yellow, opacity = 0.1] (0.000000, 0.000000) -- (5.000000, 0.000000);
\draw[color = yellow, opacity = 0.1] (0.500000, 0.866025) -- (5.000000, 0.866025);
\draw[color = yellow, opacity = 0.1] (1.000000, 1.732051) -- (5.000000, 1.732051);
\draw[color = yellow, opacity = 0.1] (1.500000, 2.598076) -- (5.000000, 2.598076);
\draw[color = yellow, opacity = 0.1] (0.000000, 0.000000) -- (1.732051, 3.000000);
\draw[color = yellow, opacity = 0.1] (0.000000, 0.000000) -- (0.000000, 0.000000);
\draw[color = yellow, opacity = 0.1] (1.000000, 0.000000) -- (0.500000, 0.866025);
\draw[color = yellow, opacity = 0.1] (2.000000, 0.000000) -- (1.000000, 1.732051);
\draw[color = yellow, opacity = 0.1] (3.000000, 0.000000) -- (1.500000, 2.598076);
\draw[color = yellow, opacity = 0.1] (4.000000, 0.000000) -- (2.267949, 3.000000);
\draw[color = yellow, opacity = 0.1] (5.000000, 0.000000) -- (3.267949, 3.000000);
\draw[color = yellow, opacity = 0.1] (5.000000, 1.732051) -- (4.267949, 3.000000);
\draw[color = yellow, opacity = 0.1] (1.000000, 0.000000) -- (2.732051, 3.000000);
\draw[color = yellow, opacity = 0.1] (2.000000, 0.000000) -- (3.732051, 3.000000);
\draw[color = yellow, opacity = 0.1] (3.000000, 0.000000) -- (4.732051, 3.000000);
\draw[color = yellow, opacity = 0.1] (4.000000, 0.000000) -- (5.000000, 1.732051);
\fill[color = yellow, opacity = 0.3] (0.000000, 0.000000) -- (5.000000, 0.000000) -- (5.000000, 3.000000) -- (1.732051, 3.000000) -- cycle;
\draw[color = black, opacity = 0.1] (0.000000, 0.000000) -- (5.000000, 0.000000);
\draw[color = black, opacity = 0.1] (0.500000, 0.866025) -- (5.000000, 0.866025);
\draw[color = black, opacity = 0.1] (1.000000, 1.732051) -- (5.000000, 1.732051);
\draw[color = black, opacity = 0.1] (1.500000, 2.598076) -- (5.000000, 2.598076);
\draw[color = black, opacity = 0.1] (0.000000, 0.000000) -- (1.732051, 3.000000);
\draw[color = black, opacity = 0.1] (0.000000, 0.000000) -- (0.000000, 0.000000);
\draw[color = black, opacity = 0.1] (1.000000, 0.000000) -- (0.500000, 0.866025);
\draw[color = black, opacity = 0.1] (2.000000, 0.000000) -- (1.000000, 1.732051);
\draw[color = black, opacity = 0.1] (3.000000, 0.000000) -- (1.500000, 2.598076);
\draw[color = black, opacity = 0.1] (4.000000, 0.000000) -- (2.267949, 3.000000);
\draw[color = black, opacity = 0.1] (5.000000, 0.000000) -- (3.267949, 3.000000);
\draw[color = black, opacity = 0.1] (5.000000, 1.732051) -- (4.267949, 3.000000);
\draw[color = black, opacity = 0.1] (1.000000, 0.000000) -- (2.732051, 3.000000);
\draw[color = black, opacity = 0.1] (2.000000, 0.000000) -- (3.732051, 3.000000);
\draw[color = black, opacity = 0.1] (3.000000, 0.000000) -- (4.732051, 3.000000);
\draw[color = black, opacity = 0.1] (4.000000, 0.000000) -- (5.000000, 1.732051);
\draw[color = yellow, opacity = 0.1] (0.000000, 0.000000) -- (0.000000, 0.000000);
\draw[color = yellow, opacity = 0.1] (-0.500000, 0.866025) -- (0.500000, 0.866025);
\draw[color = yellow, opacity = 0.1] (-1.000000, 1.732051) -- (1.000000, 1.732051);
\draw[color = yellow, opacity = 0.1] (-1.500000, 2.598076) -- (1.500000, 2.598076);
\draw[color = yellow, opacity = 0.1] (0.000000, 0.000000) -- (1.732051, 3.000000);
\draw[color = yellow, opacity = 0.1] (-0.500000, 0.866025) -- (0.732051, 3.000000);
\draw[color = yellow, opacity = 0.1] (-1.000000, 1.732051) -- (-0.267949, 3.000000);
\draw[color = yellow, opacity = 0.1] (-1.500000, 2.598076) -- (-1.267949, 3.000000);
\draw[color = yellow, opacity = 0.1] (0.000000, 0.000000) -- (-1.732051, 3.000000);
\draw[color = yellow, opacity = 0.1] (0.500000, 0.866025) -- (-0.732051, 3.000000);
\draw[color = yellow, opacity = 0.1] (1.000000, 1.732051) -- (0.267949, 3.000000);
\draw[color = yellow, opacity = 0.1] (1.500000, 2.598076) -- (1.267949, 3.000000);
\fill[color = yellow, opacity = 0.3] (0.000000, 0.000000) -- (1.732051, 3.000000) -- (-1.732051, 3.000000) -- cycle;
\draw[color = black, opacity = 0.1] (0.000000, 0.000000) -- (0.000000, 0.000000);
\draw[color = black, opacity = 0.1] (-0.500000, 0.866025) -- (0.500000, 0.866025);
\draw[color = black, opacity = 0.1] (-1.000000, 1.732051) -- (1.000000, 1.732051);
\draw[color = black, opacity = 0.1] (-1.500000, 2.598076) -- (1.500000, 2.598076);
\draw[color = black, opacity = 0.1] (0.000000, 0.000000) -- (1.732051, 3.000000);
\draw[color = black, opacity = 0.1] (-0.500000, 0.866025) -- (0.732051, 3.000000);
\draw[color = black, opacity = 0.1] (-1.000000, 1.732051) -- (-0.267949, 3.000000);
\draw[color = black, opacity = 0.1] (-1.500000, 2.598076) -- (-1.267949, 3.000000);
\draw[color = black, opacity = 0.1] (0.000000, 0.000000) -- (-1.732051, 3.000000);
\draw[color = black, opacity = 0.1] (0.500000, 0.866025) -- (-0.732051, 3.000000);
\draw[color = black, opacity = 0.1] (1.000000, 1.732051) -- (0.267949, 3.000000);
\draw[color = black, opacity = 0.1] (1.500000, 2.598076) -- (1.267949, 3.000000);
\draw[color = yellow, opacity = 0.1] (-5.000000, 0.000000) -- (0.000000, 0.000000);
\draw[color = yellow, opacity = 0.1] (-5.000000, 0.866025) -- (-0.500000, 0.866025);
\draw[color = yellow, opacity = 0.1] (-5.000000, 1.732051) -- (-1.000000, 1.732051);
\draw[color = yellow, opacity = 0.1] (-5.000000, 2.598076) -- (-1.500000, 2.598076);
\draw[color = yellow, opacity = 0.1] (0.000000, 0.000000) -- (0.000000, 0.000000);
\draw[color = yellow, opacity = 0.1] (-1.000000, 0.000000) -- (-0.500000, 0.866025);
\draw[color = yellow, opacity = 0.1] (-2.000000, 0.000000) -- (-1.000000, 1.732051);
\draw[color = yellow, opacity = 0.1] (-3.000000, 0.000000) -- (-1.500000, 2.598076);
\draw[color = yellow, opacity = 0.1] (-4.000000, 0.000000) -- (-2.267949, 3.000000);
\draw[color = yellow, opacity = 0.1] (-5.000000, 0.000000) -- (-3.267949, 3.000000);
\draw[color = yellow, opacity = 0.1] (-5.000000, 1.732051) -- (-4.267949, 3.000000);
\draw[color = yellow, opacity = 0.1] (0.000000, 0.000000) -- (-1.732051, 3.000000);
\draw[color = yellow, opacity = 0.1] (-1.000000, 0.000000) -- (-2.732051, 3.000000);
\draw[color = yellow, opacity = 0.1] (-2.000000, 0.000000) -- (-3.732051, 3.000000);
\draw[color = yellow, opacity = 0.1] (-3.000000, 0.000000) -- (-4.732051, 3.000000);
\draw[color = yellow, opacity = 0.1] (-4.000000, 0.000000) -- (-5.000000, 1.732051);
\draw[color = yellow, opacity = 0.1] (-5.000000, 0.000000) -- (-5.000000, 0.000000);
\fill[color = yellow, opacity = 0.3] (0.000000, 0.000000) -- (-1.732051, 3.000000) -- (-5.000000, 3.000000) -- (-5.000000, 0.000000) -- cycle;
\draw[color = black, opacity = 0.1] (-5.000000, 0.000000) -- (0.000000, 0.000000);
\draw[color = black, opacity = 0.1] (-5.000000, 0.866025) -- (-0.500000, 0.866025);
\draw[color = black, opacity = 0.1] (-5.000000, 1.732051) -- (-1.000000, 1.732051);
\draw[color = black, opacity = 0.1] (-5.000000, 2.598076) -- (-1.500000, 2.598076);
\draw[color = black, opacity = 0.1] (0.000000, 0.000000) -- (0.000000, 0.000000);
\draw[color = black, opacity = 0.1] (-1.000000, 0.000000) -- (-0.500000, 0.866025);
\draw[color = black, opacity = 0.1] (-2.000000, 0.000000) -- (-1.000000, 1.732051);
\draw[color = black, opacity = 0.1] (-3.000000, 0.000000) -- (-1.500000, 2.598076);
\draw[color = black, opacity = 0.1] (-4.000000, 0.000000) -- (-2.267949, 3.000000);
\draw[color = black, opacity = 0.1] (-5.000000, 0.000000) -- (-3.267949, 3.000000);
\draw[color = black, opacity = 0.1] (-5.000000, 1.732051) -- (-4.267949, 3.000000);
\draw[color = black, opacity = 0.1] (0.000000, 0.000000) -- (-1.732051, 3.000000);
\draw[color = black, opacity = 0.1] (-1.000000, 0.000000) -- (-2.732051, 3.000000);
\draw[color = black, opacity = 0.1] (-2.000000, 0.000000) -- (-3.732051, 3.000000);
\draw[color = black, opacity = 0.1] (-3.000000, 0.000000) -- (-4.732051, 3.000000);
\draw[color = black, opacity = 0.1] (-4.000000, 0.000000) -- (-5.000000, 1.732051);
\draw[color = black, opacity = 0.1] (-5.000000, 0.000000) -- (-5.000000, 0.000000);
\draw[color = blue, opacity = 0.1] (-5.000000, 0.000000) -- (0.000000, 0.000000);
\draw[color = blue, opacity = 0.1] (-2.886751, -5.000000) -- (0.000000, 0.000000);
\draw[color = blue, opacity = 0.1] (-3.886751, -5.000000) -- (-1.000000, 0.000000);
\draw[color = blue, opacity = 0.1] (-4.886751, -5.000000) -- (-2.000000, 0.000000);
\draw[color = blue, opacity = 0.1] (-5.000000, -3.464102) -- (-3.000000, 0.000000);
\draw[color = blue, opacity = 0.1] (-5.000000, -1.732051) -- (-4.000000, 0.000000);
\draw[color = blue, opacity = 0.1] (0.000000, 0.000000) -- (0.000000, 0.000000);
\draw[color = blue, opacity = 0.1] (-5.000000, -0.866025) -- (-0.500000, -0.866025);
\draw[color = blue, opacity = 0.1] (-5.000000, -1.732051) -- (-1.000000, -1.732051);
\draw[color = blue, opacity = 0.1] (-5.000000, -2.598076) -- (-1.500000, -2.598076);
\draw[color = blue, opacity = 0.1] (-5.000000, -3.464102) -- (-2.000000, -3.464102);
\draw[color = blue, opacity = 0.1] (-5.000000, -4.330127) -- (-2.500000, -4.330127);
\draw[color = blue, opacity = 0.1] (-0.500000, -0.866025) -- (-1.000000, 0.000000);
\draw[color = blue, opacity = 0.1] (-1.000000, -1.732051) -- (-2.000000, 0.000000);
\draw[color = blue, opacity = 0.1] (-1.500000, -2.598076) -- (-3.000000, 0.000000);
\draw[color = blue, opacity = 0.1] (-2.000000, -3.464102) -- (-4.000000, 0.000000);
\draw[color = blue, opacity = 0.1] (-2.500000, -4.330127) -- (-5.000000, 0.000000);
\draw[color = blue, opacity = 0.1] (-3.113249, -5.000000) -- (-5.000000, -1.732051);
\draw[color = blue, opacity = 0.1] (-4.113249, -5.000000) -- (-5.000000, -3.464102);
\fill[color = blue, opacity = 0.3] (0.000000, 0.000000) -- (-5.000000, 0.000000) -- (-5.000000, -5.000000) -- (-2.886751, -5.000000) -- cycle;

\fill[color = black] (0, 0) circle (2pt);
\end{tikzpicture}

%% file: images/apartment_perspective.tex
\begin{tikzpicture}[scale=0.5]
\draw[color = black, opacity = 0.2] (-6.000000, -3.464102) -- (6.000000, -3.464102);
\draw[color = black, opacity = 0.2] (-5.500000, -2.598076) -- (6.500000, -2.598076);
\draw[color = black, opacity = 0.2] (-5.000000, -1.732051) -- (7.000000, -1.732051);
\draw[color = black, opacity = 0.2] (-4.500000, -0.866025) -- (7.500000, -0.866025);
\draw[color = black, opacity = 0.2] (-4.000000, 0.000000) -- (8.000000, 0.000000);
\draw[color = black, opacity = 0.2] (-3.500000, 0.866025) -- (8.500000, 0.866025);
\draw[color = black, opacity = 0.2] (-3.000000, 1.732051) -- (9.000000, 1.732051);
\draw[color = black, opacity = 0.2] (-2.500000, 2.598076) -- (9.500000, 2.598076);
\draw[color = black, opacity = 0.2] (-2.000000, 3.464102) -- (10.000000, 3.464102);
\draw[color = black, opacity = 0.2] (-1.500000, 4.330127) -- (10.500000, 4.330127);
\draw[color = black, opacity = 0.2] (-6.000000, -3.464102) -- (-1.500000, 4.330127);
\draw[color = black, opacity = 0.2] (-6.000000, -3.464102) -- (-6.000000, -3.464102);
\draw[color = black, opacity = 0.2] (-5.000000, -3.464102) -- (-5.500000, -2.598076);
\draw[color = black, opacity = 0.2] (-4.000000, -3.464102) -- (-5.000000, -1.732051);
\draw[color = black, opacity = 0.2] (-3.000000, -3.464102) -- (-4.500000, -0.866025);
\draw[color = black, opacity = 0.2] (-2.000000, -3.464102) -- (-4.000000, 0.000000);
\draw[color = black, opacity = 0.2] (-1.000000, -3.464102) -- (-3.500000, 0.866025);
\draw[color = black, opacity = 0.2] (0.000000, -3.464102) -- (-3.000000, 1.732051);
\draw[color = black, opacity = 0.2] (1.000000, -3.464102) -- (-2.500000, 2.598076);
\draw[color = black, opacity = 0.2] (2.000000, -3.464102) -- (-2.000000, 3.464102);
\draw[color = black, opacity = 0.2] (3.000000, -3.464102) -- (-1.500000, 4.330127);
\draw[color = black, opacity = 0.2] (4.000000, -3.464102) -- (-0.500000, 4.330127);
\draw[color = black, opacity = 0.2] (5.000000, -3.464102) -- (0.500000, 4.330127);
\draw[color = black, opacity = 0.2] (6.000000, -3.464102) -- (1.500000, 4.330127);
\draw[color = black, opacity = 0.2] (6.500000, -2.598076) -- (2.500000, 4.330127);
\draw[color = black, opacity = 0.2] (7.000000, -1.732051) -- (3.500000, 4.330127);
\draw[color = black, opacity = 0.2] (7.500000, -0.866025) -- (4.500000, 4.330127);
\draw[color = black, opacity = 0.2] (8.000000, 0.000000) -- (5.500000, 4.330127);
\draw[color = black, opacity = 0.2] (8.500000, 0.866025) -- (6.500000, 4.330127);
\draw[color = black, opacity = 0.2] (9.000000, 1.732051) -- (7.500000, 4.330127);
\draw[color = black, opacity = 0.2] (9.500000, 2.598076) -- (8.500000, 4.330127);
\draw[color = black, opacity = 0.2] (10.000000, 3.464102) -- (9.500000, 4.330127);
\draw[color = black, opacity = 0.2] (-6.000000, -3.464102) -- (6.000000, -3.464102);
\draw[color = black, opacity = 0.2] (-6.000000, -3.464102) -- (-1.500000, 4.330127);
\draw[color = black, opacity = 0.2] (-5.000000, -3.464102) -- (-0.500000, 4.330127);
\draw[color = black, opacity = 0.2] (-4.000000, -3.464102) -- (0.500000, 4.330127);
\draw[color = black, opacity = 0.2] (-3.000000, -3.464102) -- (1.500000, 4.330127);
\draw[color = black, opacity = 0.2] (-2.000000, -3.464102) -- (2.500000, 4.330127);
\draw[color = black, opacity = 0.2] (-1.000000, -3.464102) -- (3.500000, 4.330127);
\draw[color = black, opacity = 0.2] (0.000000, -3.464102) -- (4.500000, 4.330127);
\draw[color = black, opacity = 0.2] (1.000000, -3.464102) -- (5.500000, 4.330127);
\draw[color = black, opacity = 0.2] (2.000000, -3.464102) -- (6.500000, 4.330127);
\draw[color = black, opacity = 0.2] (3.000000, -3.464102) -- (7.500000, 4.330127);
\draw[color = black, opacity = 0.2] (4.000000, -3.464102) -- (8.500000, 4.330127);
\draw[color = black, opacity = 0.2] (5.000000, -3.464102) -- (9.500000, 4.330127);
\draw[color = black, opacity = 0.2] (-6.000000, -3.464102) -- (-6.000000, -3.464102);
\fill[color = brown, opacity = 0.3] (-6.000000, -3.464102) -- (6.000000, -3.464102) -- (10.500000, 4.330127) -- (-1.500000, 4.330127) -- cycle;
\draw[color = black, opacity = 0.2] (-3.500000, -2.598076) -- (3.500000, -2.598076);
\draw[color = black, opacity = 0.2] (-3.000000, -1.732051) -- (4.000000, -1.732051);
\draw[color = black, opacity = 0.2] (-2.500000, -0.866025) -- (4.500000, -0.866025);
\draw[color = black, opacity = 0.2] (-2.000000, 0.000000) -- (5.000000, 0.000000);
\draw[color = black, opacity = 0.2] (-1.500000, 0.866025) -- (5.500000, 0.866025);
\draw[color = black, opacity = 0.2] (-1.000000, 1.732051) -- (6.000000, 1.732051);
\draw[color = black, opacity = 0.2] (-0.500000, 2.598076) -- (6.500000, 2.598076);
\draw[color = black, opacity = 0.2] (-3.500000, -2.598076) -- (-0.500000, 2.598076);
\draw[color = black, opacity = 0.2] (-3.500000, -2.598076) -- (-3.500000, -2.598076);
\draw[color = black, opacity = 0.2] (-2.500000, -2.598076) -- (-3.000000, -1.732051);
\draw[color = black, opacity = 0.2] (-1.500000, -2.598076) -- (-2.500000, -0.866025);
\draw[color = black, opacity = 0.2] (-0.500000, -2.598076) -- (-2.000000, 0.000000);
\draw[color = black, opacity = 0.2] (0.500000, -2.598076) -- (-1.500000, 0.866025);
\draw[color = black, opacity = 0.2] (1.500000, -2.598076) -- (-1.000000, 1.732051);
\draw[color = black, opacity = 0.2] (2.500000, -2.598076) -- (-0.500000, 2.598076);
\draw[color = black, opacity = 0.2] (3.500000, -2.598076) -- (0.500000, 2.598076);
\draw[color = black, opacity = 0.2] (4.000000, -1.732051) -- (1.500000, 2.598076);
\draw[color = black, opacity = 0.2] (4.500000, -0.866025) -- (2.500000, 2.598076);
\draw[color = black, opacity = 0.2] (5.000000, 0.000000) -- (3.500000, 2.598076);
\draw[color = black, opacity = 0.2] (5.500000, 0.866025) -- (4.500000, 2.598076);
\draw[color = black, opacity = 0.2] (6.000000, 1.732051) -- (5.500000, 2.598076);
\draw[color = black, opacity = 0.2] (-3.500000, -2.598076) -- (3.500000, -2.598076);
\draw[color = black, opacity = 0.2] (-3.500000, -2.598076) -- (-0.500000, 2.598076);
\draw[color = black, opacity = 0.2] (-2.500000, -2.598076) -- (0.500000, 2.598076);
\draw[color = black, opacity = 0.2] (-1.500000, -2.598076) -- (1.500000, 2.598076);
\draw[color = black, opacity = 0.2] (-0.500000, -2.598076) -- (2.500000, 2.598076);
\draw[color = black, opacity = 0.2] (0.500000, -2.598076) -- (3.500000, 2.598076);
\draw[color = black, opacity = 0.2] (1.500000, -2.598076) -- (4.500000, 2.598076);
\draw[color = black, opacity = 0.2] (2.500000, -2.598076) -- (5.500000, 2.598076);
\draw[color = black, opacity = 0.2] (3.500000, -2.598076) -- (6.500000, 2.598076);
\draw[color = black, opacity = 0.2] (-3.500000, -2.598076) -- (-3.500000, -2.598076);
\fill[color = brown, opacity = 0.3] (-3.500000, -2.598076) -- (3.500000, -2.598076) -- (6.500000, 2.598076) -- (-0.500000, 2.598076) -- cycle;

\fill[color = blue, opacity = 0.6] (-6.000000, -3.464102) -- (-2.000000, -3.464102) -- (0.000000, 0.000000) -- (-4.000000, 0.000000) -- cycle;

\draw[color = black, opacity = 0.2] (-6.000000, -3.464102) -- (-2.000000, -3.464102);
\draw[color = black, opacity = 0.2] (-5.500000, -2.598076) -- (-1.500000, -2.598076);
\draw[color = black, opacity = 0.2] (-5.000000, -1.732051) -- (-1.000000, -1.732051);
\draw[color = black, opacity = 0.2] (-4.500000, -0.866025) -- (-0.500000, -0.866025);
\draw[color = black, opacity = 0.2] (-4.000000, 0.000000) -- (0.000000, 0.000000);
\draw[color = black, opacity = 0.2] (-6.000000, -3.464102) -- (-4.000000, 0.000000);
\draw[color = black, opacity = 0.2] (-6.000000, -3.464102) -- (-6.000000, -3.464102);
\draw[color = black, opacity = 0.2] (-5.000000, -3.464102) -- (-5.500000, -2.598076);
\draw[color = black, opacity = 0.2] (-4.000000, -3.464102) -- (-5.000000, -1.732051);
\draw[color = black, opacity = 0.2] (-3.000000, -3.464102) -- (-4.500000, -0.866025);
\draw[color = black, opacity = 0.2] (-2.000000, -3.464102) -- (-4.000000, 0.000000);
\draw[color = black, opacity = 0.2] (-1.500000, -2.598076) -- (-3.000000, 0.000000);
\draw[color = black, opacity = 0.2] (-1.000000, -1.732051) -- (-2.000000, 0.000000);
\draw[color = black, opacity = 0.2] (-0.500000, -0.866025) -- (-1.000000, 0.000000);
\draw[color = black, opacity = 0.2] (-6.000000, -3.464102) -- (-2.000000, -3.464102);
\draw[color = black, opacity = 0.2] (-6.000000, -3.464102) -- (-4.000000, 0.000000);
\draw[color = black, opacity = 0.2] (-5.000000, -3.464102) -- (-3.000000, 0.000000);
\draw[color = black, opacity = 0.2] (-4.000000, -3.464102) -- (-2.000000, 0.000000);
\draw[color = black, opacity = 0.2] (-3.000000, -3.464102) -- (-1.000000, 0.000000);
\draw[color = black, opacity = 0.2] (-6.000000, -3.464102) -- (-6.000000, -3.464102);

\fill[color = red, opacity = 0.6] (3.000000, 0.000000) -- (8.000000, 0.000000) -- (10.500000, 4.330127) -- (5.500000, 4.330127) -- cycle;

\draw[color = black, opacity = 0.2] (3.000000, 0.000000) -- (8.000000, 0.000000);
\draw[color = black, opacity = 0.2] (3.500000, 0.866025) -- (8.500000, 0.866025);
\draw[color = black, opacity = 0.2] (4.000000, 1.732051) -- (9.000000, 1.732051);
\draw[color = black, opacity = 0.2] (4.500000, 2.598076) -- (9.500000, 2.598076);
\draw[color = black, opacity = 0.2] (5.000000, 3.464102) -- (10.000000, 3.464102);
\draw[color = black, opacity = 0.2] (5.500000, 4.330127) -- (10.500000, 4.330127);
\draw[color = black, opacity = 0.2] (3.000000, 0.000000) -- (5.500000, 4.330127);
\draw[color = black, opacity = 0.2] (3.000000, 0.000000) -- (3.000000, 0.000000);
\draw[color = black, opacity = 0.2] (4.000000, 0.000000) -- (3.500000, 0.866025);
\draw[color = black, opacity = 0.2] (5.000000, 0.000000) -- (4.000000, 1.732051);
\draw[color = black, opacity = 0.2] (6.000000, 0.000000) -- (4.500000, 2.598076);
\draw[color = black, opacity = 0.2] (7.000000, 0.000000) -- (5.000000, 3.464102);
\draw[color = black, opacity = 0.2] (8.000000, 0.000000) -- (5.500000, 4.330127);
\draw[color = black, opacity = 0.2] (8.500000, 0.866025) -- (6.500000, 4.330127);
\draw[color = black, opacity = 0.2] (9.000000, 1.732051) -- (7.500000, 4.330127);
\draw[color = black, opacity = 0.2] (9.500000, 2.598076) -- (8.500000, 4.330127);
\draw[color = black, opacity = 0.2] (10.000000, 3.464102) -- (9.500000, 4.330127);
\draw[color = black, opacity = 0.2] (3.000000, 0.000000) -- (8.000000, 0.000000);
\draw[color = black, opacity = 0.2] (3.000000, 0.000000) -- (5.500000, 4.330127);
\draw[color = black, opacity = 0.2] (4.000000, 0.000000) -- (6.500000, 4.330127);
\draw[color = black, opacity = 0.2] (5.000000, 0.000000) -- (7.500000, 4.330127);
\draw[color = black, opacity = 0.2] (6.000000, 0.000000) -- (8.500000, 4.330127);
\draw[color = black, opacity = 0.2] (7.000000, 0.000000) -- (9.500000, 4.330127);
\draw[color = black, opacity = 0.2] (3.000000, 0.000000) -- (3.000000, 0.000000);

\draw[color = green, line width = 2] (0.000000, 0.000000) -- (3.000000, 0.000000);

\fill[color = blue] (0, 0) circle (3pt);

\node[yshift = 4mm] at (0, 0) {$p_1^{(k)}$};

\fill[color = red] (3, 0) circle (3 pt);

\node[yshift = 4mm] at (3, 0) {$p_2^{(k)}$};

\node[yshift = -3mm] at (1.5, 0) {$\ell$};
\draw[color = black, line width = 2] (-3.500000, -2.598076) -- (3.500000, -2.598076) -- (6.500000, 2.598076) -- (-0.500000, 2.598076) -- cycle;

\fill[color = black] (-3.500000, -2.598076) circle (5pt);

\node[yshift=-3mm] at (-3.500000, -2.598076) {$x$};

\fill[color = black] (6.500000, 2.598076) circle (5pt);

\node[yshift=3mm] at (6.500000, 2.598076) {$y$};
\end{tikzpicture}

%% file: images/perspective.tex
\begin{tikzpicture}[scale=0.5]

\draw[color = brown, opacity = 0.2] (-14.480185, -5.300115) -- (3.879954, -5.300115);
\draw[color = brown, opacity = 0.2] (-10.940502, -3.510001) -- (5.271500, -3.510001);
\draw[color = brown, opacity = 0.2] (-8.142347, -2.094898) -- (6.371531, -2.094898);
\draw[color = brown, opacity = 0.2] (-5.874798, -0.948136) -- (7.262966, -0.948136);
\draw[color = brown, opacity = 0.2] (-4.000000, 0.000000) -- (8.000000, 0.000000);
\draw[color = brown, opacity = 0.2] (-2.424046, 0.797003) -- (8.619550, 0.797003);
\draw[color = brown, opacity = 0.2] (-1.080757, 1.476341) -- (9.147634, 1.476341);
\draw[color = brown, opacity = 0.2] (0.077850, 2.062280) -- (9.603114, 2.062280);
\draw[color = brown, opacity = 0.2] (1.087411, 2.572843) -- (10.000000, 2.572843);
\draw[color = brown, opacity = 0.2] (1.974949, 3.021695) -- (10.348915, 3.021695);
\draw[color = brown, opacity = 0.2] (-14.480185, -5.300115) -- (1.974949, 3.021695);
\draw[color = brown, opacity = 0.2] (-14.480185, -5.300115) -- (-14.480185, -5.300115);
\draw[color = brown, opacity = 0.2] (-12.950173, -5.300115) -- (-10.940502, -3.510001);
\draw[color = brown, opacity = 0.2] (-11.420162, -5.300115) -- (-8.142347, -2.094898);
\draw[color = brown, opacity = 0.2] (-9.890150, -5.300115) -- (-5.874798, -0.948136);
\draw[color = brown, opacity = 0.2] (-8.360139, -5.300115) -- (-4.000000, 0.000000);
\draw[color = brown, opacity = 0.2] (-6.830127, -5.300115) -- (-2.424046, 0.797003);
\draw[color = brown, opacity = 0.2] (-5.300115, -5.300115) -- (-1.080757, 1.476341);
\draw[color = brown, opacity = 0.2] (-3.770104, -5.300115) -- (0.077850, 2.062280);
\draw[color = brown, opacity = 0.2] (-2.240092, -5.300115) -- (1.087411, 2.572843);
\draw[color = brown, opacity = 0.2] (-0.710081, -5.300115) -- (1.974949, 3.021695);
\draw[color = brown, opacity = 0.2] (0.819931, -5.300115) -- (2.672780, 3.021695);
\draw[color = brown, opacity = 0.2] (2.349942, -5.300115) -- (3.370610, 3.021695);
\draw[color = brown, opacity = 0.2] (3.879954, -5.300115) -- (4.068441, 3.021695);
\draw[color = brown, opacity = 0.2] (5.271500, -3.510001) -- (4.766271, 3.021695);
\draw[color = brown, opacity = 0.2] (6.371531, -2.094898) -- (5.464102, 3.021695);
\draw[color = brown, opacity = 0.2] (7.262966, -0.948136) -- (6.161932, 3.021695);
\draw[color = brown, opacity = 0.2] (8.000000, 0.000000) -- (6.859763, 3.021695);
\draw[color = brown, opacity = 0.2] (8.619550, 0.797003) -- (7.557593, 3.021695);
\draw[color = brown, opacity = 0.2] (9.147634, 1.476341) -- (8.255424, 3.021695);
\draw[color = brown, opacity = 0.2] (9.603114, 2.062280) -- (8.953254, 3.021695);
\draw[color = brown, opacity = 0.2] (10.000000, 2.572843) -- (9.651085, 3.021695);
\draw[color = brown, opacity = 0.2] (-14.480185, -5.300115) -- (3.879954, -5.300115);
\draw[color = brown, opacity = 0.2] (-14.480185, -5.300115) -- (1.974949, 3.021695);
\draw[color = brown, opacity = 0.2] (-12.950173, -5.300115) -- (2.672780, 3.021695);
\draw[color = brown, opacity = 0.2] (-11.420162, -5.300115) -- (3.370610, 3.021695);
\draw[color = brown, opacity = 0.2] (-9.890150, -5.300115) -- (4.068441, 3.021695);
\draw[color = brown, opacity = 0.2] (-8.360139, -5.300115) -- (4.766271, 3.021695);
\draw[color = brown, opacity = 0.2] (-6.830127, -5.300115) -- (5.464102, 3.021695);
\draw[color = brown, opacity = 0.2] (-5.300115, -5.300115) -- (6.161932, 3.021695);
\draw[color = brown, opacity = 0.2] (-3.770104, -5.300115) -- (6.859763, 3.021695);
\draw[color = brown, opacity = 0.2] (-2.240092, -5.300115) -- (7.557593, 3.021695);
\draw[color = brown, opacity = 0.2] (-0.710081, -5.300115) -- (8.255424, 3.021695);
\draw[color = brown, opacity = 0.2] (0.819931, -5.300115) -- (8.953254, 3.021695);
\draw[color = brown, opacity = 0.2] (2.349942, -5.300115) -- (9.651085, 3.021695);
\draw[color = brown, opacity = 0.2] (-14.480185, -5.300115) -- (-14.480185, -5.300115);
\fill[color = brown, opacity = 0.3] (-14.480185, -5.300115) -- (3.879954, -5.300115) -- (10.348915, 3.021695) -- (1.974949, 3.021695) -- cycle;
\fill[color = blue, opacity = 0.9] (-14.480185, -5.300115) -- (-8.360139, -5.300115) -- (0.000000, 0.000000) -- (-4.000000, 0.000000) -- cycle;
\fill[color = red, opacity = 0.9] (3.000000, 0.000000) -- (8.000000, 0.000000) -- (10.348915, 3.021695) -- (6.859763, 3.021695) -- cycle;
\draw[color = brown, opacity = 0.2] (-8.238502, -3.510001) -- (1.218499, -3.510001);
\draw[color = brown, opacity = 0.2] (-5.723367, -2.094898) -- (2.743061, -2.094898);
\draw[color = brown, opacity = 0.2] (-3.685171, -0.948136) -- (3.978525, -0.948136);
\draw[color = brown, opacity = 0.2] (-2.000000, 0.000000) -- (5.000000, 0.000000);
\draw[color = brown, opacity = 0.2] (-0.583447, 0.797003) -- (5.858651, 0.797003);
\draw[color = brown, opacity = 0.2] (0.623975, 1.476341) -- (6.590536, 1.476341);
\draw[color = brown, opacity = 0.2] (1.665394, 2.062280) -- (7.221798, 2.062280);
\draw[color = brown, opacity = 0.2] (-8.238502, -3.510001) -- (1.665394, 2.062280);
\draw[color = brown, opacity = 0.2] (-8.238502, -3.510001) -- (-8.238502, -3.510001);
\draw[color = brown, opacity = 0.2] (-6.887502, -3.510001) -- (-5.723367, -2.094898);
\draw[color = brown, opacity = 0.2] (-5.536501, -3.510001) -- (-3.685171, -0.948136);
\draw[color = brown, opacity = 0.2] (-4.185501, -3.510001) -- (-2.000000, 0.000000);
\draw[color = brown, opacity = 0.2] (-2.834501, -3.510001) -- (-0.583447, 0.797003);
\draw[color = brown, opacity = 0.2] (-1.483501, -3.510001) -- (0.623975, 1.476341);
\draw[color = brown, opacity = 0.2] (-0.132501, -3.510001) -- (1.665394, 2.062280);
\draw[color = brown, opacity = 0.2] (1.218499, -3.510001) -- (2.459166, 2.062280);
\draw[color = brown, opacity = 0.2] (2.743061, -2.094898) -- (3.252938, 2.062280);
\draw[color = brown, opacity = 0.2] (3.978525, -0.948136) -- (4.046710, 2.062280);
\draw[color = brown, opacity = 0.2] (5.000000, 0.000000) -- (4.840482, 2.062280);
\draw[color = brown, opacity = 0.2] (5.858651, 0.797003) -- (5.634254, 2.062280);
\draw[color = brown, opacity = 0.2] (6.590536, 1.476341) -- (6.428026, 2.062280);
\draw[color = brown, opacity = 0.2] (-8.238502, -3.510001) -- (1.218499, -3.510001);
\draw[color = brown, opacity = 0.2] (-8.238502, -3.510001) -- (1.665394, 2.062280);
\draw[color = brown, opacity = 0.2] (-6.887502, -3.510001) -- (2.459166, 2.062280);
\draw[color = brown, opacity = 0.2] (-5.536501, -3.510001) -- (3.252938, 2.062280);
\draw[color = brown, opacity = 0.2] (-4.185501, -3.510001) -- (4.046710, 2.062280);
\draw[color = brown, opacity = 0.2] (-2.834501, -3.510001) -- (4.840482, 2.062280);
\draw[color = brown, opacity = 0.2] (-1.483501, -3.510001) -- (5.634254, 2.062280);
\draw[color = brown, opacity = 0.2] (-0.132501, -3.510001) -- (6.428026, 2.062280);
\draw[color = brown, opacity = 0.2] (1.218499, -3.510001) -- (7.221798, 2.062280);
\draw[color = brown, opacity = 0.2] (-8.238502, -3.510001) -- (-8.238502, -3.510001);
\draw[color = black, line width = 2] (-8.238502, -3.510001) -- (1.218499, -3.510001) -- (7.221798, 2.062280) -- (1.665394, 2.062280) -- cycle;

\fill[color = black] (-8.238502, -3.510001) circle (5pt);

\fill[color = black] (7.221798, 2.062280) circle (5pt);

\draw[color = black, opacity = 0.4] (-14.480185, -5.300115) -- (-8.360139, -5.300115);
\draw[color = black, opacity = 0.4] (-10.940502, -3.510001) -- (-5.536501, -3.510001);
\draw[color = black, opacity = 0.4] (-8.142347, -2.094898) -- (-3.304388, -2.094898);
\draw[color = black, opacity = 0.4] (-5.874798, -0.948136) -- (-1.495543, -0.948136);
\draw[color = black, opacity = 0.4] (-4.000000, 0.000000) -- (0.000000, 0.000000);
\draw[color = black, opacity = 0.4] (-14.480185, -5.300115) -- (-4.000000, 0.000000);
\draw[color = black, opacity = 0.4] (-14.480185, -5.300115) -- (-14.480185, -5.300115);
\draw[color = black, opacity = 0.4] (-12.950173, -5.300115) -- (-10.940502, -3.510001);
\draw[color = black, opacity = 0.4] (-11.420162, -5.300115) -- (-8.142347, -2.094898);
\draw[color = black, opacity = 0.4] (-9.890150, -5.300115) -- (-5.874798, -0.948136);
\draw[color = black, opacity = 0.4] (-8.360139, -5.300115) -- (-4.000000, 0.000000);
\draw[color = black, opacity = 0.4] (-5.536501, -3.510001) -- (-3.000000, 0.000000);
\draw[color = black, opacity = 0.4] (-3.304388, -2.094898) -- (-2.000000, 0.000000);
\draw[color = black, opacity = 0.4] (-1.495543, -0.948136) -- (-1.000000, 0.000000);
\draw[color = black, opacity = 0.4] (-14.480185, -5.300115) -- (-8.360139, -5.300115);
\draw[color = black, opacity = 0.4] (-14.480185, -5.300115) -- (-4.000000, 0.000000);
\draw[color = black, opacity = 0.4] (-12.950173, -5.300115) -- (-3.000000, 0.000000);
\draw[color = black, opacity = 0.4] (-11.420162, -5.300115) -- (-2.000000, 0.000000);
\draw[color = black, opacity = 0.4] (-9.890150, -5.300115) -- (-1.000000, 0.000000);
\draw[color = black, opacity = 0.4] (-14.480185, -5.300115) -- (-14.480185, -5.300115);

\fill[color = blue, opacity = 0.8] (-4.000000, 0.000000) -- (0.000000, 0.000000) -- (1.500000, 2.598076) -- (-2.500000, 2.598076) -- cycle;

\draw[color = black, opacity = 0.4] (3.000000, 0.000000) -- (8.000000, 0.000000);
\draw[color = black, opacity = 0.4] (4.018052, 0.797003) -- (8.619550, 0.797003);
\draw[color = black, opacity = 0.4] (4.885805, 1.476341) -- (9.147634, 1.476341);
\draw[color = black, opacity = 0.4] (5.634254, 2.062280) -- (9.603114, 2.062280);
\draw[color = black, opacity = 0.4] (6.286421, 2.572843) -- (10.000000, 2.572843);
\draw[color = black, opacity = 0.4] (6.859763, 3.021695) -- (10.348915, 3.021695);
\draw[color = black, opacity = 0.4] (3.000000, 0.000000) -- (6.859763, 3.021695);
\draw[color = black, opacity = 0.4] (3.000000, 0.000000) -- (3.000000, 0.000000);
\draw[color = black, opacity = 0.4] (4.000000, 0.000000) -- (4.018052, 0.797003);
\draw[color = black, opacity = 0.4] (5.000000, 0.000000) -- (4.885805, 1.476341);
\draw[color = black, opacity = 0.4] (6.000000, 0.000000) -- (5.634254, 2.062280);
\draw[color = black, opacity = 0.4] (7.000000, 0.000000) -- (6.286421, 2.572843);
\draw[color = black, opacity = 0.4] (8.000000, 0.000000) -- (6.859763, 3.021695);
\draw[color = black, opacity = 0.4] (8.619550, 0.797003) -- (7.557593, 3.021695);
\draw[color = black, opacity = 0.4] (9.147634, 1.476341) -- (8.255424, 3.021695);
\draw[color = black, opacity = 0.4] (9.603114, 2.062280) -- (8.953254, 3.021695);
\draw[color = black, opacity = 0.4] (10.000000, 2.572843) -- (9.651085, 3.021695);
\draw[color = black, opacity = 0.4] (3.000000, 0.000000) -- (8.000000, 0.000000);
\draw[color = black, opacity = 0.4] (3.000000, 0.000000) -- (6.859763, 3.021695);
\draw[color = black, opacity = 0.4] (4.000000, 0.000000) -- (7.557593, 3.021695);
\draw[color = black, opacity = 0.4] (5.000000, 0.000000) -- (8.255424, 3.021695);
\draw[color = black, opacity = 0.4] (6.000000, 0.000000) -- (8.953254, 3.021695);
\draw[color = black, opacity = 0.4] (7.000000, 0.000000) -- (9.651085, 3.021695);
\draw[color = black, opacity = 0.4] (3.000000, 0.000000) -- (3.000000, 0.000000);

\draw[color = black, opacity = 0.4] (-4.000000, 0.000000) -- (0.000000, 0.000000);
\draw[color = black, opacity = 0.4] (-3.500000, 0.866025) -- (0.500000, 0.866025);
\draw[color = black, opacity = 0.4] (-3.000000, 1.732051) -- (1.000000, 1.732051);
\draw[color = black, opacity = 0.4] (-2.500000, 2.598076) -- (1.500000, 2.598076);
\draw[color = black, opacity = 0.4] (-4.000000, 0.000000) -- (-2.500000, 2.598076);
\draw[color = black, opacity = 0.4] (-4.000000, 0.000000) -- (-4.000000, 0.000000);
\draw[color = black, opacity = 0.4] (-3.000000, 0.000000) -- (-3.500000, 0.866025);
\draw[color = black, opacity = 0.4] (-2.000000, 0.000000) -- (-3.000000, 1.732051);
\draw[color = black, opacity = 0.4] (-1.000000, 0.000000) -- (-2.500000, 2.598076);
\draw[color = black, opacity = 0.4] (0.000000, 0.000000) -- (-1.500000, 2.598076);
\draw[color = black, opacity = 0.4] (0.500000, 0.866025) -- (-0.500000, 2.598076);
\draw[color = black, opacity = 0.4] (1.000000, 1.732051) -- (0.500000, 2.598076);
\draw[color = black, opacity = 0.4] (-4.000000, 0.000000) -- (0.000000, 0.000000);
\draw[color = black, opacity = 0.4] (-4.000000, 0.000000) -- (-2.500000, 2.598076);
\draw[color = black, opacity = 0.4] (-3.000000, 0.000000) -- (-1.500000, 2.598076);
\draw[color = black, opacity = 0.4] (-2.000000, 0.000000) -- (-0.500000, 2.598076);
\draw[color = black, opacity = 0.4] (-1.000000, 0.000000) -- (0.500000, 2.598076);
\draw[color = black, opacity = 0.4] (-4.000000, 0.000000) -- (-4.000000, 0.000000);

\fill[color = red, opacity = 0.9] (1.500000, 2.598076) -- (6.500000, 2.598076) -- (8.000000, 0.000000) -- (3.000000, 0.000000) -- cycle;
\draw[color = black, opacity = 0.4] (1.500000, 2.598076) -- (6.500000, 2.598076);
\draw[color = black, opacity = 0.4] (1.500000, 2.598076) -- (1.500000, 2.598076);
\draw[color = black, opacity = 0.4] (3.000000, 0.000000) -- (1.500000, 2.598076);
\draw[color = black, opacity = 0.4] (4.000000, 0.000000) -- (2.500000, 2.598076);
\draw[color = black, opacity = 0.4] (5.000000, 0.000000) -- (3.500000, 2.598076);
\draw[color = black, opacity = 0.4] (6.000000, 0.000000) -- (4.500000, 2.598076);
\draw[color = black, opacity = 0.4] (7.000000, 0.000000) -- (5.500000, 2.598076);
\draw[color = black, opacity = 0.4] (1.500000, 2.598076) -- (6.500000, 2.598076);
\draw[color = black, opacity = 0.4] (2.000000, 1.732051) -- (7.000000, 1.732051);
\draw[color = black, opacity = 0.4] (2.500000, 0.866025) -- (7.500000, 0.866025);
\draw[color = black, opacity = 0.4] (3.000000, 0.000000) -- (8.000000, 0.000000);
\draw[color = black, opacity = 0.4] (1.500000, 2.598076) -- (1.500000, 2.598076);
\draw[color = black, opacity = 0.4] (2.000000, 1.732051) -- (2.500000, 2.598076);
\draw[color = black, opacity = 0.4] (2.500000, 0.866025) -- (3.500000, 2.598076);
\draw[color = black, opacity = 0.4] (3.000000, 0.000000) -- (4.500000, 2.598076);
\draw[color = black, opacity = 0.4] (4.000000, 0.000000) -- (5.500000, 2.598076);
\draw[color = black, opacity = 0.4] (5.000000, 0.000000) -- (6.500000, 2.598076);
\draw[color = black, opacity = 0.4] (6.000000, 0.000000) -- (7.000000, 1.732051);
\draw[color = black, opacity = 0.4] (7.000000, 0.000000) -- (7.500000, 0.866025);
\draw[color = black, opacity = 0.4] (8.000000, 0.000000) -- (8.000000, 0.000000);
\draw[color = black, opacity = 0.4] (3.000000, 0.000000) -- (1.500000, 2.598076);
\fill[color = green, opacity = 0.8] (1.500000, 2.598076) -- (0.000000, 0.000000) -- (3.000000, 0.000000) -- cycle;

\draw[color = black, opacity = 0.4] (1.500000, 2.598076) -- (1.500000, 2.598076);
\draw[color = black, opacity = 0.4] (0.000000, 0.000000) -- (1.500000, 2.598076);
\draw[color = black, opacity = 0.4] (3.000000, 0.000000) -- (1.500000, 2.598076);
\draw[color = black, opacity = 0.4] (1.500000, 2.598076) -- (1.500000, 2.598076);
\draw[color = black, opacity = 0.4] (1.000000, 1.732051) -- (2.000000, 1.732051);
\draw[color = black, opacity = 0.4] (0.500000, 0.866025) -- (2.500000, 0.866025);
\draw[color = black, opacity = 0.4] (0.000000, 0.000000) -- (3.000000, 0.000000);
\draw[color = black, opacity = 0.4] (0.000000, 0.000000) -- (1.500000, 2.598076);
\draw[color = black, opacity = 0.4] (1.000000, 0.000000) -- (2.000000, 1.732051);
\draw[color = black, opacity = 0.4] (2.000000, 0.000000) -- (2.500000, 0.866025);
\draw[color = black, opacity = 0.4] (3.000000, 0.000000) -- (1.500000, 2.598076);
\draw[color = black, opacity = 0.4] (2.000000, 0.000000) -- (1.000000, 1.732051);
\draw[color = black, opacity = 0.4] (1.000000, 0.000000) -- (0.500000, 0.866025);

\draw[color = black, opacity = 0.4] (6.500000, 2.598076) -- (8.000000, 0.000000);

\fill[color = black] (1.500000, 2.598076) circle (5pt);

\node[yshift = 3mm] at (1.500000, 2.598076) {$p$};

\node[yshift = -3mm] at (-8.238502, -3.510001) {$x$};
\node[yshift = 3mm] at (7.221798, 2.062280) {$y$};
\end{tikzpicture}

%% file: sections/geometric_bound.tex
\subsection{The number of equilateral triangles corresponding to a given combinatorial line segment}
\begin{proposition} \label{prop_how_many_edges}
The number of chambers in $\mathcal{B}$ containing a given edge is exactly $q+1$.
\end{proposition}
\begin{proof}
This follows from counting the number of ways to extend a line in $\mathbb{F}_q^3$ to a plane.
\end{proof}

Given a combinatorial line segment (which is composed of finitely many edges), we call any edge containing one of the endpoints of the line segment a {\it boundary edge}.

Given an equilateral triangle $\mathcal{T}$ in $\mathcal{B}$ and a choice of one of its bounding combinatorial line segments, we may partition $\mathcal{T}$ into {\it levels} similarly to the way in which we partitioned sectors into levels. Each level is a trapezoid of width one as in Figure \ref{fig_convex_hull}.
\begin{lemma} \label{lemma_count_triangles}
Suppose $\ell$ is a combinatorial line segment in $\mathcal{B}$ of length $k$. The number of equilateral triangles which have $\ell$ as one of their sides is equal to $(q+1)q^{k-1}$.
\end{lemma}
\begin{proof}
Let $e_1$ be one of the boundary edges of $\ell$. By Proposition \ref{prop_how_many_edges}, there are $q+1$ chambers containing $e_1$. Let $\mathfrak{c}_1$ be one such chamber. Then by Lemma \ref{lemma_convex_hull}, the convex hull of $\ell$ and $\mathfrak{c}_1$ is a trapezoid as in Figure \ref{fig_convex_hull}; this shape exactly corresponds to the first level of an equilateral triangle containing $e_1$. We now wish to add a second level. Let $e_2$ be one of the boundary edges of the combinatorial line bounding this trapezoid which is parallel and opposite to $\ell$. Again, $e_2$ is contained in $q+1$ chambers, one of which already occurs in the trapezoid we already constructed. For each of the $q$ other choices, we may subsequently add another level to our trapezoid. We may continue this process until we end up with an equilateral triangle; at each step after the first there are $q$ choices for how to add the next level.
\end{proof}

\subsection{Coordinatizing the relative position of triples} \label{sec_coordinatize}
Suppose $(x, y; z)$ is a triple of vertices in $\mathcal{B}$. Suppose $d_{A^+}(x, y) = (r, s)$. Let $p$ be any confluence point of $(x, y; z)$. Let $\ell$ be the branch line of $(x, y; p)$. Suppose $\ell$ has length $k$ and is in the $\alpha$-direction with respect to $x$ (i.e. $\alpha \in \{(1, 0), (0, 1)\}$). Let $w$ be the point on $\ell$ closest to $x$ (with respect to $d_c(\cdot, \cdot)$, or equivalently with respect to $d(\cdot, \cdot)$). Suppose $d_{A^+}(x, w) = (a_1, a_2)$ and $d_{A^+}(p, z) = (b_1, b_2)$. Then we assign to $(x, y; p; z)$ the ``coordinates'' $(r, s; a_1, a_2, k, \alpha; b_1, b_2)$. Note, that by Remark \ref{remark_many_confluence_points}, a given point may get assigned multiple ``coordinates'' because there may be multiple choices of confluence points. Our next goal is to bound, for a given $x, y$ with $d_{A^+}(x, y) = (r, s)$ and a give choice of $(r, s;a_1, a_2, k, \alpha; b_1, b_2)$, the number of $z$ such that there exists a confluence point $p$ for $(x, y; z)$ such that $(x, y; p; z)$ gets assigned the given coordinates.

\begin{lemma} \label{lemma_count_triples}
Suppose $x, y$ are fixed with $d_{A^+}(x, y) = (r, s)$. Then
\begin{gather}
\left\{ \begin{array}{c}
    \textnormal{$\# z$ s.t. $\exists$ confluence point $p$ of $(x, y; z)$ s.t. } \\
    \textnormal{coordinates of $(x, y; p; z)$ are $(r, s; a_1, a_2, k, \alpha; b_1, b_2)$}
  \end{array} \right\} \leq \frac{2}{\nu_3(q^{-1})} (q^2)^{\frac{k}{2} + b_1 + b_2}. \label{eqn_count_triples}
\end{gather}
\end{lemma}

\begin{proof}
The choice of coordinates uniquely determines the branch line $\ell$. Let $\mathfrak{T}$ denote the collection of all equilateral triangles which have one of their bounding line segments equal to $\ell$; the cardinality of this set is given by Lemma \ref{lemma_count_triangles}. Given $\mathcal{T} \in \mathfrak{T}$, let $a(\mathcal{T})$ denote the vertex of $\mathcal{T}$ opposite to $\ell$. Let $X$ denote the set:
\begin{gather}
X := \{(\mathcal{T}, w) : \mathcal{T} \in \mathfrak{T} \textnormal{ and } w \in G/K \textnormal{ and } d_{A^+}(a(\mathcal{T}), w) = (b_1, b_2) \}. \nonumber
\end{gather}
Given $(x, y; p; z)$ on the left hand side of \eqref{eqn_count_triples}, we can associate to it an element of $X$, namely the equilateral triangle $\mathcal{T}'$ associated to the primitive triple $(x, y; p)$ as in Lemma \ref{lemma_actual_classification} (clearly then $a(\mathcal{T}') = p$) and the point $z$ itself. This map is clearly an injection. On the other hand, for a fixed $\mathcal{T} \in \mathfrak{T}$, the number of points in $X$ whose first entry is $\mathcal{T}$ is exactly $N_{(b_1, b_2)}$ which is given by \eqref{eqn_n_lambda} (note that that formula is given in partition coordinates, but here we are using cone coordinates). Therefore the cardinality of $X$ is exactly given by the right hand side of \eqref{eqn_count_triples}.
\end{proof}

Suppose $(x, y; p; z)$ has coordinates $(r, s; a_1, a_2, k, \alpha; b_1, b_2)$. Then if $\alpha = (1, 0)$, we have $d_{A^+}(x, p) = (a_1, a_2 + k)$ and $d_{A^+}(y, p) = (s - a_2 + k, r - a_1 - k)$. If $\alpha = (0, 1)$, we have
$d_{A^+}(x, p) = (a_1 + k, a_2)$ and  
$d_{A^+}(y, p) = (s - a_2 - k, r - a_1 + k)$. This follows from Lemma \ref{lemma_actual_classification}.

\subsection{The polytope parametrizing allowable coordinates of triples of points} \label{sec_polytope_parametrization}

\subsubsection{The defining inequalities of the polytope}
\begin{proposition} \label{prop_triples_polytope}
Suppose $(x, y; p; z)$ has coordinates $(r, s; a_1, a_2, k, \alpha; b_1, b_2)$ with $z \in x E_m \cap y E_m$. Then the following inequalities must be satisfied:
\begin{align}
a_1 &\leq r \nonumber \\
a_2 &\leq s \nonumber \\
r, s, a_1, a_2, k, b_1, b_2, m &\geq 0 \label{eqn_base_inequalities}
\end{align}

\begin{enumerate}
\item Suppose $\alpha = (1, 0)$. Then the following inequalities are also satisfied:
\begin{align}
(a_1+b_1) + 2 (a_2 + k + b_2) &\leq 2 m \ \ \ \ \ \ \ \ \ \ \ \ d_{A^+}(x, z) \in P_m \nonumber \\
(s - a_2 + k + b_1) + 2 (r - a_1 - k + b_2) &\leq 2 m \ \ \ \ \ \ \ \ \ \ \ \ d_{A^+}(y, z) \in P_m \nonumber \\
a_1 + k &\leq r \ \ \ \ \ \ \ \ \ \ \ \ \ \ \ell \subset \para(x, y). \label{eqn_ineqs_alpha_10}
\end{align}
\item Suppose $\alpha = (0, 1)$. Then the following inequalities are also satisfied:
\begin{align}
(a_1 + k + b_1) + 2 (a_2 + b_2) &\leq 2 m \ \ \ \ \ \ \ \ \ \ \ \ d_{A^+}(x, z) \in P_m \nonumber\\
(s - a_2 - k + b_1) + 2 (r - a_1 + k + b_2) &\leq 2m \ \ \ \ \ \ \ \ \ \ \ \ d_{A^+}(y, z) \in P_m \nonumber \\
a_2 + k &\leq s \ \ \ \ \ \ \ \ \ \ \ \ \ \ \ell \subset \para(x, y). \label{eqn_ineqs_alpha_01}
\end{align}
\end{enumerate}
\end{proposition}

\subsubsection{The dominating term of the sum over lattice points in the polytope}
\begin{lemma} \label{lemma_maximal_vertex}
All $z \in x E_m \cap y E_m$ with $d_{A^+}(x, y) = (r, s)$ have their coordinates $(x, y; p; z)$ (with any choice of confluence point $p$) satisfying:
\begin{gather}
\frac{k}{2} + b_1 + b_2 \leq 2m - \frac{r}{2} - \frac{s}{2}. \label{eqn_dominating_term}
\end{gather}
Furthermore, when $m, r, s$ are fixed, there is at most one possible set of coordinates for such $(x, y; p; z)$ such that \eqref{eqn_dominating_term} becomes an equality.
\end{lemma}
\begin{proof}
Suppose $(x, y; z; p)$ is as in the statement of the lemma. Depending on $\alpha$, the coordinates of $z$ must satisfy all of \eqref{eqn_base_inequalities} plus either all of \eqref{eqn_ineqs_alpha_10} or all of \eqref{eqn_ineqs_alpha_01}. In all cases we refer to these inequalities collectively as the coordinates inequalities. Suppose the coordinates of $z$ also satisfy the inequality:
\begin{gather}
\frac{k}{2} + b_1 + b_2 \geq 2m - \frac{r}{2} - \frac{s}{2}. \label{eqn_opposite_dominating}
\end{gather}
The set of points in $\mathbb{R}^8$ satisfying the coordinates inequalities as well as \eqref{eqn_opposite_dominating} forms a convex polyhedron. Using the Polyhedron package in Sage \cite{sage}, one computes that, regardless of what $\alpha$ is, this polytope is the conical span of the following vectors in coordinates $(m, r, s, a_1, a_2, k, b_1, b_2)$: 
$u_1 = (1, 2, 2, 2, 0, 0, 0, 0)$ and
$u_2 = (1, 0, 0, 0, 0, 0, 2, 0)$.
Any vector in this conical span results in \eqref{eqn_opposite_dominating} becoming an equality. Hence no point in this polytope ever violates \eqref{eqn_dominating_term}. Furthermore, when $m, r, s$ are fixed, there is at most one point in the conical span whose first three coordinates are $(m, r, s)$ because the first three coordinates of $u_1$, namely $(1, 2, 2)$, and the first three coordinates of $u_2$, namely $(1, 0, 0)$, are linearly independent in $\mathbb{R}^3$. In fact for such a point we must have $r = s = a_1$, $a_2, k, b_2 = 0$ and $b_1 = 2m - r$.
\end{proof}

\subsection{Proof of Proposition \ref{prop_vol_e_lambda_m}} \label{sec_pf_e_m_lambda_bound}

\begin{proof}[Proof of Proposition \ref{prop_vol_e_lambda_m} (Upper Bound on $\card(E_m^\lambda)$)]
Recall that $\tilde{E}_m^\lambda$ is the pullback to $G$ of $E_m^\lambda = E_m \cap \varpi^{\lambda} E_m \subset G/K$ with $\lambda \in A^+$. Clearly $\vol(\tilde{E}_m^\lambda) = \card(E_m^\lambda)$. Suppose $\lambda = (r, s)$ in cone coordinates. We can associate to each $z \in E_m^\lambda$ the coordinates $(m; r, s; a_1, a_2, k, \alpha; b_1, b_2)$ (after choosing some confluence point $p$). Let $(\beta, \cdot)$ be the functional defined by
\begin{gather}
    (\beta, (a_1, a_2, k, b_1, b_2)) = \frac{k}{2} + b_1 + b_2. \label{eqn_beta}
\end{gather}
By Lemma \ref{lemma_count_triples}, the number of $z$ which map to a given set of coordinates is at most 
\begin{gather}
    (q^2)^{\frac{k}{2} + b_1 + b_2} = (q^2)^{(\beta, (a_1, a_2, k, b_1, b_2))}. \label{eqn_to_sum}
\end{gather}
Let $P_\alpha(m, r, s)$ be the polytope defined by the inequalities in \eqref{eqn_base_inequalities} and either \eqref{eqn_ineqs_alpha_10} (if $\alpha = (1, 0)$) or \eqref{eqn_ineqs_alpha_01} (if $\alpha = (0, 1)$). We can then bound $\card(E_m^\lambda)$ by summing \eqref{eqn_to_sum} over all integer points (in coordinates $(a_1, a_2, k, b_1, b_2)$) in the polytopes $\mathcal{P}_{(1, 0)}(m, r, s)$ and $\mathcal{P}_{(0, 1)}(m, r, s)$. 

For the remainder of the proof, we assume that we have fixed an $\alpha$ (the subsequent discussion does not depend on $\alpha$). Notice that all of the defining inequalities of $P_{\alpha}(m, r, s)$ are of the form
$f_i (a_1, a_2, k, b_1, b_2) \leq g_i (m, r, s)$
where $f_i$ and $g_i$ are linear functionals. We have that $(m, r, s)$ parametrizes some vector space which we call $V_{(m, r, s)}$.

We also have each of our polytopes $\mathcal{P}_{\alpha}(m, r, s)$ living inside of some five-dimensional space which we call $W$ and which is naturally coordinatized via $(a_1, a_2, k, b_1, b_2)$. We can partition up $V_{(m, r, s)}$ according to the type of the underlying polytope in $W$ together with the information about which vertices maximize $(\beta, \cdot)$. We call each such component a {\it $\beta$-region}. Clearly there are only finitely many such regions.

Because all polytopes associated to a given $\beta$-region have the same type, it makes sense to discuss the ``same'' vertex for different polytopes in a given $\beta$-region. Some $\beta$-regions are such that their associated polytopes have a vertex $v^*$ which, at least for some polytopes in the $\beta$-region, satisfies $(\beta, v^*) = 2m - r/2 - s/2$; we shall refer to such $\beta$-regions as {\it extremal}. By Lemma \ref{lemma_maximal_vertex}, in such cases $v^*$ is the unique vertex satisfying this equation and is the unique vertex maximizing $(\beta, \cdot)$. In general, at the vertex $w$ maximizing $(\beta, \cdot)$ in any polytope, we have that $\beta$ is in the polar cone of the cone generated at $w$. Hence $\beta$ is in the polar cone of $\textnormal{Cone}_{\mathcal{P}_{\alpha}(m, r, s)}(v^*)$ (which we view as a cone based at the origin, not based at $v^*$). Furthermore, because $v^*$ is the unique vertex maximizing $(\beta, \cdot)$, we must have $\beta$ in the interior of the polar cone.

Recall that the structure of $\textnormal{Cone}_{\mathcal{P}_{\alpha}(m, r, s)}(v^*)$ only depends on the type of the underlying polytope (see Section \ref{sec_polytope_type}). Hence this cone does not depend on $(m, r, s)$ (as long as we stay in the same $\beta$-region). Because $\beta$ is in the interior of the polar cone, we have that the sum of \eqref{eqn_to_sum} over the lattice points in $\textnormal{Cone}_{\mathcal{P}_{\alpha}(m, r, s)}(v^*)$ has a finite value; in fact this is exactly what is referred to as
\begin{gather}
    \varsigma(\textnormal{Cone}_{\mathcal{P}_{\alpha}(m, r, s)}(v^*); 2\beta) \label{eqn_sum_over_this_cone}
\end{gather}
in Section \ref{brion_section}. If we translate this cone to be based at $v^*$, then the sum of \eqref{eqn_to_sum} over the lattice points in this cone gives an upper bound on the sum of \eqref{eqn_to_sum} over the lattice points in $\mathcal{P}_{\alpha}(m, r, s)$. On the other hand, this quantity is at most $(q^2)^{2m - r/s - s/2}$ times \eqref{eqn_sum_over_this_cone}. Hence we obtain that for polytopes in extremal $\beta$-regions:
\begin{gather}
    \card(E_m^{\lambda}) \leq \varsigma(\textnormal{Cone}_{\mathcal{P}_{\alpha}(m, r, s)}(v^*); 2 \beta) (q^2)^{2m - r/2 - s/2}. \nonumber
\end{gather}

We now consider $\beta$-regions which are adjacent to extremal ones, namely ones whose closure intersects the closure of an extremal $\beta$-region. If we have a sequence of points $(m_i, r_i, s_i)$ in such a $\beta$-region which approaches the boundary of an extremal $\beta$-region, any vertex $w_i$ in the associated polytopes which maximizes $(\beta, \cdot)$ must converge to $v^*$ and hence the collection of coordinate inequalities which become equalities at $w_i$ must be a subset of the inequalities that become equalities at $v^*$. 

From the analysis in the proof of Lemma \ref{lemma_maximal_vertex}, we must have $v^* = (r, 0, 0, 2m - r, 0)$ with $r = s$. Therefore the collection of inequalities which become equalities at $v^*$ is (supposing, for the moment, that $\alpha = (1, 0)$):

\begin{minipage}{0.4\textwidth}
        \begin{align*}
            a_1 & \leq r \\
            a_2 & \geq 0 \\
            k & \geq 0 \\
            b_2 & \geq 0
        \end{align*}
\end{minipage}
\begin{minipage}{0.5\textwidth}
        \begin{align*}
            a_1 + 2 a_2 + 2k + b_1 + 2 b_2 &\leq 2m \\
            -2a_1 - a_2 - k + b_1 + 2b_2 &\leq 2m - s - 2r \\
            a_1 + k &\leq r
        \end{align*}
\end{minipage}
\begin{minipage}{0.1\textwidth}
\end{minipage}

\vspace{2mm}
The cone of any such $w_i$ in a region adjacent to an extremal $\beta$-region must be contained in the cone obtained from at least 5 of these inequalities (as the relevant polytope lives in a five-dimensional space). Furthermore this cone must contain $\beta$ in its polar cone. Using the Sage Polyhedron package \cite{sage}, we can enumerate all such possible cones; we in fact find that there are 4 possibilities and in all cases we have that $\beta$ is in the interior of the polar cone. We may then use the same style of analysis as for the extremal $\beta$-regions to conclude that there exists a $C_1$ such that the following holds (we can take $C_1$ to be the maximum over all possible $\varsigma (\textnormal{Cone}_{\mathcal{P}_{\alpha}(m, r, s)}(w_i); 2 \beta)$; by the preceding analysis this set is finite and hence the maximum is finite):
\begin{gather}
    \card(E_m^\lambda) \leq C_1 (q^2)^{(\beta, w_i)} \leq C_1 (q^2)^{2m - r/2 - s/2}, \nonumber
\end{gather}
for all $(m, r, s)$ in a $\beta$-region adjacent to an extremal one. A similar analysis can be performed for $\alpha = (0, 1)$. 

We are left now with analyzing the $\beta$-regions which are not extremal nor adjacent to an extremal one; we call such $\beta$-regions {\it tame}. Let $B$ be a tame $\beta$-region. We seek to apply the degenerate case of Brion's formula for polytopes associated to points in $B$. However, we first note that Brion's formula requires that all vertices of the underlying polytope lie in some lattice. We are in particular interested in $\mathbb{Z}^5 \subset W$, but the vertices of the polytopes might not always be integer points. However, because all of the defining coordinate inequalities have integer coefficients, and because we are only really interested in the case when $(m, r, s) \in \mathbb{Z}^3$, there exists $L \in \mathbb{N}$ such that the vertices of the underlying polytopes always lie in $(\mathbb{Z}/L)^5$ (assuming $m, r, s$ are integers). If we sum up \eqref{eqn_to_sum} over the lattice points in this bigger lattice lying inside of the underlying polytope, then we still obtain an upper bound for $\card(E_m^\lambda)$.

We are now in a position to use degenerate Brion's formula, i.e. Theorem \ref{prop_degenerate_brion}. This tells us that for all polytopes associated to points in $B$, we have
\begin{gather}
    \sum_{x \in \mathcal{P}_{\alpha}(m, r, s) \cap (\mathbb{Z}/L)^5} = \sum_{\textnormal{vertices $v$ of $\mathcal{P}_{\alpha}(m, r, s)$}} R_v(m, r, s) (q^2)^{h_v(m, r, s)}, \nonumber
\end{gather}
where $R_v(m, r, s)$ is a polynomial in $m, r, s$ and $h_v(m, r, s)$ is some linear functional in $m, r, s$. More specifically $h_v(m, r, s)$ is obtained by dotting $\beta$ with a given vertex $v$ in the family of polytopes associated to points in $B$; the coordinates of that vertex are in turn linear functionals in $m, r, s$ so long as we are in a given $\beta$-region. For any fixed degree $n$ polynomial in $m, r, s$ may be bounded in absolute value on the locus $m, r, s \geq 0$ by $C (m + r + s)^n$ for some choice of $C$.

Consider the following locus in $V_{(m, r, s)}$:
\begin{gather}
    A = \{(m, r, s): m + r + s = 1 \textnormal{ and } m, r, s \geq 0\}. \nonumber
\end{gather}
Consider the set $A \cap B$. Elements in this set are uniformly far away (in, say, the Euclidean metric) from the intersection of the extremal $\beta$-regions with $A$ because all adjacent $\beta$-regions to the extremal ones are also in the complement of $B$. On $A \cap B$ consider the function $h_v(m, r, s) - (2m - r/2 - s/2)$. By Lemma \ref{lemma_maximal_vertex}, this function is negative on $A \cap B$, but since $A$ is compact we in fact must have that there exists a $C_2$ such that 
\begin{gather}
    h_v(m, r, s) - (2m - r/2 - s/2) \leq C_2 < 0 \nonumber
\end{gather}
on $A \cap B$.

Suppose $m + r + s \geq 1$ (since we assume $(m, r, s) \in \mathbb{Z}^3$, this only excludes the case of $m = r = s = 0$; in that case $\card(E_m^\lambda) = 1$). We therefore get that there exists a $C_3, C_4 > 0$ and a $p \geq 0$ such that:
\begin{align*}
R_v(m, r, s) (q^2)^{h_v(m, r, s)} &\leq C_3 (m + r + s)^p \frac{(q^2)^{h_v(m, r, s)}}{(q^2)^{2m - r/2 - s/2}} (q^2)^{2m - r/2 - s/2} \\
    & \leq C_3 (m + r + s)^p (q^2)^{C_2 (m + r + s)} (q^2)^{2m - r/2 - s/2} \\
    & \leq C_4 (q^2)^{2m - r/2 - s/2}.
\end{align*}

By combining this bound over all vertices in $v$, we get that there exists a $C_5 > 0$ such that for all $(m, r, s)$ in $B$:
\begin{gather}
    \card(E_m^\lambda) \leq C_5 (q^2)^{2m - r/2 - s/2}. \nonumber
\end{gather}

Finally, by combining the bounds we obtained for the extremal, adjacent to extremal, and tame $\beta$-regions (which altogether consist of only finitely many regions), we get that there exists a $C_6$ such that for all $(m, r, s) \in \mathbb{N}^3$ we have:
\begin{gather}
    \card(E_m^\lambda) \leq C_6 (q^2)^{2m - r/2 - s/2} = C_6 (q^2)^{(\delta, m \cdot p^\dag - \frac{\lambda}{2})}. \nonumber
\end{gather}
\end{proof}

%% file: sections/final_brion.tex
\subsection{Proof of Proposition \ref{prop_final_brion}}
\begin{proof}[Proof of Proposition \ref{prop_final_brion} (Bounding the Sum over $H_M^\Lambda$)]
We now seek to prove:
\begin{gather}
\sum_{\lambda \in H_M^\Lambda} (q^2)^{\theta ( \delta, \lambda - 2|\lambda|_H \cdot p^\dag )} \lesssim M. \label{eqn_final_brion}
\end{gather}
We wish to use Brion's formula. However, the exponent in each summand is not quite a linear functional.

We now examine $|\lambda|_H$ more closely. Recall that $H$ is defined by 
$2r + s \leq 6$ and 
$r + 2s \leq 6$. Therefore
$|\lambda|_H = \max \Big\{ \frac{2r + s}{6}, \frac{r + 2s}{6} \Big\}$,
where $\lambda = (r, s)$ in cone coordinates. The locus where $2r +s = r + 2s$ is exactly the locus $r = s$. We may partition $H$ into two pieces: $H^{r \leq s}$ and $H^{r \geq s}$ by adding in the relevant constraint to the definition of $H$. See Figure \ref{fig_split_h}. On $H^{r \leq s}$, we have that $|\lambda|_H = \frac{r}{6} + \frac{s}{3}$, and on $H^{r \geq s}$, we have $|\lambda|_H = \frac{r}{3} + \frac{s}{6}$. Hence we define:
$\alpha_{r \leq s} := \frac{r}{6} + \frac{s}{3}$ and
$\alpha_{r \geq s} := \frac{r}{3} + \frac{s}{6}$.
We may thus split up the left hand side of \eqref{eqn_final_brion} into a sum over $(H^{r \leq s})_M^\Lambda$ and $(H^{r \geq s})_M^\Lambda$, and on each piece the summands become exactly of the form $(q^2)^{f(r, s)}$ where $f$ is a linear functional in $r, s$.

\begin{figure}[h] 
\centering
\subfile{../images/h_norm}
\caption{}\label{fig_split_h}
\end{figure}

Let's focus now on $H^{r \leq s}$. On this polytope we have:
\begin{align*}
\theta ( \delta, \lambda - 2|\lambda|_H \cdot p^\dag ) = \theta \Big( \delta, \lambda - 2 ( \alpha_{r \leq s}, \lambda ) \cdot p^\dag \Big)  = \theta \Big( \lambda, \delta - 2 ( \delta, p^\dag ) \alpha_{r \leq s} \Big)
\end{align*}
In cone coordinates, we have
\begin{gather*}
    \delta = (1, 1), \ \ p^\dag = (2, 0), \ \ \alpha_{r \leq s} = \Big(\frac{1}{3}, \frac{1}{6} \Big) \implies \delta - 2 ( \delta, p^\dag ) \alpha_{r \leq s} = \Big(\frac{1}{3}, -\frac{1}{3} \Big) =: \beta_{r \leq s}.
\end{gather*}
Similarly we have
\begin{gather*}
\alpha_{r \geq s} = \Big(\frac{1}{6}, \frac{1}{3} \Big) \implies \delta - 2 ( \delta, p^\dag ) \alpha_{r \geq s} = \Big(-\frac{1}{3}, \frac{1}{3} \Big) =: \beta_{r \geq s}.
\end{gather*}

The vertices of $H$ are $h_1 = (0, 0)$, $h_2 = (3, 0)$, $h_3 = (0, 3)$, and $h_4 = (2, 2) =: h^\dag$. We have that the vertices of $H^{r \leq s}$ are $h_1, h_3, h_4$, and the vertices of $H^{r \geq s}$ are $h_1, h_2, h_4$. 

On $H^{r \leq s}$ we have $( h_1, \beta_{r \leq s} ) = ( h_4, \beta_{r \leq s} ) = 0$ and $( h_3, \beta_{r \leq s} )  = - \frac{1}{9}$. Using degenerate Brion's formula, i.e. Theorem \ref{prop_degenerate_brion}, we get that there exists a constant $C_1$ and a degree one polynomial $f(M)$ such that
\begin{gather*}
\sum_{ \lambda \in (H^{r \leq s})_M^\Lambda} (q^2)^{\theta ( \delta, \lambda - 2 |\lambda|_H \cdot p^\dag )} = f(M) + C_1 (q^2)^{-\frac{M \cdot \theta}{9}} \leq C_2 \cdot M 
\end{gather*}
for some $C_2$ and for all $M \geq 1$. An analogous analysis may be carried out on $(H^{r \geq s})_M^\Lambda$. 
\end{proof}

%% file: images/h_norm.tex
\begin{tikzpicture}[scale=1.5]
\draw[color = orange, line width = 2]  (0.000000, 0.000000) --  (3.000000, 0.000000) --  (3.000000, 1.732051) --  (1.500000, 2.598076) -- cycle;
\fill[color = red, opacity = 0.3]  (0.000000, 0.000000) --  (3.000000, 1.732051) --  (3.000000, 0.000000) -- cycle;
\fill[color = blue, opacity = 0.3]  (0.000000, 0.000000) --  (3.000000, 1.732051) --  (1.500000, 2.598076) -- cycle;
\draw[color = black, opacity = 0.1] (0.000000, 0.000000) -- (4.000000, 0.000000);
\draw[color = black, opacity = 0.1] (0.500000, 0.866025) -- (4.000000, 0.866025);
\draw[color = black, opacity = 0.1] (1.000000, 1.732051) -- (4.000000, 1.732051);
\draw[color = black, opacity = 0.1] (1.500000, 2.598076) -- (4.000000, 2.598076);
\draw[color = black, opacity = 0.1] (2.000000, 3.464102) -- (4.000000, 3.464102);
\draw[color = black, opacity = 0.1] (0.000000, 0.000000) -- (2.309401, 4.000000);
\draw[color = black, opacity = 0.1] (0.000000, 0.000000) -- (0.000000, 0.000000);
\draw[color = black, opacity = 0.1] (1.000000, 0.000000) -- (0.500000, 0.866025);
\draw[color = black, opacity = 0.1] (2.000000, 0.000000) -- (1.000000, 1.732051);
\draw[color = black, opacity = 0.1] (3.000000, 0.000000) -- (1.500000, 2.598076);
\draw[color = black, opacity = 0.1] (4.000000, 0.000000) -- (2.000000, 3.464102);
\draw[color = black, opacity = 0.1] (4.000000, 1.732051) -- (2.690599, 4.000000);
\draw[color = black, opacity = 0.1] (4.000000, 3.464102) -- (3.690599, 4.000000);
\draw[color = black, opacity = 0.1] (1.000000, 0.000000) -- (3.309401, 4.000000);
\draw[color = black, opacity = 0.1] (2.000000, 0.000000) -- (4.000000, 3.464102);
\draw[color = black, opacity = 0.1] (3.000000, 0.000000) -- (4.000000, 1.732051);
\draw[color = black, opacity = 0.1] (4.000000, 0.000000) -- (4.000000, 0.000000);
\draw[color = black, opacity = 0.1] (0.000000, 0.000000) -- (0.000000, 0.000000);
\draw[color = black, opacity = 0.1] (-0.500000, 0.866025) -- (0.500000, 0.866025);
\draw[color = black, opacity = 0.1] (-1.000000, 1.732051) -- (1.000000, 1.732051);
\draw[color = black, opacity = 0.1] (-1.000000, 2.598076) -- (1.500000, 2.598076);
\draw[color = black, opacity = 0.1] (-1.000000, 3.464102) -- (2.000000, 3.464102);
\draw[color = black, opacity = 0.1] (0.000000, 0.000000) -- (2.309401, 4.000000);
\draw[color = black, opacity = 0.1] (-0.500000, 0.866025) -- (1.309401, 4.000000);
\draw[color = black, opacity = 0.1] (-1.000000, 1.732051) -- (0.309401, 4.000000);
\draw[color = black, opacity = 0.1] (-1.000000, 3.464102) -- (-0.690599, 4.000000);
\draw[color = black, opacity = 0.1] (0.000000, 0.000000) -- (-1.000000, 1.732051);
\draw[color = black, opacity = 0.1] (0.500000, 0.866025) -- (-1.000000, 3.464102);
\draw[color = black, opacity = 0.1] (1.000000, 1.732051) -- (-0.309401, 4.000000);
\draw[color = black, opacity = 0.1] (1.500000, 2.598076) -- (0.690599, 4.000000);
\draw[color = black, opacity = 0.1] (2.000000, 3.464102) -- (1.690599, 4.000000);
\draw[color = black, opacity = 0.1] (-1.000000, 0.000000) -- (0.000000, 0.000000);
\draw[color = black, opacity = 0.1] (-1.000000, 0.866025) -- (-0.500000, 0.866025);
\draw[color = black, opacity = 0.1] (-1.000000, 1.732051) -- (-1.000000, 1.732051);
\draw[color = black, opacity = 0.1] (0.000000, 0.000000) -- (0.000000, 0.000000);
\draw[color = black, opacity = 0.1] (-1.000000, 0.000000) -- (-0.500000, 0.866025);
\draw[color = black, opacity = 0.1] (-1.000000, 1.732051) -- (-1.000000, 1.732051);
\draw[color = black, opacity = 0.1] (0.000000, 0.000000) -- (-1.000000, 1.732051);
\draw[color = black, opacity = 0.1] (-1.000000, 0.000000) -- (-1.000000, 0.000000);
\draw[color = black, opacity = 0.1] (-1.000000, 0.000000) -- (0.000000, 0.000000);
\draw[color = black, opacity = 0.1] (-0.577350, -1.000000) -- (0.000000, 0.000000);
\draw[color = black, opacity = 0.1] (-1.000000, 0.000000) -- (-1.000000, 0.000000);
\draw[color = black, opacity = 0.1] (0.000000, 0.000000) -- (0.000000, 0.000000);
\draw[color = black, opacity = 0.1] (-1.000000, -0.866025) -- (-0.500000, -0.866025);
\draw[color = black, opacity = 0.1] (-0.500000, -0.866025) -- (-1.000000, 0.000000);
\draw[color = black, opacity = 0.1] (0.000000, 0.000000) -- (0.000000, 0.000000);
\draw[color = black, opacity = 0.1] (-0.577350, -1.000000) -- (0.000000, 0.000000);
\draw[color = black, opacity = 0.1] (0.577350, -1.000000) -- (0.000000, 0.000000);
\draw[color = black, opacity = 0.1] (-0.500000, -0.866025) -- (0.500000, -0.866025);
\draw[color = black, opacity = 0.1] (0.422650, -1.000000) -- (0.500000, -0.866025);
\draw[color = black, opacity = 0.1] (-0.422650, -1.000000) -- (-0.500000, -0.866025);
\draw[color = black, opacity = 0.1] (0.000000, 0.000000) -- (4.000000, 0.000000);
\draw[color = black, opacity = 0.1] (0.000000, 0.000000) -- (0.000000, 0.000000);
\draw[color = black, opacity = 0.1] (0.577350, -1.000000) -- (0.000000, 0.000000);
\draw[color = black, opacity = 0.1] (1.577350, -1.000000) -- (1.000000, 0.000000);
\draw[color = black, opacity = 0.1] (2.577350, -1.000000) -- (2.000000, 0.000000);
\draw[color = black, opacity = 0.1] (3.577350, -1.000000) -- (3.000000, 0.000000);
\draw[color = black, opacity = 0.1] (4.000000, 0.000000) -- (4.000000, 0.000000);
\draw[color = black, opacity = 0.1] (0.500000, -0.866025) -- (4.000000, -0.866025);
\draw[color = black, opacity = 0.1] (0.500000, -0.866025) -- (1.000000, 0.000000);
\draw[color = black, opacity = 0.1] (1.422650, -1.000000) -- (2.000000, 0.000000);
\draw[color = black, opacity = 0.1] (2.422650, -1.000000) -- (3.000000, 0.000000);
\draw[color = black, opacity = 0.1] (3.422650, -1.000000) -- (4.000000, 0.000000);

\node at (2, 0.5) {$H^{r \geq s}$};
\node at (1.5, 1.5) {$H^{r \leq s}$};

\fill[color = black] (0, 0) circle (2pt);
\fill[color = black] (3, 0) circle (2pt);
\fill[color = black] (3.000000, 1.732051) circle (2pt);
\fill[color = black] (1.500000, 2.598076) circle (2pt);

\node[yshift = -3mm] at (0, 0) {$h_1$};
\node[yshift = -3mm] at (3, 0) {$h_2$};
\node[yshift = 3mm] at (1.500000, 2.598076) {$h_3$};
\node[yshift = 3mm, xshift = 5mm] at (3.000000, 1.732051) {$h_4 = h^\dag$};
\end{tikzpicture}